\documentclass[reqno]{amsart}

\usepackage{amsmath}
\usepackage{hyperref}
\hypersetup
{
colorlinks,
citecolor=blue,
filecolor=red,
linkcolor=blue,
urlcolor=black
}
\usepackage{amssymb}
\usepackage{amsthm}
\usepackage{extarrows}
\usepackage{enumitem}
\usepackage{mathrsfs}
\usepackage[numbers,sort&compress]{natbib}
\usepackage{color}
\usepackage
[
a4paper,% other options: a3paper, a5paper, etc
vmargin=2.7cm,
hmargin=2.5cm 
]{geometry}

\def\R{\mathbb R}
\def\Z{\mathbb Z}
\def\T{\mathbb T}
\def\N{\mathbb N}
\def\E{\mathbb E}
\def\e{\varepsilon}
\def\p{\mathbb P}
\def\F{\mathcal F}
\def\U {\mathfrak U}
\def\W{\mathcal W}

\def\X{\mathcal X}
\def\EE{\mathcal E}
\def\ZZ{\mathcal Z}
\def\Y{\mathcal Y}
\def\LL{\mathcal L}

\def\dd{\, {\rm d}}
\def\Wlip{W^{1,\infty}}
\def\pas{{\mathbb P}\text{-}{\rm a.s.}}

\theoremstyle{plain}
\numberwithin{equation}{section}
\newtheorem{Theorem}{Theorem}[section]
\newtheorem{Definition}{Definition}[section]
\newtheorem{Proposition}{Proposition}[section]
\newtheorem{Lemma}{Lemma}[section]

\theoremstyle{definition}
\newtheorem{Remark}{Remark}[section]
\newtheorem{Hypothesis}{Hypothesis}

\begin{document}

\title[Stochastic Camassa-Holm Type Equation]{Well-posedness for a stochastic Camassa-Holm type equation with higher order nonlinearities}

\author[Y. Miao]{Yingting Miao}
\address{School of Mathematics, South China University of Technology, 
Guangzhou 510641, PR China}
\email{yingtmiao2-c@my.cityu.edu.hk}

\author[C. Rohde]{Christian Rohde}
\address{Institut f\"{u}r Angewandte Analysis und Numerische Simulation, Universit\"{a}t Stuttgart, Pfaffenwaldring 57, 70569 Stuttgart, Germany}
\email{christian.rohde@mathematik.uni-stuttgart.de}
\thanks{C. R. is funded  by Deutsche Forschungsgemeinschaft
(DFG, German Research Foundation) under
Germany's Excellence Strategy --EXC 2075--390740016}

\author[H. Tang]{Hao Tang}
\address{Department of Mathematics,
University of Oslo, P.O. Box 1053, Blindern, N-0316 Oslo, Norway}
\email{haot@math.uio.no; haotang4-c@my.cityu.edu.hk}
\thanks{The major part of this work was carried out when H. T. was supported by the Alexander von Humboldt Foundation.}

%    General info

%\date{January 8, 2020.}

%\date{\today}



\begin{abstract}
This paper aims at studying a generalized Camassa--Holm equation under random perturbation. We establish a local well-posedness result in the sense of Hadamard, i.e., existence, uniqueness and continuous dependence on initial data, as well as  blow-up criteria for pathwise solutions in the Sobolev spaces $H^s$ with $s>3/2$ for $x\in\R$.  The analysis on continuous dependence on initial data for nonlinear stochastic partial differential equations has gained less attention in the literature so far. In this work, we first show that the solution map is continuous. Then we introduce a notion of stability of exiting time. We provide  an example showing that  one cannot improve the stability of the exiting time and simultaneously improve the continuity of the dependence on initial data. Finally, we analyze the regularization effect of nonlinear noise in preventing blow-up. Precisely, we demonstrate that global existence holds true almost surely provided that the noise is strong enough.

\end{abstract}

\maketitle

%\section*{This is an unnumbered first-level section head}
%This is an example of an unnumbered first-level heading.
%
%%% The correct journal style for \specialsection is all uppercase; a known bug
%%% in amsart.cls prevents this, so input must be uppercase until it is fixed.
%%\specialsection*{This is a Special Section Head}
%\specialsection*{THIS IS A SPECIAL SECTION HEAD}
%This is an example of a special section head%
%%%%%%%%%%%%%%%%%%%%%%%%%%%%%%%%%%%%%%%%%%%%%%%%%%%%%%%%%%%%%%%%%%%%%%%%%
%\footnote{Here is an example of a footnote. Notice that this footnote
%text is running on so that it can stand as an example of how a footnote
%with separate paragraphs should be written.
%\par
%And here is the beginning of the second paragraph.}%
%%%%%%%%%%%%%%%%%%%%%%%%%%%%%%%%%%%%%%%%%%%%%%%%%%%%%%%%%%%%%%%%%%%%%%%%
%\tableofcontents

\section{Introduction and main results}
We consider  the following stochastic generalized  Camassa--Holm (CH) equation on $\R$:
\begin{equation}\label{gCH random dissipation}
	u_t-u_{xxt}+(k+2)u^ku_x-(1-\partial^2_{x})h(t,u)\dot{\mathcal W}=(k+1)u^{k-1}u_xu_{xx}+u^ku_{xxx},\ \ k\in\N_{>0}. 
\end{equation} 
In \eqref{gCH random dissipation}, $\W$  is a cylindrical Wiener process.

For  $h= 0$ and $k=1$, equation \eqref{gCH random dissipation} reduces to the deterministic 
CH equation given by
\begin{eqnarray}
	u_{t}-u_{xxt}+3uu_{x}=2u_{x}u_{xx}+uu_{xxx}.\label{CH}
\end{eqnarray}
Equation \eqref{CH} was introduced by
Fokas \& Fuchssteiner \cite{Fuchssteiner-Fokas-1981-PhyD}  to study completely integrable 
generalizations of the Korteweg-de Vries equation with bi-Hamiltonian structure. In \cite{Camassa-Holm-1993-PRL}, Camassa \& Holm proved that \eqref{CH} can be connected to the unidirectional propagation of shallow water waves over a flat bottom.  
Since then,  \eqref{CH} has been studied intensively, and we only mention a few related results here. 
The CH equation   exhibits both phenomena of soliton interaction (peaked soliton solutions) and wave breaking (the solution remains bounded while its slope becomes unbounded in finite time \cite{Constantin-Escher-1998-Acta}).

When $h= 0$ and $k=2$, equation \eqref{gCH random dissipation}  becomes the  so-called Novikov equation
\begin{equation}\label{Novikov}
	u_t-u_{xxt}+4u^2u_x=3uu_xu_{xx}+u^2u_{xxx},
\end{equation} 
which was derived  in \cite{Novikov-2009-JPA}. Equation \eqref{Novikov} also possesses a bi-Hamiltonian structure 
with an  infinite sequence of conserved quantities, and it admits peaked solutions \cite{Geng-Xue-2009-Nonlinearity}, as well as   multipeakon solutions  with explicit formulas \cite{Hone-etal-2009-DPDE}.
For the study of other deterministic instances of \eqref{gCH random dissipation},  
we refer to \cite{Himonas-Holliman-2014-ADE,
	Zhou-Mu-2013-JNS}.

When additional noise is included, 
as in \cite{Rohde-Tang-2021-JDDE}, the noise term can be used to account for the randomness arising from the energy exchange mechanisms. Indeed, in \cite{Lenells-Wunsch-2013-JDE,Wu-Yin-2008-SIAM}, the weakly dissipative term $(1-\partial^2_{x})(\lambda u)$ with $\lambda>0$ was added to the governing equations. In \cite{Rohde-Tang-2021-JDDE}, such weakly dissipative term is assumed to be time-dependent, nonlinear in $u$ and random. Therefore,   $(1-\partial^2_{x})h(t,u)\dot{\mathcal W} $ is proposed to describe random energy exchange mechanisms.

In this work,  we  consider the Cauchy problem for \eqref{gCH random dissipation} on the whole space $\R$.  Applying the operator $(1-\partial_{x}^2)^{-1}$ to  \eqref{gCH random dissipation}, we reformulate  the equation as  
\begin{equation} \label{SGCH problem}
	\left\{
	\begin{aligned}
		{\rm d}u+\left[u^k\partial_xu+F(u)\right]  \!\dd t&=h(t,u)\dd\W,\quad x\in\R, \ t>0,\  k\in\N_{>0},\\
		u(\omega,0,x)&=u_0(\omega,x),\quad x\in\R,
	\end{aligned} 
	\right.
\end{equation}
with
\begin{equation}\label{F decomposition}
	F(u):=F_1(u)+F_2(u)+F_3(u)\ \text{and}\ 
	\left\{
	\begin{aligned}
		&F_1(u):=(1-\partial_{x}^2)^{-1}\partial_x\left(u^{k+1}\right),\\
		&F_2(u):=\frac{2k-1}{2}(1-\partial_{x}^2)^{-1}\partial_x\left(u^{k-1}u_x^2\right),\\
		&F_3(u):=\frac{k-1}{2}(1-\partial_{x}^2)^{-1}\left(u^{k-2}u_x^3\right).
	\end{aligned} 
	\right.
\end{equation}
Here we remark that $F_3(u)$ in \eqref{F decomposition} will disappear for the CH case, (i.e., when $k=1$). The operator $(1-\partial_{x}^2)^{-1}$ in $F(\cdot)$ is understood as
\begin{align*} 
	\left[(1-\partial_{x}^2)^{-1}f\right](x)=\Big[\frac{1}{2} {\rm e}^{-|\cdot|}\star f\Big](x),
\end{align*}
where $\star$ stands for the convolution.

In this paper, regarding \eqref{SGCH problem}, we focus on the following issues:
\begin{itemize}
	\item Local well-posedness,  in the sense of Hadamard (existence, uniqueness and continuous dependence on initial data), and blow-up criterion of \eqref{SGCH problem}.
	
	\item Understanding the dependence on initial data, and in particular how continuous the solution map $u_0\mapsto u$ is.
	
	\item Analyzing the effect of noise vs blow-up of the deterministic counterpart of \eqref{SGCH problem}.
\end{itemize}

For the first and second issue, we refer to Theorems \ref{Local-WP} and \ref{Weak instability}, respectively. Extended remarks, explanations of difficulties,  and a review of literature are given in Remarks \ref{remark on LW-1}, \ref{remark on LW-2},  \ref{remark on instability 1} and \ref{remark on instability 2}.

The third question in our targets is on the impact of noise, which is one of the central questions in the study of stochastic partial differential equations (SPDEs). Regularization effects of noise have been  observed for many different models. For example, it is known that the well-posedness of linear stochastic transport equations with noise can be established under weaker hypotheses than its deterministic counterpart, cf. \cite{Flandoli-Gubinelli-Priola-2010-Invention}. Particularly, 
for the impact of linear noise in different models, we refer to \cite{Kroker-Rohde-2012-ANM,GlattHoltz-Vicol-2014-AP,Chen-Gao-2016-PA,Chen-Miao-Shi-2023-JMP,Rohde-Tang-2021-NoDEA,Tang-Wang-2023-CCM,Alonso-etal-2019-NODEA}.  

Notably, the existing results on regularization by noise are largely restricted to linear equations or linear noise. Hence we have particular interest in the nonlinear noise case.  Finding such noise is important
as it helps us to understand the stabilizing mechanisms of noise. This is the first step to characterize relevant noise which provides regularization effects for the CH-type equations.  In order to emphasize our ideas in a simple way, we only consider the noise as a 1-D Brownian motion in the current setting. That is, we  consider the case that 
$h(t,u) \dd \W
=q(t,u) \, {\rm d}W$, where $W$ is a standard 1-D Brownian motion and $q:[0,\infty)\times H^s\rightarrow H^s$ is a nonlinear function. Here we use the notation $q$ rather than $h$ because $h$ needs to be a Hilbert--Schmidt operator (see \eqref{stochastic integral}) to define the stochastic integral with respect to a cylindrical  Wiener process $\W$. Then we will focus on
\begin{equation}\label{SGCH non blow up Eq}
	\left\{\begin{aligned}
		&{\rm d}u+\left[u^ku_x+F(u)\right]\!   \dd t
		= q(t,u) \dd W,\quad x\in\R, \ t>0,\ k\in\N_{>0},\\
		&u(\omega,0,x)=u_{0}(\omega,x), \quad  x\in \R.
	\end{aligned} \right.
\end{equation}
In Theorem \ref{Non breaking}, we provide a sufficient condition on $q$ such that global existence can be guaranteed. We refer to Remark \ref{remark on global existence} for further remarks on Theorem \ref{Non breaking}.

Before we introduce the notations, definitions and assumptions, we recall some recent results on stochastic CH-type equations. For the stochastic CH type equation with multiplicative noise, we refer to  \cite{Tang-2018-SIMA,Rohde-Tang-2021-JDDE,Rohde-Tang-2021-NoDEA}, where global existence and wave breaking were studied in the periodic case, i.e., $x\in\T$.  In particular, when the noise is of transport type, we refer to \cite{Albeverio-etal-2021-JDE,Galimberti-etal-2022-arXiv,Holden-Karlsen-Pang-2021-JDE,Holden-Karlsen-Pang-2022-DCDS,Alonso-Rohde-Tang-2021-JNLS}.  We also refer to 
\cite{Chen-Duan-Gao-2021-CMS,Chen-Duan-Gao-2021-PhyD,Ren-Tang-Wang-2020-Arxiv} for more results in stochastic CH type equations.

\subsection{Notations}
We begin by introducing some notations.  Let $({\rm \Omega},\{\F_t\}_{t\ge 0}, \p)$ be  a right-continuous complete filtration probability space.   Formally, we consider a separable Hilbert space $\U$ and let $\{e_n\}$  be a complete orthonormal basis of $\U$. Let  $\{W_n\}_{n\geq1}$ be a sequence of mutually independent standard 1-D Brownian motions on $({\rm \Omega},\{\mathcal{F}_t\}_{t\geq0},\p)$. Then we define the cylindrical Wiener process $\W$ as
\begin{equation}\label{define cylindrical process}
	\W:=\sum_{n=1}^\infty W_ne_n.
\end{equation}
Let $\X$ be a separable Hilbert space.  $\LL_2(\U; \X)$ stands for the Hilbert–Schmidt operators from $\U$ to $\X$. If $Z\in L^2 ({\rm \Omega}; L^2_{\rm loc}([0,\infty);\LL_2(\U; \X)))$ is progressively measurable, then the integral
\begin{equation}\label{stochastic integral}
	\int_0^t Z \dd \W:=\sum_{n=1}^\infty\int_0^t Z e_n \dd W_n
\end{equation}
is a well-defined $\X$-valued continuous
square-integrable martingale (see \cite{Breit-Feireisl-Hofmanova-2018-Book,Gawarecki-Mandrekar-2010-Springer} for example). Throughout the paper, when a stopping time is defined, we set $\inf\emptyset:=\infty$ by convention.

For $s\in\R$, the differential operator  $D^s:=(1-\partial_{x}^2)^{s/2}$ is defined by
$\widehat{D^sf}(\xi)=(1+\xi^2)^{s/2}\widehat{f}(\xi)$, where $\widehat{f}$ denotes the  Fourier transform  of $f$. The Sobolev space $H^s(\R)$ is defined as
\begin{align*}
	H^s(\R):=\left\{f\in L^2(\R):\|f\|_{H^s(\R)}^2:=\int_{\R}(1+|\xi|^2)^s|\widehat{f}(\xi)|^2  \dd \xi<+\infty\right\},
\end{align*}
and the inner product    on $H^s(\R)$   
is $(f,g)_{H^s}:=(D^sf,D^sg)_{L^2}.$
In the sequel, for simplicity, we will drop $\R$ if there is no ambiguity. We will use $\lesssim $ to denote estimates that hold up to some universal \textit{deterministic} constant which may change from line to line but whose meaning is clear from the context. For linear operators $A$ and $B$,  $[A,B]:=AB-BA$ is the commutator of $A$ and $B$.

\subsection{Definitions and assumptions}
We first  make the precise notion of a  solution to \eqref{SGCH problem}.
\begin{Definition}\label{Definition of solution}
	Let $({\rm \Omega}, \{\mathcal{F}_t\}_{t\geq0},\p, \W)$ be a fixed in advance. Let $s>3/2$, $k\in\N_{>0}$ and $u_0$ be an $H^s$-valued $\mathcal{F}_0$-measurable random variable.
	\begin{enumerate}[label={$\bf (\arabic*)$}]
		\item A local  solution to \eqref{SGCH problem} is a pair $(u,\tau)$, where $\tau$ is a stopping time satisfying $\p\{\tau>0\}=1$ and
		$u:{\rm \Omega}\times[0,\infty)\rightarrow H^s$  is an $\mathcal{F}_t$-predictable $H^s$-valued process satisfying
		\begin{equation*} 
			u(\cdot\wedge \tau)\in C([0,\infty);H^s)\ \ \pas,
		\end{equation*}
		and for all $t>0$,
		\begin{equation*}
			u(t\wedge \tau)-u(0)+\int_0^{t\wedge \tau}
			\left[u^k\partial_xu+F(u)\right]      \dd t'
			=\int_0^{t\wedge \tau}h(t',u) \dd \W\ \ \pas
		\end{equation*}
		\item The local solutions are said to be unique, if given any two pairs of local solutions $(u_1,\tau_1)$ and $(u_2,\tau_2)$ with $\p\left\{u_1(0)=u_2(0)\right\}=1,$ we have
		\begin{equation*}
			\p\left\{u_1(t,x)=u_2(t,x),\  (t,x)\in[0,\tau_1\wedge\tau_2]\times \R\right\}=1.
		\end{equation*}
		\item Additionally, $(u,\tau^*)$ is called a maximal solution to \eqref{SGCH problem} if $\tau^*>0$ almost surely and if there is an increasing sequence $\tau_n\rightarrow\tau^*$ such that for any $n\in\N$, $(u,\tau_n)$ is a solution to \eqref{SGCH problem} and on the set $\{\tau^*<\infty\}$, we have
		\begin{equation*} 
			\sup_{t\in[0,\tau_n]}\|u\|_{H^s}\geq n.
		\end{equation*}
		\item If $(u,\tau^*)$ is a maximal solution and
		$\tau^*=\infty$ almost surely, then we say that the solution exists globally.
	\end{enumerate}
\end{Definition}

Motivated by \cite{Rohde-Tang-2021-JDDE,Tang-2022-arXiv}, we introduce the  concept on 
stability of  exiting time in Sobolev spaces. Exiting time, as its name would suggest, is defined as the time when solution leaves a certain range. 

\begin{Definition}[Stability of exiting time]\label{Definition stability of exiting time}
	Let $({\rm \Omega}, \{\mathcal{F}_t\}_{t\geq0},\p, \W)$ be fixed, $s>3/2$ and $k\in\N_{>0}$. Let $u_0$ be an $H^s$-valued $\mathcal{F}_0$-measurable random variable such that $\E\|u_0\|_{H^s}^2<\infty$. Assume that $\{u_{0,n}\}$ is a sequence of $H^s$-valued $\mathcal{F}_0$-measurable random variables satisfying $\E\|u_{0,n}\|_{H^s}^2<\infty$. For each $n$, let $u$ and $u_n$ be the unique solutions to \eqref{SGCH problem}, as in Definition \ref{Definition of solution}, with initial values $u_0$ and $u_{0,n}$, respectively. For any $R>0$, define the $R$-exiting times 
	\begin{equation*}
		\tau_n^R:=\inf\{t\geq 0:\|u_n\|_{H^s}>R\},\ \ \tau^R:=\inf\{t\geq 0:\|u\|_{H^s}>R\}.
	\end{equation*}
	Now we define the following properties on stability:
	\begin{enumerate}[label={$\bf (\arabic*)$}]
		\item If $u_{0,n}\rightarrow u_0$ in $H^{s}$ $\pas$  implies that
		\begin{align}
			\lim_{n\rightarrow\infty}\tau^R_{n}=\tau^R\ \  \pas,\label{hitting time convergence}
		\end{align}
		then the $R$-exiting time of $u$ is said to be stable.
		
		\item If $u_{0,n}\rightarrow u_0$ in $H^{s'}$ for all $s'<s$ almost surely  implies that \eqref{hitting time convergence} holds true,
		the $R$-exiting time of $u$ is said to be strongly stable.
	\end{enumerate}
\end{Definition}

Our main results rely on the following assumptions concerning 
the noise coefficient $h(t,u)$ in \eqref{gCH random dissipation}. 

\begin{Hypothesis}\label{Assumption-1}
	For $s>1/2$, we assume that $h:[0,\infty)\times H^s\ni (t,u)\mapsto h(t,u)\in \LL_2(\U; H^s)$  is measurable and satisfies the following conditions:
	
	\begin{enumerate}[label={$\bf H_1(\arabic*)$}]
		
		\item\label{Ass1-growth assumption} There is a non-decreasing   function $f(\cdot):[0,+\infty)\rightarrow[0,+\infty)$ such that for any  $u\in H^s$ with $s>3/2$, we have the following  growth condition
		\begin{align*}
			\sup_{t\ge0}\|h(t,u)\|_{\LL_2(\U; H^s)}\leq f(\|u\|_{W^{1,\infty}}) (1+\|u\|_{H^s}).
		\end{align*} 
		
		\item \label{Ass1-Lip s>3/2} There is a non-decreasing  function $g_1(\cdot):[0,\infty)\rightarrow[0,\infty)$ such that for all $N\ge1$,
		
		\begin{equation*}
			\sup_{t\ge0,\,\|u\|_{H^s},\,\|v\|_{H^s}\le N}\left\{{\bf 1}_{\{u\ne v\}}  \frac{\|h(t,u)-h(t,v)\|_{\LL_2(\U, H^s)}}{\|u-v\|_{H^s}}\right\} \le g_1(N),\ \ s>3/2.
		\end{equation*}

		\item \label{Ass1-Lip s>1/2} There is a non-decreasing function $g_2(\cdot):[0,\infty)\rightarrow[0,\infty)$ such that  for all $N\ge1$ and $3/2\geq s>1/2$,
		\begin{equation*} 
			\sup_{t\ge0,\,\|u\|_{H^{s+1}},\,\|v\|_{H^{s+1}}\le N}\left\{{\bf 1}_{\{u\ne v\}}  \frac{\|h(t,u)-h(t,v)\|_{\LL_2(\U, H^s)}}{\|u-v\|_{H^s}}\right\} \le g_2(N).
		\end{equation*}
	\end{enumerate}
\end{Hypothesis}
Here we outline \ref{Ass1-Lip s>3/2} is the classical local Lipschitz condition. \ref{Ass1-Lip s>1/2} is needed to prove uniqueness in Lemma \ref{Uniqueness bounded data}. Indeed, if one finds two solutions $u,v\in H^s$ to \eqref{SGCH problem}, one can only estimate $u-v$ in $H^{s'}$ for $s'\leq s-1$ because the term $u^ku_x$ loses one derivative. We refer to Remark \ref{remark on LW-1} for more details.

\begin{Hypothesis}\label{Assumption-2} 
	When we consider \eqref{SGCH problem} in Sect. \ref{sect:weak instability}, we assume that there is a real number $\rho_0\in(1/2,1)$ such that for $s\ge\rho_0$, $h:[0,\infty)\times H^s\ni (t,u)\mapsto h(t,u)\in \LL_2(\U; H^s)$ is  measurable. Besides, we suppose the following:

	\begin{enumerate}[label={$\bf H_2(\arabic*)$}]
		
		\item\label{Ass2-growth-Lip} There exists a non-decreasing function $l(\cdot):[0,+\infty)\rightarrow[0,+\infty)$ such that for any  $u\in H^s$ with $s>3/2$,
		\begin{align*}
			\sup_{t\ge0}\|h(t,u)\|_{\LL_2(\U; H^s)}\leq l(\|u\|_{W^{1,\infty}})\|u\|_{H^s},
		\end{align*} 
		and  \ref{Ass1-Lip s>3/2} holds.
		\item\label{Ass2-rho-0} There is a non-decreasing   function $g_3(\cdot):[0,+\infty)\rightarrow[0,+\infty)$ such that for all $N\ge1$,
		\begin{align}
			\sup_{t\ge0,\,\|u\|_{H^s}\le N}\|h(t,u)\|_{\LL_2(\U; H^{\rho_0})}\leq g_3(N) {\rm e}^{-\frac{1}{\|u\|_{H^{\rho_0}}}},\ \ s>3/2,\label{Ass2-small assumption}
		\end{align}
		and
		\begin{equation*} 
			\sup_{t\ge0,\,\|u\|_{H^{\rho_0}},\,\|v\|_{H^{\rho_0}}\le N}\left\{{\bf 1}_{\{u\ne v\}}  \frac{\|h(t,u)-h(t,v)\|_{\LL_2(\U, H^{\rho_0})}}{\|u-v\|_{H^{\rho_0}}}\right\} \le g_3(N).
		\end{equation*}
	\end{enumerate}

\end{Hypothesis}

We remark here that \eqref{Ass2-small assumption} means that there is a $\rho_0\in(1/2,1)$ such that, if $u_n$ is bounded in $H^s$ and $u_n$ tends to zero in the topology of $H^{\rho_0}$ as $n$ tends to $\infty$, then $\|h(t,u_n)\|_{\LL_2(\U; H^{\rho_0})}$ tends to zero exponentially as $n$ tends to $\infty$. Examples of such noise structure can be found in Sect. \ref{Example:Weak instability}.

As for the regularization effect of noise, we impose the following condition on $q$ in \eqref{SGCH non blow up Eq}:

\begin{Hypothesis}\label{Assumption-q}
	We assume that when $s>3/2$,  $q:[0,\infty)\times H^s\ni (t,u)\mapsto q(t,u)\in H^s$ is measurable. Define the set $\mathcal{V}$ as a subset of $ C^2([0,\infty);[0,\infty))$ such that
	\begin{equation*}
		\mathcal{V}:=\left\{V(0)=0,\  V'(x)>0, \ V''(x)\le 0 \ \text{and} \ \lim_{x\to\infty} V(x)=\infty\right\}.
	\end{equation*} 
	Then we assume  the following:
	
	\begin{enumerate}[label={$\bf H_3(\arabic*)$}]
		\item \label{Ass3-cut-off}
		There is a non-decreasing  function $g_4(\cdot):[0,+\infty)\rightarrow[0,+\infty)$ such that for any  $u\in H^s$ with $s>3/2$, we have the following  growth condition
		\begin{align*}
			\sup_{t\ge0}\|q(t,u)\|_{H^s}\leq g_4(\|u\|_{W^{1,\infty}}) (1+\|u\|_{H^s}).
		\end{align*}

		\item\label{Ass3-Lip}  $q(\cdot,u)$ is bounded for all $u\in H^s$ and there is a non-decreasing  function $g_4(\cdot):[0,\infty)\rightarrow[0,\infty)$, such that 
		\begin{equation*}
			\sup_{t\ge0,\,\|u\|_{H^{s}},\,\|v\|_{H^{s}}\le N}\left\{{\bf 1}_{\{u\ne v\}}  \frac{\|q(t,u)-q(t,v)\|_{H^s}}{\|u-v\|_{H^s}}\right\} \le g_4(N),\ \ N\ge 1,\  s>3/2.
		\end{equation*}

		\item\label{Ass3-growth} 
		There is a $V\in\mathcal{V}$ and constants $N_1, N_2>0$ such that for all $(t,u)\in [0,\infty)\times H^s$ with $s>3/2$,
		\begin{align*} 
			\mathcal{H}_{s}(t,u)
			\leq \, &   N_1
			-N_2 \frac{\left\{V'(\|u\|^2_{H^s})\left|\left(q(t,u), u\right)_{H^s}\right|\right\}^2}{ 1+V(\|u\|^2_{H^s})},\ \ 
		\end{align*}
		where
		\begin{align*}
			\,&\mathcal{H}_{s}(t,u)\\
			:=\, &V'(\|u\|^2_{H^{s}})\Big\{
			2\lambda_s \|u\|^k_{W^{1,\infty}}\|u\|^2_{H^{s}}+\|q(t,u)\|^2_{H^{s}}
			\Big\} +
			2V''(\| u\|^2_{H^{s}})\left|\left( q(t,u), u\right)_{H^{s}}\right|^2
		\end{align*}
		and $\lambda_s>0$ is the constant given in Lemma \ref{uux+F u Hs inner product} below.

	\end{enumerate}

\end{Hypothesis}

Examples of the noise structure satisfying Hypothesis {\rm \ref{Assumption-q}} can be found in Sect. \ref{Example:large noise}.

\subsection{Main results and remarks} Now we summarize our major contributions providing  proofs later in the remainder of the paper.

\begin{Theorem}\label{Local-WP} Let $s>3/2$, $k\geq1$ and let $h(t,u)$ satisfy Hypothesis {\rm\ref{Assumption-1}}.   Assume that $u_0$ is an $H^s$-valued $\mathcal{F}_0$-measurable random variable satisfying $\E\|u_0\|^2_{H^s}<\infty$. Then

	\begin{enumerate} [label={\bf (\roman*)}]
		
		\item \label{Local-WP-existence} $($Existence and uniqueness$)$ There is a unique local solution $(u,\tau)$ to \eqref{SGCH problem} in the sense of Definition \ref{Definition of solution} with
		\begin{equation}\label{L2 moment bound}
			\E\sup_{t\in[0,\tau]}\|u(t)\|^2_{H^s}<\infty.
		\end{equation} 
		
		\item\label{Local-WP-blow-up} $($Blow-up criterion$)$ The local solution $(u,\tau)$ can be extended to a unique maximal solution $(u,\tau^*)$ with
		\begin{equation}\label{Blow-up criterion}
			\textbf{1}_{\left\{\limsup_{t\rightarrow \tau^*}\|u(t)\|_{H^{s}}=\infty\right\}}=\textbf{1}_{\left\{\limsup_{t\rightarrow \tau^*}\|u(t)\|_{W^{1,\infty}}=\infty\right\}}\ \pas
		\end{equation}
		
		\item\label{Local-WP-dependence}  $($Stability for almost surely bounded initial data$)$ 
		Assume additionally that $u_0\in L^\infty({\rm \Omega};H^s)$.
		Let $v_0\in L^\infty({\rm \Omega};H^s)$ be another $H^s$-valued $\mathcal{F}_0$-measurable random variable.   For any $T>0$ and any $\epsilon>0$, there is a $\delta=\delta(\epsilon,u_0,T)>0$ such that if
		\begin{equation}\label{stability condition L-infty}
			\|u_0-v_0\|_{L^\infty({\rm \Omega};H^s)}<\delta,
		\end{equation}
		then there is a stopping time $\tau\in(0,T]$ $\pas$  and
		\begin{equation}\label{stability result}
			\E\sup_{t\in[0,\tau]}\|u(t)-v(t)\|^2_{H^s}<\epsilon,
		\end{equation}
		where $u$ and $v$ are the solutions to \eqref{SGCH problem} with initial data $u_0$ and $v_0$, respectively.
		
	\end{enumerate}
\end{Theorem}

\begin{Remark} \label{remark on LW-1}
	Existence and uniqueness  have been studied for abundant SPDEs. 
	In many works,  the authors did not address the continuous dependence on initial data. In this work, our Theorem \ref{Local-WP} provides a local well-posedness result in the sense of Hadamard including the continuous dependence on initial data. Moreover, a blow-up criterion is also obtained.  We refer to \cite{Fedrizzi-Flandoli-2011-Sto,Marinelli-Prevot-Rockner-2010-JFA,Chen-Pang-2021-JFA} for the study about the dependence on the initial data for cases that solutions to the target problems exist globally. However, it is necessary to point out that almost nothing is known on  the analysis for dependence on initial data for SPDEs whose solutions may blow up in finite time.
	
	The key difficulty for such a case is as follows: on one hand, if solutions to a nonlinear stochastic partial differential equation (SPDE) blow up in finite time, it is usually very difficult  to obtain the lifespan estimates.
	On the other hand, we have to find a positive time $\tau$ to obtain an inequality like  \eqref{stability result}. In addition,  the target problem \eqref{SGCH problem} is more difficult because the classical It\^{o} formulae are \textit{not} applicable. Indeed, for $u_0\in H^s$, we can \textit{only} know $u\in H^s$ because this is a transport type equation, then $u^ku_x\in H^{s-1}$. However,  the inner product $ \left(u^ku_x, u\right)_{H^s}$ appears if one uses the It\^{o} formula in a Hilbert space  (cf. \cite[Theorem 2.10]{Gawarecki-Mandrekar-2010-Springer})  and the dual product ${}_{H^{s-1}}\langle u^ku_x, u\rangle_{H^{s+1}}$  appears in the It\^{o} formula  under a  Gelfand triplet (cf.  \cite[Theorem I.3.1]{Krylov-Rozovskiui-1979-chapter}). Since we only have $u\in H^s$ and $u^ku_x\in H^{s-1}$, neither of them  are well-defined. Likewise,  when we consider the $H^s$-norm for the difference between two solutions $u,v\in H^s$ to \eqref{SGCH problem}, we will have to handle $(u^ku_x-v^kv_x,u-v)_{H^s}$, which gives rise to control either  $\|u\|_{H^{s+1}}$ or $\|v\|_{H^{s+1}}$.

\end{Remark}

\begin{Remark} \label{remark on LW-2}
	
	Now we list some  technical remarks on the statements of Theorem \ref{Local-WP}.

	\begin{enumerate}[label={ $\bf (\arabic*)$}]
		\item Our proof for \ref{Local-WP-existence} in Theorem \ref{Local-WP} is motivated by the recent results in \cite{Tang-Yang-2022-AIHP}. For the convenience of the reader,  here we also give a brief comparison between our approach  and the framework employed in many previous works. 
		
		\begin{itemize}
			\item We first briefly review the martingale approach used to prove existence of nonlinear SPDEs. Roughly speaking, in searching for a solution to a nonlinear SPDE in some space $\X$, the martingale approach, as its name would suggest,  includes obtaining martingale solution first and then establishing (pathwise) uniqueness to obtain the (pathwise)  solution. To begin with,  one needs to approximate the equation and establish uniform estimate. For nonlinear problems, one may have to add a \textit{cut-off} function to cut the nonlinear parts growing in  some space $\mathcal{Z}$ with $\X\hookrightarrow\mathcal{Z}$
			(such choice of $\mathcal{Z}$  depends on concrete problems). As far as we know, the technique of \textit{cut-off} first appears in \cite{deBouard-Debussche-1999-CMP} for the stochastic {S}chr\"{o}dinger equation. This \textit{cut-off} enables us to split the expectation of nonlinear terms, and then the  $L^2({\rm \Omega}; \X)$ estimate can be closed. For example, for \eqref{SGCH problem}, the estimate for $\E\|u\|^2_{H^s}$ will give rise to   $\E\left(\|u\|_{\Wlip}\|u\|^2_{H^s}\right)$,  and hence we need to add a function to cut $\|\cdot\|_{\Wlip}$. 
			With this additional \textit{cut-off}, we need to consider the \textit{cut-off} version of the problem first and remove it then. The first main step in the martingale approach is finding a martingale solution. Usually, this can be done by first obtaining tightness of the measures defined by the approximative solutions in some space $\Y$, and then using Prokhorov's Theorem and Skorokhod's Theorem to obtain the convergence in $\Y$. Since $\X$ is usually  infinite dimensional (usually, $\X$ is a Sobolev space), to obtain tightness, it is required that $\X$ is compactly embedded into $\Y$, i.e, $\X\hookrightarrow\hookrightarrow\Y$. This brings another requirement to specify $\mathcal{Z}$, that is, $\Y\hookrightarrow\mathcal{Z}$.
			Otherwise, taking limits will \textit{not} bring us back to the \textit{cut-off} problem due to the additional \textit{cut-off} term $\|\cdot\|_{\mathcal{Z}}$ (in some cases, the choice of $\mathcal{Z}$ may only give rise to a semi-norm and here we use this notation $\|\cdot\|_{\mathcal{Z}}$ only for simplicity). Usually, in bounded domains, it is not difficult to pick $\Y$ and $\mathcal{Z}$ such that $\X\hookrightarrow\hookrightarrow\Y\hookrightarrow\mathcal{Z}$ (Sobolev spaces enjoy compact embeddings in bounded domains), see for example \cite{Brzezniak-Ondrejat-2007-JFA,
				Debussche-Glatt-Temam-2011-PhyD,GlattHoltz-Vicol-2014-AP,Tang-2018-SIMA,Alonso-Rohde-Tang-2021-JNLS}.  
			In unbounded domains, the difficulty lies in the choice of $\Y$ and $\mathcal{Z}$ such that $\X\hookrightarrow\hookrightarrow\Y\hookrightarrow\mathcal{Z}$. We refer to
			\cite{Brzezniak-Motyl-2013-JDE,Brzezniak-Motyl-2019-SIMA} for fluid models with certain cancellation properties (for example, divergence free) and linear growing noise. However,  it is difficult to achieve this for SPDEs with 
			general nonlinear terms and nonlinear noise. For instance, the \textit{cut-off}
			in our case  will have to involve $\|\cdot\|_{\mathcal{Z}}=\|\cdot\|_{W^{1,\infty}}$ (see \ref{Ass1-growth assumption} and \eqref{approximate problem}). Even though we can get the convergence in $H_{\rm loc}^{s'}$ with some $\frac32<s'<s$, it is still not clear whether the convergence holds true in $W^{1,\infty}$, and this is because \textit{local} convergence can \textit{not} control a \textit{global}
			object $\|\cdot\|_{\Wlip}$. Therefore, technically speaking, nonlinear SPDEs are more non-local than its deterministic counterpart.

			\item  Due to the above unsolved technical issue, the martingale approach is difficult to apply in our problem and we will try to prove convergence
			directly, which is motivated by \cite{Li-Liu-Tang-2021-SPA,Tang-Yang-2022-AIHP} (see also \cite{Tang-Wang-2023-CCM,Tang-2022-arXiv,Tang-Wang-2022-arXiv} for recent developments). Generally speaking,  we will analyze the difference between two approximative  solutions and \textit{directly} find a space $\Y$ such that $\X\hookrightarrow\Y\hookrightarrow\mathcal{Z}$ and convergence (up to a subsequence) holds true in $\Y$. The difficult part is finding convergence in $\Y$ \textit{without} compactness $\X\hookrightarrow\hookrightarrow\Y$ (compared to the martingale approach, tightness comes from the compact embedding $\X\hookrightarrow\hookrightarrow\Y$).
			In this paper, the target path space is $C([0,T];H^{s})=\X$, and we are able to prove convergence  (up to a subsequence) in $C([0,T];H^{s-\frac32})=\Y$ directly. 
			After taking limits to obtain a solution, one can improve the regularity to $H^s$ again, and the technical difficulty in this step  is to prove the time continuity of the solution because the classical It\^{o} formula is \textit{not} applicable (see in Remark \ref{remark on LW-1}).  To overcome this difficulty, we apply a mollifier $J_\e$ to equation and estimate $\E \|J_\e u\|^2_{H^s}$ first (see \eqref{TnX 2}). 
			We also remark that the techniques in removing the \textit{cut-off} have been used in \cite{GlattHoltz-Ziane-2009-ADE,Breit-Feireisl-Hofmanova-2018-Book,Tang-Wang-2023-CCM}. Here we formulate such a technical result in Lemma \ref{cut-combine} in an abstract way. 
		\end{itemize}

		\item Now we give a remark on \ref{Local-WP-dependence} in Theorem \ref{Local-WP}. For the question on dependence on initial data, there are some delicate differences between the stochastic and the deterministic case.
		In the deterministic counterpart of \eqref{SGCH problem}, due to the lifespan estimate (see \eqref{solution bound T u_0} for instance), for given $u_0\in H^s$, it can be shown that if $\|u_0-v_0\|_{H^s}$ is small enough, then there is a $T>0$ depending on $u_0$ such that $\sup_{t\in[0,T]}\|u(t)-v(t)\|^2_{H^s}$ is also small. In stochastic setting, since existence and uniqueness are obtained in the framework of $L^2({\rm \Omega};H^s)$, it is therefore very natural to expect that, for given $u_0\in L^2({\rm \Omega};H^s)$, if $\E\|u_0-v_0\|^2_{H^s}$ is small enough, then for some almost surely positive $\tau$ depending $u_0$, $\E\sup_{t\in[0,\tau]}\|u(t)-v(t)\|^2_{H^s}$ is also small. However, so far we have only proved it with assuming the smallness of $\|u_0-v_0\|_{L^\infty({\rm \Omega};H^s)}$.  Since $L^\infty({\rm \Omega};H^s)$ can be viewed as being \textit{less} random than $L^2({\rm \Omega};H^s)$,  one may roughly conclude that  what the solution map needs to be continuous/stable (the initial data and its perturbation are $L^\infty({\rm \Omega};H^s)$) is more ``picky" in determinism than what the existence of such a solution map requires (existence and uniqueness guarantee that a solution map can be defined). For the technical difficulties involved, we have the following explanations:

		\begin{itemize}
			\item As is mentioned in Remark \ref{remark on LW-1}, when we estimate the $H^s$-norm for the difference between two solutions $u$ and $v$, $H^{s+1}$-norm will appear. Hence, we have to use smooth approximations to make the analysis valid. More precisely, we approximate $u$ and $v$ by smooth process $u_\e$ and $v_\e$ and consider
			\begin{equation}\label{continuity 3 terms}
				\|u-v\|_{H^s}\leq
				\|u-u_{\varepsilon}\|_{H^s}+\|u_{\varepsilon}-v_{\varepsilon}\|_{H^s}+\|v_{\varepsilon}-v\|_{H^s}.
			\end{equation}
			Then all terms can be estimated because $\|u_\e\|_{H^{s+1}}$ and $\|v_\e\|_{H^{s+1}}$ make sense. Here we refer to Remark \ref{Continuous dependence reamrk} for more details on the construction of such an approximation.
			
			\item In dealing with the above three terms in the stochastic case,  two sequences of stopping times (exiting times) are needed to control $\|u_\e\|_{H^s}$ and $\|v_\e\|_{H^s}$ (see \eqref{tau u v stability} below).
			Since we aim at obtaining $\tau>0$ almost surely in \eqref{stability result} (otherwise the difference between two solutions on the set $\{\tau=0\}$ can \textit{not} be measured), we will have to guarantee that those stopping times used in bounding $\|u_\e\|_{H^s}$ and $\|v_\e\|_{H^s}$ have positive lower bounds almost surely. Up to now, we have only achieved this for initial values belonging to $L^\infty({\rm \Omega};H^s)$.  We also remark that this is different from the proof for existence. In the proof for existence,  $u_\e$  exists on a common interval $[0,T]$ for all $\e$ and enjoys a uniform-in-$\e$ estimate \eqref{E r e}, hence we can get rid of stopping times in convergence (from \eqref{Cauchy H s-3/2} to \eqref{convergence ue a.s.}). Here we do \textit{not} have such common existence interval due to the lack of a 
			lifespan estimate, which is a significant difference between the stochastic and the deterministic cases. Indeed, we can easily find the lifespan estimate for the deterministic counterpart of \eqref{SGCH problem} (see \eqref{solution bound T u_0} below).

			\item Moreover, even if the above issue can be handled, in dealing with the three terms in \eqref{continuity 3 terms}, we are confronted with $\frac{1}{\e^2}\E\|u_0- v_0\|^2_{H^{s'}}\|u_0\|^2_{H^{s}}$ for some suitably chosen $s'$ (cf. \eqref{u-en v-en difference}). After $\e$ is fixed, the smallness of $\E\|u_0-v_0\|^2_{H^s}$ is \textit{not} enough to control  $\frac{1}{\e^2}\E\|u_0- v_0\|^2_{H^{s'}}\|u_0\|^2_{H^{s}}$, either. We use the $L^\infty({\rm \Omega};H^s)$ condition  to take $\|u_0\|^2_{H^{s}}$ out of $\frac{1}{\e^2}\E\|u_0- v_0\|^2_{H^{s'}}\|u_0\|^2_{H^{s}}$. In deterministic case, no expectation is involved, $\frac{1}{\e^2}\|u_0- v_0\|^2_{H^{s'}}\|u_0\|^2_{H^{s}}$ can be controlled by $\|u_0- v_0\|^2_{H^{s}}$.
		\end{itemize}

	\end{enumerate}

\end{Remark}

Roughly speaking, 
\ref{Local-WP-dependence} in Theorem \ref{Local-WP} means that for any fixed $u_0\in L^\infty({\rm \Omega};H^s)$ and any $T>0$, if $\|u_0-v_0\|_{L^\infty({\rm \Omega};H^s)}\to 0$, then
$$
\exists\, \tau\in(0,T] \ \pas  \ \text{such that}\  \E\|u(\cdot\wedge\tau)-v(\cdot\wedge\tau)\|^2_{C([0,T];H^s)}\to 0,
$$
where $u,v$ are solutions corresponding to $u_0$, $v_0$, respectively.
Below we will study this issue quantitatively. The next result addresses at least   a partially negative answer.
\begin{Theorem}[Weak instability]\label{Weak instability}
	Let $s>5/2$ and $k\ge1$. If $h$ satisfies Hypothesis {\rm\ref{Assumption-2}}, then at least one of the following properties holds true:
	\begin{enumerate}[label={\bf (\roman*)}]
		\item For any $R\gg 1$, the $R$-exiting time is not strongly stable for the zero solution to \eqref{SGCH problem} in the sense of Definition \ref{Definition stability of exiting time};
		\item There is a $T>0$ such that the solution map $u_0\mapsto u$ defined by \eqref{SGCH problem} is not uniformly continuous as a map from $L^\infty({\rm \Omega};H^s)$ into $L^1({\rm \Omega};C([0,T];H^s))$. 
		More precisely, there exist two sequences of solutions $u^{1,n}$ and $u^{2,n}$, and two sequences of stopping time $\tau_{1,n}$ and $\tau_{2,n}$, such that
		\begin{itemize}
			\item For $i=1,2$, $\p\{\tau_{i,n}>0\}=1$ for each $n>1$. Besides,
			\begin{equation}\label{tau 1 2 n}
				\lim_{n\rightarrow\infty}\tau_{1,n}=\lim_{n\rightarrow\infty}\tau_{2,n}=\infty\ \ \pas
			\end{equation}
			\item For $i=1,2$, $u^{i,n}\in C([0,\tau_{i,n}];H^s)$ $\pas$, and
			\begin{equation}\label{sup u}
				\E\left(\sup_{t\in[0,\tau_{1,n}]}\|u^{1,n}(t)\|_{H^s}+\sup_{t\in[0,\tau_{2,n}]}\|u^{2,n}(t)\|_{H^s}\right)\lesssim 1.
			\end{equation}
			\item At initial time $t=0$, for any $p\in[1,\infty]$,
			\begin{equation}\label{same initail data}
				\lim_{n\rightarrow\infty}\|u^{1,n}(0)-u^{2,n}(0)\|_{L^p({\rm \Omega};H^s)}=0.
			\end{equation}
			\item When $t>0$,
			\begin{align}\label{sup sin t}
				\liminf_{n\rightarrow\infty}\E\sup_{t\in[0,T\wedge\tau_{1,n}\wedge\tau_{2,n}]}
				\|u^{1,n}(t)&-u^{2,n}(t)\|_{H^s}\notag\\
				\gtrsim &
				\left\{\begin{aligned}
					&\sup_{t\in[0,T]}|\sin(t)|,\ \text{if}\ k\ \text{is odd},\\
					& \sup_{t\in[0,T]}\big|\sin\big(\frac{t}{2}\big)\big|,\ \text{if}\ k\ \text{is even}.
				\end{aligned}\right.
			\end{align}
		\end{itemize}
	\end{enumerate}
\end{Theorem}

\begin{Remark}\label{remark on instability 1}
	We first briefly outline the  main difficulties encountered in the proof for Theorem \ref{Weak instability} and the main strategies we used.
	\begin{enumerate}[label={ $\bf (\arabic*)$}]
		
		\item Because we can \textit{not} get an explicit expression of the solution to \eqref{SGCH problem}, to obtain \eqref{sup sin t}, we will construct two sequences of  approximative solutions $\{u_{m,n}\}$ ($m\in\{1,2\}$)  such that the actual solutions $\{u^{m,n}\}$ with $u^{m,n}(0)=u_{m,n}(0)$  satisfy
		\begin{align}
			\lim_{n\rightarrow\infty}\E\sup_{[0,\tau_{m,n}]}\|u^{m,n}-u_{m,n}\|_{H^s}
			=0,\label{non uniform remark equ}
		\end{align}
		where $u^{m,n}$ exists at least on $[0,\tau_{m,n}]$.  Then, one can establish \eqref{sup sin t} by estimating $\{u_{m,n}\}$ rather than $\{u^{m,n}\}$.
		We also remark that  the construction of approximative solution $u_{m,n}$ for $x\in\R$ is more difficult than the construction of approximative solution for  $x\in\T$ (see \cite{Rohde-Tang-2021-JDDE}) since the approximative solution involves both high and low frequency parts (high frequency part is already enough for  the case $x\in\T$, cf. \cite{Rohde-Tang-2021-JDDE,Tang-Yang-2022-AIHP}). The key point is that we need to guarantee $\inf_{n}\tau_{m,n}>0$ almost surely in dealing with \eqref{non uniform remark equ}. Hence we are confronted with a common difficulty in SPDEs again, that is, the lack of   lifespan estimate. In
		deterministic cases, one can easily obtain the lifespan estimate, which enables us to find a common interval $[0,T]$ such that all actual solutions exist on $[0,T]$ (see for example Lemma  \ref{lemma ul Hr}). In the  stochastic case, so far we have \textit{not} been able to prove this.

		\item To settle the above difficulty, we observe that  the bound $\inf_{n}\tau_{m,n}>0$ can be connected to the stability property of the exiting time (see Definition \ref{Definition stability of exiting time}). 
		The condition that the $R_0$-exiting time is strongly stable at the zero solution will be used to provide a common existence time $T>0$ such that for all $n$, $u^{m,n} $ exists up to $T$ (see Lemma \ref{exiting time infty lemma} below).  Therefore, to prove Theorem \ref{Weak instability}, we will show that, if the $R_0$-exiting time is strongly stable at the zero solution  for some $R_0\gg1$,  then the solution map $u_0\mapsto u$ defined by \eqref{SGCH problem} can \textit{not} be uniformly continuous. 
		To get \eqref{non uniform remark equ}, we estimate the error in $H^{2s-\rho_0}$ and $H^{\rho_0}$, respectively, where $\rho_0$ is given in \ref{Assumption-2}. Then  \eqref{non uniform remark equ} is a consequence of the interpolation.  We remark that
		\eqref{same initail data} holds because the approximative solutions are constructed deterministically.

	\end{enumerate}

\end{Remark}

\begin{Remark}\label{remark on instability 2}
	With regard to  similar results in the literature and further  hypotheses, we  give some more  remarks on Theorem \ref{Weak instability}.
	
	\begin{enumerate}[label={ $\bf (\arabic*)$}]
		\item In deterministic cases, the issue of the (optimal) initial-data dependence of solutions has been extensively investigated for various nonlinear dispersive and integrable equations. We refer to \cite{Kato-1975-ARMA} for the inviscid Burgers equation  and to \cite{Koch-Tzvetkov-2005-IMRN} for the Benjamin–Ono equation.
		For the CH equation
		we refer the readers to \cite{Himonas-Kenig-2009-DIE,Himonas-Kenig-Misiolek-2010-CPDE} concerning  the non-uniform dependence on initial data in Sobolev spaces $H^s$ with $s>3/2$. For the first results of this type in
		Besov spaces, we refer to
		\cite{Tang-Zhao-Liu-2014-AA,Tang-Liu-2014-JMP}. Particularly, non-uniform dependence on initial data in critical Besov space first appears in \cite{Tang-Shi-Liu-2015-MM,Tang-Liu-2015-ZAMP}. In this work,  Theorem \ref{Weak instability} and \ref{Local-WP-dependence} in Theorem \ref{Local-WP} demonstrate that the continuity of the solution map $u_{0}\mapsto u$ is almost an optimal result in the sense that, when the growth of the  noise coefficient satisfies certain conditions (cf.\ Hypothesis {\rm\ref{Assumption-2}}), the map $u_{0}\mapsto u$ is continuous, but one can \textit{not} improve the stability of the exiting time and simultaneously the continuity of the map $u_{0}\mapsto u$. Up to our knowledge,
		results of this type for SPDEs  first appeared in \cite{Rohde-Tang-2021-JDDE,Tang-2022-arXiv}. We also refer to  \cite{Alonso-Miao-Tang-2022-JDE,Tang-Yang-2022-AIHP,Miao-Wang-Zhao-CPAA} for recent developments.
		
		\item  It is worthwhile mentioning that, as noted in \textbf{(1)} of Remark \ref{remark on instability 1},  the strong stability of exiting times is used as a technical ``assumption"  to handle the lower bound of a sequence of stopping times. So far we have not been able to verify the non-emptyness of this strong stability  assumption for the current model. However, if the transport noise $u_x\circ \dd W$ is considered ($W$ is a standard 1-D Brownian motion and $\circ \dd W$ means the Stratonovich stochastic differential), we might conjecture that either  the notion of strong stability of exiting times can be captured, or the solution map $u_0\mapsto u$ can become more regular than being continuous. Indeed, if $h(t,u) \dd \W$ is replaced by $u_x\circ \dd W$ in \eqref{SGCH problem}, one can rewrite the equation into It\^o's form with an additional viscous term $-\frac{1}{2}u_{xx}$ on the left hand side of the equation. Therefore, it is reasonable to expect that in this case, either the strong stability of exiting times or the continuity of the solution map $u_0\mapsto u$ can be improved.  We refer to \cite{Himonas-Misiolek-2010-CMP} and \cite{Henry-1981-book} for deterministic examples on the continuity of the solution map.

	\end{enumerate}

\end{Remark}

\begin{Theorem}[Noise prevents blow-up]\label{Non breaking}  Let $s>5/2$, $k\ge1$ and $u_0\in H^s$ be an $\mathcal{F}_0$-measurable random variable with $\E\|u_0\|^2_{H^s}<\infty$.  If 
	Hypothesis {\rm \ref{Assumption-q}} holds true, then
	the corresponding maximal solution
	$(u,\tau^*)$  to \eqref{SGCH non blow up Eq} satisfies
	$$\p\left\{\tau^*=\infty\right\}
	=1.$$
\end{Theorem}

\begin{Remark}\label{remark on global existence}
	
	We notice that many of the existing results on regularization effects by noise are  essentially  restricted to linear equations or linear growing noise. In Theorem \ref{Non breaking},  both the drift and diffusion term are nonlinear. We also remark that the blow-up can actually occur in the  deterministic counterpart of \eqref{SGCH non blow up Eq}. For example, when $k=1$, blow-up (as wave breaking) of solutions to the CH equation can be found in \cite{Constantin-Escher-1998-Acta}.   Therefore, Theorem \ref{Non breaking} demonstrates that large enough noise can  prevent singularities. Indeed, \ref{Ass3-growth} means that the growth of $u^ku_x+F(u)$ can be controlled provided that the noise grows fast enough in terms of 
	a Lyapunov type function $V$. In contrast to \ref{Ass1-Lip s>3/2} and \ref{Ass1-Lip s>1/2}, we require $s>3/2$ in both \ref{Ass3-Lip} and \ref{Ass3-growth}. As is stated in Hypothesis {\rm\ref{Assumption-1}}, \ref{Ass3-Lip} implies that uniqueness holds true for solutions in $H^s$ with $s>5/2$. It seems that one can require $s>1/2$ in \ref{Ass3-Lip} to guarantee   uniqueness in $H^\rho$ with $\rho>3/2$,  but at present we can only construct examples for the case that $s>3/2$ is required in both \ref{Ass3-Lip} and \ref{Ass3-growth}.
	
\end{Remark}

We outline the remainder of the paper. In Sect. \ref{section:cut-off}, we study the cut-off version of \eqref{SGCH problem} and then we remove the cut-off and prove 
Theorem \ref{Local-WP} in Sect. \ref{sect:Local-WP}. We prove Theorem \ref{Weak instability} in Sect. \ref{sect:weak instability}.  Concerning the interplay of noise vs  blow-up, we prove Theorem \ref{Non breaking} in Sect. \ref{sect:noise vs blow-up}.

\section{Cut-off version: Regular solutions}\label{section:cut-off}

We first consider a cut-off version of \eqref{SGCH problem}. To this end,
for any $R>1$, we let $\chi_R(x):[0,\infty)\rightarrow[0,1]$ be a $C^{\infty}$-function such that $\chi_R(x)=1$ for $x\in[0,R]$ and $\chi_R(x)=0$ for $x>2R$.
Then we consider the following cut-off problem
\begin{equation} \label{cut-off problem}
	\left\{\begin{aligned}
		&\, {\rm d}u+\chi_R(\|u\|_{W^{1,\infty}})\left[u^k\partial_xu+F(u)\right]\dd t=\chi_R(\|u\|_{W^{1,\infty}})h(t,u) \dd \W,\\
		&u(\omega,0,x)=u_0(\omega,x)\in H^{s}.
	\end{aligned} \right.
\end{equation}

In this section, we aim at proving the following result:

\begin{Proposition}\label{global solution to cut-off problem}
	Let $s>3$, $k\ge1$, $R>1$ and   Hypothesis {\rm\ref{Assumption-1}} be satisfied. Assume that $u_0\in L^2({\rm \Omega};H^{s})$ is an $H^{s}$-valued $\mathcal{F}_0$-measurable random variable. Then, for  any $T>0$,  \eqref{cut-off problem} has a solution $u \in  L^2\left({\rm \Omega}; C\left([0,T];H^{s}\right)\right)$.
	More precisely, there is a constant $C(R,T,u_0)>0$ such that
	\begin{equation}
		\E\sup_{t\in[0,T]}\|u\|^2_{H^{s}}\leq C(R,T,u_0).\label{solution bound cut-off}
	\end{equation}
\end{Proposition}

The proof for Proposition \ref{global solution to cut-off problem} is given in the following subsections.

\subsection{The approximation scheme}

The first step is to construct a suitable approximation scheme. 
From Lemma \ref{F lemma}, we see that the nonlinear term $F(u)$ preserves the $H^s$-regularity of $u\in H^s$ for any $s>3/2$. However, to apply the theory of SDEs in Hilbert space to \eqref{cut-off problem}, we will have to mollify the transport term $u^k\partial_xu$ since the product $u^k\partial_xu$ loses one regularity. To this end, we consider the following approximation scheme:
\begin{equation} \label{approximate problem}
	\left\{\begin{aligned}
		\, {\rm d}u+H_{1,\e}(u)\dd t&=H_{2}(t,u) \dd \W,\\
		H_{1,\e}(u)&=\chi_R(\|u\|_{W^{1,\infty}})\left[J_{\varepsilon}
		\left((J_{\varepsilon}u)^k\partial_xJ_{\varepsilon}u\right)+F(u)\right],\\
		H_{2}(t,u)&=\chi_R(\|u\|_{W^{1,\infty}})h(t,u),\\
		u(0,x)&=u_0(x)\in H^{s},
	\end{aligned} \right.
\end{equation}
where $J_{\e}$ is the Friedrichs mollifier defined in Appendix  \ref{Section:Preliminaries}. After mollifying the transport term $u^k\partial_xu$, it follows from \ref{Ass1-Lip s>3/2} and Lemmas \ref{mollifier properties} and \ref{F lemma} that for any $\e\in(0,1)$, $H_{1,\e}(\cdot)$ and $H_{2}(t,\cdot)$ are locally Lipschitz continuous in $H^s$ with $s>\frac32$. Besides, we notice that the cut-off function $\chi_R(\|\cdot\|_{\Wlip})$ guarantees the linear growth condition (cf. Lemma \ref{F lemma} and \ref{Ass1-growth assumption}). Thus, for fixed $({\rm \Omega}, \{\mathcal{F}_t\}_{t\geq0},\p, \W)$ and for $u_0\in L^2({\rm \Omega};H^s)$ with $s>3/2$, the existence theory of SDE in Hilbert space  (see for example \cite{Gawarecki-Mandrekar-2010-Springer}) means that  \eqref{approximate problem} admits a unique solution $u_{\e}\in C([0,\infty);H^s)$ $\pas$.

\subsection{Uniform estimates}
Now we establish some uniform-in-$\e$ estimates for $u_\e$.

\begin{Lemma}\label{global solution to appro and estimates} 
	Let $k\geq1$, $s>3/2$, $R>1$ and $\e\in(0,1)$. Assume that $h$ satisfies Hypothesis {\rm\ref{Assumption-1}} and $u_0\in L^2({\rm \Omega};H^s)$ is an $H^s$-valued $\mathcal{F}_0$-measurable random variable. Let $u_{\e}\in C([0,\infty);H^{s})$ be the unique solution to \eqref{approximate problem}. Then  for any $T>0$, 
	there is a constant $C=C(R,T,u_0)>0$ such that
	\begin{align}
		\sup_{\e>0}\E\sup_{t\in[0,T]}\|u_{\varepsilon}(t)\|^2_{H^s}
		\leq C.\label{E r e}
	\end{align}
\end{Lemma}

\proof
Using the It\^{o} formula for $\|u_\e\|^2_{H^s}$, we have that for any $t>0$,
\begin{align*}
	\, {\rm d}\|u_\e(t)\|^2_{H^s}
	=\,  &  2\chi_R\left(\|u_\e\|_{W^{1,\infty}}\right)\left( h(t,u_\e) \dd \W,u_\e \right)_{H^s}\notag\\
	&-2\chi_R\left(\|u_\e\|_{W^{1,\infty}}\right)\left(D^sJ_{\varepsilon}\left[
	(J_{\varepsilon}u_\e)^k\partial_xJ_{\varepsilon}u_\e\right],D^s u_\e \right)_{L^2}\dd t\notag\\
	&-2\chi_R\left(\|u_\e\|_{W^{1,\infty}}\right)\left(D^s F(u_\e),D^s u_\e \right)_{L^2}\dd t\notag\\
	&+ \chi^2_R(\|u_\e\|_{W^{1,\infty}})\|h(t,u_\e)\|_{\LL_2(\U; H^s)}^2\dd t.\notag
\end{align*}
On account of Lemmas  \ref{mollifier properties} and \ref{Kato-Ponce commutator estimate}, we derive
\begin{align*}
	&\,\left|\left(D^sJ_{\varepsilon}\left[
	(J_{\varepsilon}u_\e)^k\partial_xJ_{\varepsilon}u_\e\right],
	D^s u_\e \right)_{L^2}\right|
	\leq   C\|u_\e\|^{k}_{W^{1,\infty}}\|u_\e\|^2_{H^s},
\end{align*}
Therefore, one can infer from the BDG inequality, \ref{Ass1-growth assumption}, Lemma \ref{F lemma} and the above estimate that
\begin{align*}
	\E\sup_{t\in[0,T]}\|u_\e(t)\|^2_{H^s}-\E\|u_0\|^2_{H^s} 
	\leq  \, &    \frac12\E\sup_{t\in[0,T]}\|u_\e\|_{H^s}^2+
	C_R\E \int_0^{T}\left(1+\|u_\e\|^2_{H^s}\right) \dd t,
\end{align*}
which implies
\begin{align}
	\E\sup_{t\in[0,T]}\|u_\e(t)\|^2_{H^s}
	\leq 2\E\|u_0\|^2_{H^s}+ C_R\int_0^{T} \left(1+\E\sup_{t'\in[0,t]}\|u_\e(t')\|_{H^s}^2\right)\dd t.\label{E 2 u-e}
\end{align}
Using Gr\"{o}nwall's inequality in \eqref{E 2 u-e} implies   \eqref{E r e}.
\qed

\subsection{Convergence of approximative solutions}

Now we are going to show that the family $\{u_\e\}$ contains a convergent subsequence.  For different layers $u_\e$ and $u_\eta$, we see that $v_{\e,\eta}:=u_\e-u_\eta$ satisfies the following problem:
\begin{equation}
	{\rm d}v_{\e,\eta}+\Big(\sum_{i=1}^{8}q_i\Big) \dd t=\,\Big(\sum_{i=9}^{10}q_i\Big)  \dd \W,\ \ v_{\e,\eta}(0,x)=0,\label{G12-G12}
\end{equation}
where
\begin{align*}
	q_1:=\,&
	\left[\chi_R(\|u_{\e}\|_{W^{1,\infty}})-\chi_R\left(\|u_{\eta}\|_{W^{1,\infty}}\right)\right]
	J_\e[(J_\e u_\e)^k\partial_xJ_\e u_\e],\\
	q_2:=\,&\chi_R\left(\|u_{\eta}\|_{W^{1,\infty}}\right)
	(J_\e-J_\eta)[(J_\e u_\e)^k\partial_xJ_\e u_\e],\\
	q_3:=\,&\chi_R\left(\|u_{\eta}\|_{W^{1,\infty}}\right)
	J_\eta[((J_\e u_\e)^k-(J_\eta u_\e)^k)\partial_xJ_\e u_\e],\\
	q_4:=\,&\chi_R\left(\|u_{\eta}\|_{W^{1,\infty}}\right)
	J_\eta[((J_\eta u_\e)^k-(J_\eta u_\eta)^k)\partial_xJ_\e u_\e],\\
	q_5:=\,&\chi_R\left(\|u_{\eta}\|_{W^{1,\infty}}\right)
	J_\eta[(J_\eta u_{\eta})^k\partial_x(J_\e-J_\eta) u_\e],\\
	q_6:=\,&\chi_R\left(\|u_{\eta}\|_{W^{1,\infty}}\right)
	J_\eta[(J_\eta u_{\eta})^k\partial_xJ_\eta (u_\e-u_\eta)],\\
	q_7:=\,&
	\left[\chi_R(\|u_{\e}\|_{W^{1,\infty}})-\chi_R\left(\|u_{\eta}\|_{W^{1,\infty}}\right)\right]
	F(u_\e),\\
	q_8:=\,&\chi_R\left(\|u_{\eta}\|_{W^{1,\infty}}\right)[F(u_\e)-F(u_\eta)],\\
	q_9:=\,&\left[\chi_R(\|u_{\e}\|_{W^{1,\infty}})-\chi_R\left(\|u_{\eta}\|_{W^{1,\infty}}\right)\right]
	h(t,u_\e),\\
	q_{10}:=\,&\chi_R\left(\|u_{\eta}\|_{W^{1,\infty}}\right)[h(t,u_\e)-h(t,u_\eta)].
\end{align*}

\begin{Lemma}\label{Qi Lemma}
	Let $s>3$ and $k\ge1$ and let $\mathcal{G}(x):=x^{2k+2}+1$. For any $\e,\eta\in(0,1)$, we find a constant $C>0$ such that
	\begin{align*}
		\sum_{i=1}^{8}\left|(q_i,
		v_{\e,\eta})_{H^{s-\frac32}}\right|
		\leq\,& C\mathcal{G}\left(\|u_\e\|_{H^{s}}+\|u_\eta\|_{H^{s}}\right)
		\left(\|v_{\e,\eta}\|^2_{H^{s-\frac32}}+\max\{\e,\eta\}\right).
	\end{align*}
\end{Lemma}
\proof
Using Lemmas \ref{mollifier properties}, \ref{Kato-Ponce commutator estimate} and \ref{F lemma}, the mean value theorem for $\chi_R(\cdot)$, and the embedding $H^{s-\frac32}\hookrightarrow W^{1,\infty}$, we have that for some $C>0$,
\begin{align*}
	\left\|D^{s-\frac32}q_1\right\|_{L^2},\ 
	\left\|D^{s-\frac32}q_7\right\|_{L^2}
	\leq \, &   C\|v_{\e,\eta}\|_{H^{s-\frac32}}\|u_\e\|^{k+1}_{H^{s}},
\end{align*}
and
\begin{align*}
	\left\|D^{s-\frac32}q_8\right\|_{L^2}
	\leq C\left(\|u_\e\|_{H^{s}}+\|u_\eta\|_{H^{s}}\right)^{k}\|v_{\e,\eta}\|_{H^{s-\frac32}}.
\end{align*}
Using Lemma \ref{mollifier properties}, we see that
\begin{align*}
	\left\|D^{s-\frac32}q_i\right\|_{L^2}\leq \, &  C\max\{\e^{1/2},\eta^{1/2}\}\|u_\e\|^{k+1}_{H^{s}},\ i=2,3,\\
	\left\|D^{s-\frac32}q_4\right\|_{L^2}\leq \, &  C\left(\|u_\e\|_{H^{s}}+\|u_\eta\|_{H^{s}}\right)^{k-1}\|v_{\e,\eta}\|_{H^{s-\frac32}}\|u_\e\|_{H^{s}},\\
	\left\|D^{s-\frac32}q_5\right\|_{L^2}\leq \, &  C\max\{\e^{1/2},\eta^{1/2}\}\|u_\e\|_{H^{s}}\|u_\eta\|^k_{H^{s}}.
\end{align*}
For $q_6$, using Lemma \ref{mollifier properties} and then integrating by part, we have
\begin{align*}
	\left(D^{s-\frac32}q_6,D^{s-\frac32}v_{\e,\eta}\right)_{L^2}
	=\,  & \chi_R\left(\|u_{\eta}\|_{W^{1,\infty}}\right)\int_\R [D^{s-\frac32},(J_\eta u_{\eta})^k]\partial_xJ_\eta v_{\e,\eta}\cdot D^{s-\frac32}J_\eta v_{\e,\eta}\, {\rm d}x\\
	&-\frac12\chi_R\left(\|u_{\eta}\|_{W^{1,\infty}}\right)\int_\R \partial_x (J_\eta u_{\eta})^k(D^{s-\frac32}J_\eta v_{\e,\eta})^2\, {\rm d}x.
\end{align*}
Via the embedding $H^{s-\frac32}\hookrightarrow W^{1,\infty}$ and Lemmas \ref{mollifier properties} and \ref{Kato-Ponce commutator estimate},  we obtain
\begin{align*}
	\left|\left(D^{s-\frac32}q_6,D^{s-\frac32}v_{\e,\eta}\right)_{L^2}\right|
	\lesssim \, &   \|u_\eta\|^k_{H^{s}}\|v_{\e,\eta}\|^{2}_{H^{s-\frac32}}.
\end{align*}
Therefore, we can put this all together to find
\begin{align*}
	&\,\sum_{i=1}^{8}\left|(q_i,
	v_{\e,\eta})_{H^{s-\frac32}}\right|\\
	\leq  \, &     C\left((\|u_\e\|_{H^{s}}+\|u_\eta\|_{H^{s}})^{k+1}+1\right)\|v_{\e,\eta}\|^2_{H^{s-\frac32}}+C(\|u_\e\|_{H^{s}}+\|u_\eta\|_{H^{s}})^{2k+2}\max\{\e,\eta\},
\end{align*}
which gives rise to the desired estimate.
\qed

\begin{Lemma}
	Let $s>3$, $R>1$ and $\e\in(0,1)$.   For any $T>0$ and $K>1$, we define
	\begin{equation} \label{tau-e,eta,T}
		\tau^{T}_{\e,K}:=\inf\left\{t\geq0:\|u_\e(t)\|_{H^{s}}\geq K\right\}\wedge T,\ \ 
		\tau^T_{\e,\eta,K}
		:=\tau^{T}_{\e,K}\wedge\tau^{T}_{\eta,K}.
	\end{equation}
	Then we have
	\begin{equation} \label{Cauchy H s-3/2}
		\lim_{\e\rightarrow0}\sup_{\eta\leq\e}
		\E\sup_{t\in[0,\tau^T_{\e,\eta,K}]}\|u_\e-u_\eta\|_{H^{s-\frac32}}=0.
	\end{equation}
\end{Lemma}
\proof
By employing the BDG inequality to \eqref{G12-G12},  for some constant $C>0$, we arrive at
\begin{align*}
	&\,\E\sup_{t\in[0,\tau^T_{\e,\eta,K}]}\|v_{\e,\eta}(t)\|^{2}_{H^{s-\frac32}}\\
	\leq  \, &    
	\frac12\E\sup_{t\in[0,\tau^T_{\e,\eta,K}]}\|v_{\e,\eta}\|^{2}_{H^{s-\frac32}}
	+C\E \int_0^{\tau^T_{\e,\eta,K}}\sum_{i=1}^{8}\left|(q_i,
	v_{\e,\eta})_{H^{s-\frac32}}\right| \dd t\\
	&+C\E \int_0^{\tau^T_{\e,\eta,K}}\sum_{i=9}^{10}
	\|q_i\|_{\LL_2\big(\U;H^{s-\frac{3}{2}}\big)}^2\dd t.
\end{align*}
For  $q_9$ and $q_{10}$, we 
use \eqref{tau-e,eta,T}, the mean value theorem for $\chi_R(\cdot)$,  \ref{Ass1-growth assumption} and \ref{Ass1-Lip s>3/2} to  find a constant $C=C(K)>0$ such that
\begin{align*}
	\E \int_0^{\tau^T_{\e,\eta,K}}\sum_{i=9}^{10}
	\|q_i\|_{\LL_2\big(\U;H^{s-\frac{3}{2}}\big)}^2\dd t
	\leq\, & C(K)\int_0^T\E\sup_{t'\in[0,\tau^t_{\e,\eta,K}]}\|v_{\e,\eta}(t')\|^{2}_{H^{s-\frac32}}\dd t.
\end{align*}
On account of Lemma \ref{Qi Lemma} and the above estimate, we find
\begin{align*}
	&\,\E\sup_{t\in[0,\tau^T_{\e,\eta,K}]}\|v_{\e,\eta}(t)\|^{2}_{H^{s-\frac32}}\\
	\leq  \, &    
	C(K)\int_0^{T}\E\sup_{t'\in[0,\tau^t_{\e,\eta,K}]}\|v_{\e,\eta}(t')\|^{2}_{H^{s-\frac32}}
	\dd t+C(K)T\max\{\e,\eta\}.
\end{align*}
Therefore, \eqref{Cauchy H s-3/2} holds true. 
\qed

\begin{Lemma}\label{Convergence of ue}
	For any fixed $s>3$ and $T>0$, there is an $\{\mathcal{F}_t\}_{t\geq0}$ progressive measurable $H^{s-3/2}$-valued process $u$ and a countable subsequence of $\{u_\e\}$ $($still denoted as $\{u_\e\})$ such that 
	\begin{equation}\label{convergence ue a.s.}
		u_\e\xrightarrow[]{\e\rightarrow 0}u \  {\rm in}\  C\left([0,T];H^{s-\frac32}\right)\ \ \pas
	\end{equation}
\end{Lemma}

\proof
We first let $\e$ be discrete, i.e., $\e=\e_n (n\ge1)$ such that $\e_n\rightarrow0$ as $n\rightarrow\infty$. In this way, for all $n$, $u_{\e_n}$ can be defined on the same set $\widetilde{{\rm \Omega}}$ with $\p\{\widetilde{{\rm \Omega}}\}=1$. For brevity, $u_{\e_n}$ is still denoted as $u_\e$.  For any $\epsilon>0$, by using  \eqref{tau-e,eta,T}, Lemma \ref{global solution to appro and estimates} and Chebyshev's inequality, we see that
\begin{align*}
	&\, \p\left\{\sup_{t\in[0,T]}\|u_\e-u_\eta\|_{H^{s-\frac32}}>\epsilon\right\} \\
	\leq  \, &     \p\left\{\tau^T_{\e,K}<T\right\}+\p\left\{\tau^T_{\eta,K}<T\right\}
	+\p\left\{
	\sup_{t\in[0,\tau^T_{\e,\eta,K}]}\|u_\e-u_\eta\|_{H^{s-\frac32}}>\epsilon\right\}\\
	\leq \, &  \frac{2C(R,T,u_0)}{K^2}
	+\p\left\{
	\sup_{t\in[0,\tau^T_{\e,\eta,K}]}\|u_\e-u_\eta\|_{H^{s-\frac32}}>
	\epsilon\right\}.
\end{align*}
Now \eqref{Cauchy H s-3/2} clearly forces
$$\lim_{\e\rightarrow0}\sup_{\eta\leq\e}\p\left\{\sup_{t\in[0,T]}\|u_\e-u_\eta\|_{H^{s-\frac32}}>\epsilon\right\}
\leq
\frac{2C(R,T,u_0)}{K^2}.$$
Letting $K\rightarrow\infty$, we see that $u_\e$ converges in probability in $C\left([0,T];H^{s-\frac32}\right)$. 
Therefore, up to a further subsequence, \eqref{convergence ue a.s.} holds. 
\qed

\subsection{Proof for Proposition \ref{global solution to cut-off problem}}
By \eqref{convergence ue a.s.}, since for each $\e\in(0,1)$, $u_\e$ is $\{\mathcal{F}_t\}_{t\geq0}$ progressive measurable, so is $u$. 
Notice that $H^{s-3/2}\hookrightarrow \Wlip$. Then one can send $\e\rightarrow0$ in \eqref{approximate problem} to prove that  $u$ solves \eqref{cut-off problem}.  Furthermore, it follows from Lemma \ref{global solution to appro and estimates} and  Fatou's lemma that
\begin{align}
	\E\sup_{t\in[0,T]}\|u(t)\|^2_{H^s}
	<C(R,u_0,T).\label{u L2 bound}
\end{align} 
With \eqref{u L2 bound}, to prove \eqref{solution bound cut-off}, we  only need to prove $u\in C([0,T];H^{s})$, $\pas$ 
Due to  Lemma \ref{Convergence of ue} and \eqref{L2 moment bound},  $u\in C([0,T];H^{s-3/2})\cap L^\infty\left(0,T;H^s\right) $ almost surely. Since $H^s$ is dense in $H^{s-3/2}$, we see that $u\in C_w\left([0,T];H^s\right)$  ($C_w\left([0,T];H^s\right)$ is the set of weakly continuous functions with values in $H^s$).  Therefore, we only need to  prove the continuity of $[0,T]\ni t\mapsto \|u(t)\|_{H^s}$. As is mentioned in Remark \ref{remark on LW-1}, we first consider the following mollified version with $J_\e$ being defined in \eqref{Friedrichs mollifier}:
\begin{align}
	\, {\rm d}\|J_\e u(t)\|^2_{H^s}
	=\,  & 2\chi_R(\|u\|_{\Wlip}) \left(J_\e h(t,u) \dd \W,J_\e u\right)_{H^s}\nonumber\\
	&-2\chi_R(\|u\|_{\Wlip}) \left(J_{\varepsilon} \left[u^ku_x+F(u)\right], J_\e u\right)_{H^s}\dd t\notag\\
	&+\chi^2_R(\|u\|_{\Wlip})\|J_\e h(t',u)\|_{\LL_2(\U;H^s)}^2\dd t.\label{TnX 2}
\end{align}
By \eqref{u L2 bound}, 
\begin{equation} 
	\tau_N:=\inf\{t\ge0:\|u(t)\|_{H^s}>N\}\rightarrow \infty\ \text{as}\ N\rightarrow\infty \ \ \pas\label{tau N continuous in t}
\end{equation}
Then we only need to prove the continuity up to time $\tau_N\wedge T$ for each $N\ge 1$. 
Let $[t_2,t_1] \subset [0,T]$ with $t_1-t_2<1$. We use Lemma \ref{uux+F u Hs inner product}, the BDG inequality, Hypothesis {\rm\ref{Assumption-1}} and \eqref{tau N continuous in t} to find
\begin{align*}
	\E\left[\left(\|J_\e u(t_1\wedge\tau_N)\|^2_{H^s}-\|J_\e u(t_2\wedge\tau_N)\|^2_{H^s}\right)^4\right]
	\leq  \, &     C(N,T)|t_1-t_2|^{2}.
\end{align*}
We notice that for any $T>0$, $J_\e u$ tends to $u$ in $C\left([0,T];H^{s}\right)$ as $\e\rightarrow 0$. 
This, together with Fatou's lemma,  implies
\begin{align*}
	\E\left[\left(\|u(t_1\wedge\tau_N)\|^2_{H^s}-\| u(t_2\wedge\tau_N)\|^2_{H^s}\right)^4\right]
	\leq  \, &     C(N,T)|t_1-t_2|^{2}.
\end{align*}
This and Kolmogorov's continuity theorem ensure the continuity of $t\mapsto\|u(t\wedge\tau_N)\|_{H^{s}}$.

\section{Proof for Theorem \ref{Local-WP}}\label{sect:Local-WP}

Now we can prove Theorem \ref{Local-WP}. For the sake of clarity, we provide the proof in several subsections.

\subsection{Proof for \ref{Local-WP-existence} in Theorem \ref{Local-WP}: Existence and uniqueness}

\subsubsection{Uniqueness}

Before we prove the existence of a solution in $H^s$ with $s>3/2$, we first prove uniqueness since some estimates here will be used later.
\begin{Lemma}\label{Uniqueness bounded data}
	Let $s>3/2$, $k\ge1$, and  Hypothesis {\rm\ref{Assumption-1}} hold. Suppose that $u_0$ and $v_0$ are two $H^s$-valued $\mathcal{F}_0$-measurable random variables satisfying $ u_0,v_0 \in L^2({\rm \Omega};H^s)$. Let $(u,\tau_1)$ and $(v,\tau_2)$ be two local solutions to \eqref{SGCH problem} in the sense of Definition \ref{Definition of solution} such that $u(0)=u_0$, $v(0)=v_0$ almost surely. 
	For any $N>0$ and $T>0$, we denote
	\begin{align*} 
		\tau_{u}:=\inf\left\{t\geq0: \|u(t)\|_{H^s}>N\right\},\ \
		\tau_{v}:=\inf\left\{t\geq0: \|v(t)\|_{H^s}>N\right\},
	\end{align*}
	and $\tau^T_{u,v}:=\tau_{u}\wedge \tau_{v}\wedge T$.
	Then for $s'\in\left(\frac{1}{2},\min\left\{s-1,\frac{3}{2}\right\}\right)$,    we have that 
	\begin{equation}\label{local time Lip estimate}
		\displaystyle\E\sup_{t\in[0,\tau^T_{u,v}]}\|u(t)-v(t)\|^2_{H^{s'}}\leq C(N,T)\E\|u_0-v_0\|^2_{H^{s'}}.
	\end{equation}
\end{Lemma}
\proof
Let $w(t)=u(t)-v(t)$ for $t\in[0,\tau_1\wedge\tau_2]$. We have
\begin{align*}
	\, {\rm d}w+
	\frac{1}{k+1}\partial_x\left[u^{k+1}-v^{k+1}\right]\dd t+\left[F(u)-F(v)\right]\dd t=
	\left[h(t,u)-h(t,v)\right] \dd \W.
\end{align*}
Then we use the It\^{o} formula for $\|w\|^2_{H^{s'}}$ with $s'\in\left(\frac{1}{2},\min\left\{s-1,\frac{3}{2}\right\}\right)$ to find that
\begin{align*}
	\, {\rm d}\|w\|^2_{H^{s'}}
	=\,  & 2 \left(\left[h(t,u)- h(t,v)\right] \dd \W, w\right)_{H^{s'}}
	-\frac{2}{k+1}\left(\partial_x(P_k w), w\right)_{H^{s'}}\dd t\nonumber\\
	&-2\left(\left[F(u)-F(v)\right],w\right)_{H^{s'}}\dd t
	+ \|h(t,u)-h(t,v)\|_{\LL_2(\U; H^{s'})}^2\dd t\notag\\
	:=\,  &  R_{1}+\sum_{i=2}^4R_i\dd t,
\end{align*}
where $P_k=u^{k}+u^{k-1}v+\cdots+u v^{k-1}+v^k$. 
Taking the supremum over $t\in[0,\tau^T_{u,v}]$ and using the BDG inequality, \ref{Ass1-Lip s>1/2} and the Cauchy--Schwarz inequality yield
\begin{align*}
	&\,\E\sup_{t\in[0,\tau^T_{u,v}]}\|w(t)\|^2_{H^{s'}} -\E\|w(0)\|^2_{H^{s'}} \notag\\
	\leq \, &  \frac12\E\sup_{t\in[0,\tau^T_{u,v}]}\|w\|_{H^{s'}}^2
	+Cg_2(N)\int_0^{T}
	\E\sup_{t'\in[0,\tau^t_{u,v}]}\|w(t')\|^2_{H^{s'}}\dd t
	+\sum_{i=2}^4\E\int_0^{\tau^T_{u,v}}|R_i| \dd t.
\end{align*}
Using Lemma \ref{Taylor}, integration by parts and $H^s\hookrightarrow\Wlip$, we have
\begin{align*}
	|R_2|\lesssim\,&
	\left|\left([D^{s'}\partial_x,P_k]w,D^{s'}w\right)_{L^2}\right|
	+\left|\left(P_kD^{s'}\partial_xw,D^{s'}w\right)_{L^2}\right|\\
	\lesssim\,& \|w\|^2_{H^{s'}}\left(\|u\|_{H^s}+\|v\|_{H^s}\right)^k.
\end{align*}
Therefore, for some constant $C(N)>0$, we have that
\begin{align*}
	\E\int_0^{\tau^T_{u,v}}|R_2| \dd t
	\leq
	C(N)\int_0^{T} 
	\E\sup_{t'\in[0,\tau^t_{u,v}]}\|w(t')\|^2_{H^{s'}}\dd t.
\end{align*}
Similarly, Lemma \ref{F lemma} and \ref{Ass1-Lip s>1/2} yield
\begin{align*}
	\sum_{i=3}^4\E\int_0^{\tau^T_{u,v}}|R_i| \dd t
	\leq
	C(N)\int_0^{T} 
	\E\sup_{t'\in[0,\tau^t_{u,v}]}\|w(t')\|^2_{H^{s'}}\dd t.
\end{align*}
Therefore, we combine the above estimates to  find  
\begin{align*}
	\E\sup_{t\in[0,\tau^T_{u,v}]}\|w(t)\|^2_{H^{s'}}
	\leq 2\E\|w(0)\|^2_{H^{s'}}+
	C(N)\int_0^{T} 
	\E\sup_{t'\in[0,\tau^t_{u,v}]}\|w(t')\|^2_{H^{s'}}
	\dd t.
\end{align*}
Using the Gr\"{o}nwall  inequality in the above estimate leads to \eqref{local time Lip estimate}.
\qed

Similarly, one can obtain  the  following uniqueness result for the original problem \eqref{SGCH problem}, and we omit the details for simplicity.
\begin{Lemma}\label{pathwise uniqueness s>3/2}
	Let $s>3/2$, and let Hypothesis {\rm\ref{Assumption-1}} be true. Let $u_0$ be an $H^s$-valued $\mathcal{F}_0$-measurable random variable such that $u_0\in L^2({\rm \Omega};H^s)$. If $(u_1,\tau_1)$ and $(u_2,\tau_2)$ are two local solutions to \eqref{SGCH problem}   
	satisfying $u_i(\cdot\wedge \tau_i)\in L^2\left({\rm \Omega};C([0,\infty);H^s)\right)$ for $i=1,2$ and
	$\p\{u_1(0)=u_2(0)=u_0(x)\}=1$, then
	\begin{equation*}
		\p\left\{u_1(t,x)=u_2(t,x), \ \ 
		(t,x)\in[0,\tau_1\wedge\tau_2]\times \R\right\}=1.
	\end{equation*}
\end{Lemma}

\subsubsection{The case  $s>3$.}

To begin with, we first state the following existence and uniqueness results in $H^s$ with $s>3$ for the Cauchy problem \eqref{SGCH problem}:

\begin{Proposition}\label{Local regular pathwise solution}
	Let $s>3$, $k\geq1$, and $h(t,u)$ satisfy Hypothesis {\rm\ref{Assumption-1}}. If $u_0$ is an $H^s$-valued $\mathcal{F}_0$-measurable random variable satisfying $\E\|u_0\|^2_{H^s}<\infty$, then there is a unique local  solution $(u,\tau)$ to \eqref{SGCH problem} in the sense of Definition \ref{Definition of solution} with
	\begin{equation}\label{regualr L2 moment bound}
		u(\cdot\wedge \tau)\in L^2\left({\rm \Omega}; C\left([0,\infty);H^s\right)\right).
	\end{equation}
\end{Proposition}

\proof
Since uniqueness has been obtained in Lemma \ref{pathwise uniqueness s>3/2}, 
via Proposition \ref{global solution to cut-off problem}, we only need to remove the cut-off function.
For $u_0(\omega,x)\in L^2({\rm \Omega}; H^s)$, we let  
$$
{\rm \Omega}_m:=\{m-1\leq\|u_0\|_{H^s}<m\},\ m\ge 1.
$$
Let $u_{0,m}:=u_0\textbf{1}_{\{m-1\leq\|u_0\|_{H^s}<m\}}$.
For any $R>0$,
on account of Proposition \ref{global solution to cut-off problem}, we let $u_{m,R}$ be the global solution to the cut-off problem \eqref{cut-off problem} with initial value $u_{0,m}$ and cut-off function $\chi_R(\cdot)$. Define
\begin{equation*}
	\tau_{m,R}:=\inf\left\{t>0:\sup_{t'\in[0,t]}\|u_{m,R}(t')\|^2_{H^s}>\|u_{0,m}\|^2_{H^s}+2\right\}.
\end{equation*}
Then for any $R>0$ and $m\in\N$, it follows from the time continuity  of the solution that $\p\{\tau_{m,R}>0\}=1$. Particularly, for any $m\in\N$, we assign $R=R_m$  such that $R^2_m>c^2m^2+2c^2$, where $c>0$ is the embedding constant such that $\|\cdot\|_{\Wlip} \leq c \|\cdot\|_{H^s}$ for $s>3$. For simplicity, we denote $(u_m,\tau_m):=(u_{m,R_m},\tau_{m,R_m})$.  Then we have
$$\p\left\{\|u_m\|^2_{W^{1,\infty}}\leq c^2\|u_m\|^2_{H^{s}}\leq c^2\|u_{0,m}\|^2_{H^{s}}+2c^2<R^2_m,\ t\in[0,\tau_{m}],\  m\ge1\right\}=1,$$
which means $\p\left\{\chi_{R_m}(\|u_m\|_{W^{1,\infty}})=1,\  t\in[0,\tau_m],\    m\ge1\right\}=1.$
Therefore, $(u_m,\tau_m)$ is the solution  to \eqref{SGCH problem} with initial value $u_{0,m}$. 
Since $\E\|u_0\|^2_{H^s}<\infty$, the condition \eqref{Omega i condition} is satisfied with $I=\N^+$.
Applying Lemma \ref{cut-combine} means that
\begin{equation*}
	\left(u=\sum_{m\geq1}\textbf{1}_{\{m-1\leq\|u_0\|_{H^s}<m\}}u_m,\ \
	\tau=\sum_{m\geq1}\textbf{1}_{\{m-1\leq\|u_0\|_{H^s}<m\}}\tau_m\right)
\end{equation*}
is a solution to \eqref{SGCH problem} corresponding to the initial condition $u_0$. Besides, 
\begin{align*}
	\sup_{t\in[0,\tau]}\|u\|_{H^s}^2
	=\, & \sum_{m\geq1}\textbf{1}_{\{m-1\leq\|u_0\|_{H^s}<m\}}\sup_{t\in[0,\tau_m]}\|u_m\|_{H^s}^2
	\leq 2\|u_{0}\|^2_{H^s}+4\ \ \pas
\end{align*}
Taking expectation gives rise to \eqref{regualr L2 moment bound}. 
\qed

\subsubsection{The case  $s>3/2$}

When $s>3/2$, we first consider the  following problem
\begin{equation} \label{smooth initial problem}
	\left\{\begin{aligned}
		&\, {\rm d}u+\left[u^k\partial_xu+F(u)\right]\dd t=h(t,u) \dd \W,\ \ k\ge1,\ \ x\in\R,\ t>0,\\
		&u(\omega,0,x)=J_{\e} u_0(\omega,x)\in H^{\infty},\ \ x\in\R,\ \e\in(0,1),
	\end{aligned} \right.
\end{equation}
where $J_\e$ is the mollifier defined in \eqref{Friedrichs mollifier}. Proposition \ref{Local regular pathwise solution} implies that for each $\e\in(0,1)$, \eqref{smooth initial problem} has a local pathwise solution $(u_\e,\tau_\e)$ such that $u_\e \in  L^2\left({\rm \Omega}; C\left([0,\tau_\e];H^{s}\right)\right)$.

\begin{Lemma}\label{Cauchy estimates lemma}
	Assume  $u_0$ is an $H^s$-valued $\mathcal{F}_0$-measurable random variable such that $\|u_0\|_{H^s}\leq M$ for some $M>0$.
	For any $T>0$ and $s>3/2$, we define
	\begin{equation} \label{tau-e eta T s>3/2}
		\tau^T_{\e}:=\inf\left\{t\geq0:\|u_\e\|_{H^s}\geq \|J_{\e} u_0\|_{H^s}+2\right\}\wedge T,
		\ \ 
		\tau^{T}_{\e,\eta}:=\tau^T_{\e}\wedge \tau^T_{\eta},\ \ \e,\eta\in(0,1).
	\end{equation}
	Let $K\ge 2M+5$ be fixed and let $s'\in\left(\frac{1}{2},\min\left\{s-1,\frac{3}{2}\right\}\right)$.
	Then, there is a constant $C(K,T)>0$ such that $w_{\e,\eta}=u_\e-u_\eta$ satisfies
	\begin{align}
		&\,\E\sup_{t\in[0,\tau^{T}_{\e,\eta}]}\|w_{\e,\eta}(t)\|^2_{H^s}\notag\\
		\leq \, &      C(K,T)\E\left\{\|w_{\e,\eta}(0)\|^2_{H^s}+\|w_{\e,\eta}(0)\|^2_{H^{s'}}\|u_\e(0)\|^2_{H^{s+1}}\right\}\notag\\
		&+C(K,T)\E\sup_{t\in[0,\tau^{T}_{\e,\eta}]}\|w_{\e,\eta}(t)\|^2_{H^{s'}}.\label{Cauchy lemma step 1}
	\end{align}
\end{Lemma}

\proof
To start with, we notice that Lemma \ref{mollifier properties} implies
\begin{equation}
	\|J_{\e} u_0\|_{H^s}\leq M,\ \ \e\in(0,1)\ \ \pas\label{bound mollify initial}
\end{equation}
Since \eqref{tau-e eta T s>3/2} and \eqref{bound mollify initial} are used frequently in the following, they will be used without further notice. 
Let \begin{equation}\label{q k n m}
	P_{l}=P_{l}(u_\e,u_\eta):=
	\begin{cases}
		u_{\e}^{l}+u_{\e}^{l-1}u_{\eta}+\cdots+u_{\e}u_{\eta}^{l-1}+u_{\eta}^{l},\ \ \text{if}\ \ l\ge1,\\
		1,\ \ \text{if}\ \ l=0.
	\end{cases}
\end{equation}
Applying the It\^{o} formula to $\|w_{\e,\eta}\|^2_{H^s}$ gives rise to
\begin{align}
	\, {\rm d}\|w_{\e,\eta}\|^2_{H^s}
	=\,  & 2
	\left(\left[h(t,u_\e)-h(t,u_\eta)\right] \dd \W, w_{\e,\eta}\right)_{H^s} 
	-2\left(w_{\e,\eta} P_{k-1}\partial_x u_{\e},w_{\e,\eta}\right)_{H^s}\dd t \notag\\
	&-2\left(u^k_{\eta}\partial_x w_{\e,\eta},w_{\e,\eta}\right)_{H^s}\dd t -2\left(\left[F(u_\e)-F(u_\eta)\right],w_{\e,\eta}\right)_{H^s}\dd t\notag\\
	&+\left\|h(t,u_\e)-h(t,u_\eta)\right\|_{\LL_2(\U;H^s)}^2\dd t  
	:=\,  Q_{1,s}+\sum_{i=2}^{5}Q_{i,s}\dd t.\label{w Qis equation}
\end{align}
Since $H^{s'}\hookrightarrow L^\infty$ and $H^{s}\hookrightarrow \Wlip$, we can use Lemmas \ref{Kato-Ponce commutator estimate} and \ref{F lemma} to find
\begin{align*}
	\left|Q_{2,s}\right|
	\lesssim \, &   \left(\|w_{\e,\eta}P_{k-1}\|_{H^s}\|\partial_x u_\e\|_{L^\infty}
	+\|w_{\e,\eta}P_{k-1}\|_{L^\infty}\|\partial_x u_\e\|_{H^s}\right)
	\|w_{\e,\eta}\|_{H^s}\\
	\lesssim \, &   \|w_{\e,\eta}\|^2_{H^s}\left(\|u_\e\|_{H^{s}}+\|u_\eta\|_{H^{s}}\right)^{k}
	+\|w_{\e,\eta}\|^2_{H^{s'}}
	\| u_\e\|^2_{H^{s+1}}\\
	&+\left(\|u_\e\|_{H^{s}}+\|u_\eta\|_{H^{s}}\right)^{2k-2}
	\|w_{\e,\eta}\|^2_{H^s}
\end{align*}
\begin{align*}
	\left|Q_{3,s}\right|\lesssim \, &   
	\left|\left([D^s,u^k_\eta]\partial_x w_{\e,\eta},D^s w_{\e,\eta}\right)_{L^2}\right|
	+\left|\left(u^k_\eta\partial_x D^s w_{\e,\eta},D^s w_{\e,\eta}\right)_{L^2}\right|\\
	\lesssim \,&\|w_{\e,\eta}\|^2_{H^s}\|u_\eta\|_{H^{s}}^k,
\end{align*}
and
\begin{align*}
	\left|Q_{4,s}\right|
	\lesssim \, &  \|w_{\e,\eta}\|^2_{H^s}\left(\|u_\e\|_{H^{s}}+\|u_\eta\|_{H^{s}}\right)^{k}.
\end{align*}
The above estimates and \ref{Ass1-Lip s>3/2} imply that there is a constant $C(K)>0$ such that
\begin{align*}
	&\,\sum_{i=2}^5\E\int_0^{\tau^{T}_{\e,\eta}}|Q_{i,s}| \dd t \\
	\lesssim \, &   \E\int_0^{\tau^{T}_{\e,\eta}}
	\left[\left(\left(\|u_\e\|_{H^s}+\|u_\eta\|_{H^s}\right)^{2k}+1\right)\|w_{\e,\eta}\|^2_{H^s}
	+\|w_{\e,\eta}\|^2_{H^{s'}}\|u_\e\|^2_{H^{s+1}}\right]\dd t\notag\\
	&+\E\int_0^{\tau^{T}_{\e,\eta}}
	g_1^2(K)\|w_{\e,\eta}\|^2_{H^s} \dd t\notag\\
	\leq \, &  C(K)\int_0^{T}
	\E\sup_{t'\in[0,\tau^{t}_{\e,\eta}]}\|w_{\e,\eta}(t')\|^2_{H^s}\dd t
	+C(K)T\E\sup_{t\in[0,\tau^{T}_{\e,\eta}]}\|w_{\e,\eta}(t)\|^2_{H^{s'}}\|u_\e(t)\|^2_{H^{s+1}}.
\end{align*}
For $Q_{1,s}$, applying the BDG inequality and \ref{Ass1-Lip s>3/2}, we derive
\begin{align*}
	&\, \E\left(\sup_{t\in[0,\tau^{T}_{\e,\eta}]}\int_0^{t}\Big(\left[h(t,u_\e)-h(t,u_\eta)\right] \dd \W, w_{\e,\eta}\Big)_{H^s}\right)\\
	\leq  \, &     \frac12\E\sup_{t\in[0,\tau^{T}_{\e,\eta}]}\|w_{\e,\eta}\|_{H^s}^2
	+Cg_1^2(K)\int_0^{T}
	\E\sup_{t'\in[0,\tau^{t}_{\e,\eta}]}\|w_{\e,\eta}(t')\|^2_{H^s}\dd t.
\end{align*}
Summarizing the above estimates and then using Gr\"{o}nwall's inequality,  we find some constant $C=C(K,T)>0$ such that
\begin{align}
	\,&\E\sup_{t\in[0,\tau^{T}_{\e,\eta}]}\|w_{\e,\eta}(t)\|^2_{H^s} \notag\\
	\leq  \,   
	&C\left(\E\|w_{\e,\eta}(0)\|^2_{H^s}
	+\E\sup_{t\in[0,\tau^{T}_{\e,\eta}]}\|w_{\e,\eta}(t)\|^2_{H^{s'}}\|u_\e(t)\|^2_{H^{s+1}}\right).
	\label{Verify Cauchy 1}
\end{align}
Now we estimate $\E\sup_{t\in[0,\tau^{T}_{\e,\eta}]}\|w_{\e,\eta}(t)\|^2_{H^{s'}}\|u_\e(t)\|^2_{H^{s+1}}$.  To this end, we first recall \eqref{define cylindrical process} and then apply the It\^{o} formula to deduce that for any $\rho>0$,
\begin{align}
	\, {\rm d}\|u_\e\|^2_{H^{\rho}}
	=\,  & 2\sum_{l=1}^{\infty}\left(h(t,u_\e)e_l,u_\e\right)_{H^{\rho}} \dd W_l
	-2\left(D^\rho\left[(u_\e)^k\partial_x u_\e\right],D^\rho u_\e\right)_{L^2}\dd t\notag\\
	&-2\left(D^\rho F(u_\e),D^\rho u_\e\right)_{L^2}\dd t
	+\|h(t,u_\e)\|_{\LL_2(\U;H^\rho)}^2\dd t\nonumber\\
	:=\, & \sum_{l=1}^{\infty}Z_{1,\rho,l} \dd W_l+\sum_{i=2}^4Z_{i,\rho}\dd t.\label{u Zis equation}
\end{align}
In the same way,  we also rewrite $Q_{1,s}$ in \eqref{w Qis equation} as
\begin{equation}
	Q_{1,s}=2\sum_{j=1}^{\infty}
	\left(\left[h(t,u_\e)- h(t,u_\eta)\right]e_j,w_{\e,\eta}\right)_{H^s} \dd W_j:=\sum_{j=1}^{\infty}Q_{1,s,j} \dd W_j.\label{Q1s sum form}
\end{equation}
With the summation form \eqref{Q1s sum form} at hand,  applying the It\^{o}  product rule to \eqref{w Qis equation}  and \eqref{u Zis equation}, we derive
\begin{align*}
	\, {\rm d}\|w_{\e,\eta}\|^2_{H^{s'}}\|u_\e\|^2_{H^{s+1}}
	=\, & \sum_{j=1}^{\infty}\left(\|w_{\e,\eta}\|^2_{H^{s'}}Z_{1,s+1,j}
	+\|u_\e\|^2_{H^{s+1}}Q_{1,s',j}\right) \dd W_j\nonumber\\
	&+\sum_{i=2}^4\|w_{\e,\eta}\|^2_{H^{s'}}Z_{i,s+1}\dd t
	+\sum_{i=2}^5\|u_\e\|^2_{H^{s+1}}Q_{i,s'}\dd t \\
	&+\sum_{j=1}^{\infty}Q_{1,s',j}Z_{1,s+1,j}\dd t. 
\end{align*}
We first notice that
$$Q_{2,s'}+Q_{3,s'}=\frac{2}{k+1}\left(\partial_x (P_{k}w_{\e,\eta}),w_{\e,\eta}\right)_{H^{s'}},$$
where $P_{k}$ is defined by \eqref{q k n m}. As a result,
Lemma \ref{Taylor}, integration by parts and $H^s\hookrightarrow\Wlip$ give rise to
\begin{align*}
	|Q_{2,s'}+Q_{3,s'}|\lesssim \, &  
	\|w_{\e,\eta}\|^2_{H^{s'}}\left(\|u_\e\|_{H^s}+\|u_\eta\|_{H^s}\right)^k.
\end{align*}
Using Lemma \ref{Kato-Ponce commutator estimate}, Hypothesis {\rm\ref{Assumption-1}}, Lemma \ref{F lemma} as well as the embedding of $H^s\hookrightarrow W^{1,\infty}$ for $s>3/2$, we obtain that for some $C(K)>0$,
\begin{align*}
	\sum_{i=2}^4\|w_{\e,\eta}\|^2_{H^{s'}}|Z_{i,s+1}|
	\lesssim  \|w_{\e,\eta}\|^2_{H^{s'}}
	\left[
	\|u_\e\|^k_{H^s}\|u_\e\|^2_{H^{s+1}}+f^2(\|u_\e\|_{H^s})(1+\|u_\e\|^2_{H^{s+1}})
	\right],
\end{align*}
\begin{align*}
	\sum_{i=4}^5\E\int_0^{\tau^{T}_{\e,\eta}}\|u_\e\|^2_{H^{s+1}}|Q_{i,s'}| \dd t 
	\leq \,    C(K)\int_0^{T}
	\E\sup_{t'\in[0,\tau^{t}_{\e,\eta}]}\|u_\e(t')\|^2_{H^{s+1}}\|w_{\e,\eta}(t')\|^2_{H^{s'}}
	\dd t.
\end{align*}
Then one can infer from
the above three inequalities, the BDG inequality and Hypothesis {\rm\ref{Assumption-1}}  that for some constant $C(K)>0$,
\begin{align}
	&\,\E\sup_{t\in[0,\tau^{T}_{\e,\eta}]}\|w_{\e,\eta}\|^2_{H^{s'}}\|u_\e\|^2_{H^{s+1}}
	-\E\|w_{\e,\eta}(0)\|^2_{H^{s'}}\|u_\e(0)\|^2_{H^{s+1}}\notag\\
	\lesssim \, &   \E\left(\int_0^{\tau^{T}_{\e,\eta}}
	\|w_{\e,\eta}\|^4_{H^{s'}}\|h(t,u_\e)\|^2_{\LL_2(\U; H^{s+1})}\|u_\e\|^2_{H^{s+1}}\dd t\right)^{\frac12}\notag\\
	&+\E\left(\int_0^{\tau^{T}_{\e,\eta}}
	\|u_\e\|^4_{H^{s+1}}\|h(t,u_\e)-h(t,u_\eta)\|^2_{\LL_2(\U; H^{s'})}\|w_{\e,\eta}\|^2_{H^{s'}}
	\dd t\right)^{\frac12}\notag\\
	&+\sum_{i=2}^4\E\int_0^{\tau^{T}_{\e,\eta}}\|w_{\e,\eta}\|^2_{H^{s'}}|Z_{i,s+1}| \dd t
	+\E\int_0^{\tau^{T}_{\e,\eta}}\|u_\e\|^2_{H^{s+1}}|Q_{2,s'}+Q_{3,s'}| \dd t\notag\\
	&+\sum_{i=4}^5\E\int_0^{\tau^{T}_{\e,\eta}}\|u_\e\|^2_{H^{s+1}}|Q_{i,s'}| \dd t
	+\E\int_0^{\tau^{T}_{\e,\eta}}\sum_{j=1}^{\infty}|Q_{1,s',j}Z_{1,s+1,j}| \dd t\notag\\
	\leq \, &   
	\frac12\E\sup_{t\in[0,\tau^{T}_{\e,\eta}]}
	\|w_{\e,\eta}\|^2_{H^{s'}}\|u_\e\|^2_{H^{s+1}}\notag\\
	&+C(K)\int_0^{T}
	\E\sup_{t'\in[0,\tau^{t}_{\e,\eta}]}\|w_{\e,\eta}(t')\|^2_{H^{s'}}\|u_\e(t')\|^2_{H^{s+1}}
	\dd t\notag\\
	&+C(K)T
	\E\sup_{t\in[0,\tau^{T}_{\e,\eta}]}\|w_{\e,\eta}(t)\|^2_{H^{s'}}
	+\E\int_0^{\tau^{T}_{\e,\eta}}\sum_{j=1}^{\infty}|Q_{1,s',j}Z_{1,s+1,j}| \dd t.
	\label{w u product estimate}
\end{align}
For the last term, we proceed as follows:
\begin{align*}
	&\,\E\int_0^{\tau^{T}_{\e,\eta}}\sum_{j=1}^{\infty}\left|Q_{1,s',j}Z_{1,s+1,j}\right| \dd t\\
	\lesssim \, &   \E\int_0^{\tau^{T}_{\e,\eta}}\|h(t,u_\e)-h(t,u_\eta)\|_{\LL_2(\U;H^{s'})}\|w_{\e,\eta}\|_{H^{s'}}
	\|h(t,u_\e)\|_{\LL_2(\U;H^{s+1})}\|u_\e\|_{H^{s+1}}\dd t.\notag\\
	\leq \, &   C(K)T
	\E\sup_{t\in[0,\tau^{T}_{\e,\eta}]}\|w_{\e,\eta}(t)\|^2_{H^{s'}}\\
	&+C(K)\int_0^{T}
	\E\sup_{t'\in[0,\tau^{t}_{\e,\eta}]}\|w_{\e,\eta}(t')
	\|^2_{H^{s'}}\|u_\e(t')\|^2_{H^{s+1}}\dd t.
\end{align*}
Consequently, \eqref{w u product estimate} reduces to
\begin{align*}
	&\E\sup_{t\in[0,\tau^{T}_{\e,\eta}]}\|w_{\e,\eta}\|^2_{H^{s'}}\|u_\e\|^2_{H^{s+1}}-2\E\|w_{\e,\eta}(0)\|^2_{H^{s'}}\|u_\e(0)\|^2_{H^{s+1}}\\
	\leq \, &   C(K)T
	\E\sup_{t\in[0,\tau^{T}_{\e,\eta}]}\|w_{\e,\eta}(t)\|^2_{H^{s'}}\\
	&+C(K)\int_0^{T}
	\E\sup_{t'\in[0,\tau^{t}_{\e,\eta}]}\|w_{\e,\eta}(t')
	\|^2_{H^{s'}}\|u_\e(t')\|^2_{H^{s+1}}\dd t,
\end{align*}
which means that for some $C(K,T)>0$,
\begin{align}
	&\,\E\sup_{t\in[0,\tau^{T}_{\e,\eta}]}\|w_{\e,\eta}\|^2_{H^{s'}}\|u_\e\|^2_{H^{s+1}}\notag\\
	\leq \, &   C\left(\E\|w_{\e,\eta}(0)\|^2_{H^{s'}}\|u_\e(0)\|^2_{H^{s+1}}
	+\E\sup_{t\in[0,\tau^{T}_{\e,\eta}]}\|w_{\e,\eta}(t)\|^2_{H^{s'}}\right).
	\label{to check Cauchy 2 2 terms}
\end{align}
Combining \eqref{Verify Cauchy 1} and \eqref{to check Cauchy 2 2 terms}, we obtain \eqref{Cauchy lemma step 1}.
\qed

To proceed further, we state the following lemma in \cite{GlattHoltz-Ziane-2009-ADE} as a form which is
convenient for our purposes.
\begin{Lemma}[Lemma 5.1, \cite{GlattHoltz-Ziane-2009-ADE}]\label{Cauchy lemma GZ}
	Let all the conditions in Lemma \ref{Cauchy estimates lemma} hold true. 
	Assume
	\begin{equation}
		\lim_{\e\rightarrow0}\sup_{\eta\leq \e}
		\E\sup_{t\in[0,\tau^T_{\e,\eta}]}\|u_\e-u_\eta\|_{H^s}=0\label{Cauchy 1}
	\end{equation}
	and
	\begin{equation} \label{Cauchy 2}
		\lim_{a\rightarrow0}\sup_{\e\in(0,1)}
		\p\left\{\sup_{t\in[0,\tau^T_{\e}\wedge a]}\|u_\e\|_{H^s}\geq \|J_{\e}u_0\|_{H^s}+1\right\}=0
	\end{equation}
	hold true. Then we have:
	
	\begin{enumerate} [label={\bf (\alph*)}]
		\item \label{GZ lemma-time} There exits 
		a sequence of stopping times $\xi_{\e_n}$, for some countable sequence $\{\e_n\}$ with $\e_n\rightarrow 0$ as $n\rightarrow\infty$, and a stopping time $\tau$ such that
		\begin{equation*}
			\xi_{\e_n}\leq \tau^T_{\e_n},\ \ 
			\lim_{n\rightarrow\infty}\xi_{\e_n}=\tau\in(0,T] \ \ \pas
		\end{equation*}

		\item \label{GZ lemma-convergence} There is a process $u\in C([0,\tau];H^s)$ such that
		\begin{equation*}
			\lim_{n\rightarrow\infty}\sup_{t\in[0,\tau]}\|u_{\e_n}-u\|_{H^s}=0,\ 
			\sup_{t\in[0,\tau]}\|u\|_{H^s}\leq \|u_0\|_{H^s}+2\ \ \pas
		\end{equation*}

		\item\label{GZ lemma-moment} There is a sequence of sets ${\rm \Omega}_n \uparrow {\rm \Omega}$ such that for any $p\in[1,\infty)$,
		\begin{equation*}
			{\bf 1}_{{\rm \Omega}_n }\sup_{t\in[0,\tau]}\|u_{\e_n}\|_{H^s}\leq \|u_0\|_{H^s}+2\ \ \pas,
			\ \text{and}\ \sup_{n}\E \left({\bf 1}_{{\rm \Omega}_n }
			\sup_{t\in[0,\tau]}\|u_{\e_n}\|^p_{H^s}\right)<\infty.
		\end{equation*}

	\end{enumerate}
	
\end{Lemma}

%\begin{center}
%{there exi\textbf{s}ts in (a), ``exists", check the whole paper}
%\end{center}

\begin{Remark}
	In the original form of \cite[Lemma 5.1]{GlattHoltz-Ziane-2009-ADE}, the authors only emphasize the existence of stopping time $\tau\in(0,T]$ such that \ref{GZ lemma-convergence} and \ref{GZ lemma-moment} in Lemma \ref{Cauchy lemma GZ} hold true. However, here we point out that they obtained such $\tau$ by constructing 
	stopping times $\xi_{\e_n}$. We refer to  (5.2), (5.12), (5.15), (5.20) and (5.24) in \cite{GlattHoltz-Ziane-2009-ADE} for the details. The properties \ref{GZ lemma-time} and \ref{GZ lemma-moment} in  Lemma \ref{Cauchy lemma GZ}  will be used in the proof for \ref{Local-WP-dependence} in Theorem \ref{Local-WP}.
\end{Remark}

\begin{Proposition}\label{pathwise solution s>3/2 bounded initial data}
	Let Hypothesis {\rm\ref{Assumption-1}} hold. Assume that $s>3/2$, $k\ge1$ and let $u_0$ is an $H^s$-valued $\mathcal{F}_0$-measurable random variable such that $\|u_0\|_{H^s}\leq M$ for some $M>0$. Then \eqref{SGCH problem} has a unique pathwise solution $(u,\tau)$ in the sense of Definition \ref{Definition of solution} such that
	\begin{equation*}
		\sup_{t\in[0,\tau]}\|u\|_{H^s}\leq \|u_0\|_{H^s}+2\ \ \pas
	\end{equation*}
\end{Proposition}
\proof
We  first prove that $\{u_\e\}$ satisfies the estimates \eqref{Cauchy 1} and \eqref{Cauchy 2}. 

\smallskip

\textbf{(i) \eqref{Cauchy 1} is satisfied.}
Lemma \ref{mollifier properties} tells us that
\begin{align}
	\lim_{\e\rightarrow0}\sup_{\eta\leq \e}\E\|w_{\e,\eta}(0)\|^2_{H^s}=0.
	\label{Claim 1}
\end{align}
Since $\|J_\e u_0\|_{H^s}\leq M$,   as in Lemma \ref{Uniqueness bounded data}, we have
\begin{align}
	\lim_{\e\rightarrow0}\sup_{\eta\leq \e}
	\E\sup_{t\in[0,\tau^T_{\e,\eta}]}\|w_{\e,\eta}(t)\|^2_{H^{s'}}
	\leq C(M,T)\lim_{\e\rightarrow0}\sup_{\eta\leq \e}
	\E\|w_{\e,\eta}(0)\|^2_{H^{s}}=0.\label{Claim 2}
\end{align}
Moreover, it follows from Lemma \ref{mollifier properties} that
\begin{align}
	\lim_{\e\rightarrow0}\sup_{\eta\leq \e}\|w_{\e,\eta}(0)\|^2_{H^{s'}}\|u_\e(0)\|^2_{H^{s+1}}
	\lesssim \lim_{\e\rightarrow0}\sup_{\eta\leq \e}o\left(\e^{2s-2s'}\right)O\left(\frac{1}{\e^2}\right)=0.\label{Claim 3}
\end{align}
Summarizing \eqref{Claim 1}, \eqref{Claim 2}, \eqref{Claim 3} and Lemma \ref{Cauchy estimates lemma},  \eqref{Cauchy 1} holds true.

\smallskip

\textbf{(ii) \eqref{Cauchy 2} is satisfied.}

Recall \eqref{u Zis equation} and let $a>0$. We have
\begin{align*}
	\sup_{t\in[0,\tau_\e^T\wedge a]}\|u_{\e}(t)\|^2_{H^s}
	\leq  \, &    \|J_{\e}u_0\|^2_{H^s}
	+\sup_{t\in[0,\tau_\e^T\wedge a]}\left|\int_0^{t}\sum_{j=1}^{\infty}Z_{1,s,j} \dd W_j\right|
	+\sum_{i=2}^4\int_0^{\tau_\e^T\wedge a}|Z_{i,s}| \dd t,
\end{align*}
which clearly forces that
\begin{align*}
	&\p\left\{\sup_{t\in[0,\tau_\e^T\wedge a]}\|u_{\e}(t)\|^2_{H^s}>\|J_{\e}u_0\|^2_{H^s}+1\right\}
	\notag\\
	\leq \, &  
	\p\left\{\sup_{t\in[0,\tau_\e^T\wedge a]}\left|\int_0^{t}\sum_{j=1}^{\infty}Z_{1,s,j} \dd W_j\right|>\frac12\right\}
	+\p\left\{\sum_{i=2}^4\int_0^{\tau_\e^T\wedge a}|Z_{i,s}| \dd t>\frac12\right\}.
\end{align*}
Due to the Chebyshev inequality, Lemmas \ref{Kato-Ponce commutator estimate} and \ref{F lemma}, Hypothesis {\rm\ref{Assumption-1}}, the embedding of $H^s\hookrightarrow W^{1,\infty}$ for $s>3/2$, \eqref{tau-e eta T s>3/2} and \eqref{bound mollify initial}, we have
\begin{align*}
	\p\left\{\sum_{i=2}^4\int_0^{\tau_\e^T\wedge a}|Z_{i,s}| \dd t>\frac12\right\}
	\leq \, &   C\sum_{i=2}^4\E\int_0^{\tau_\e^T\wedge a}|Z_{i,s}| \dd t\notag\\
	\leq \, &   C\E\int_0^{\tau_\e^T\wedge a}
	\left[ \|u_\e\|^{k+2}_{H^{s}}+f^2(\|u_\e\|_{H^s})(1+\|u_\e\|^2_{H^{s}})\right]
	\dd t\notag\\
	\leq \, &   C\E\int_0^{\tau_\e^T\wedge a}C(M,T)\dd t\leq C(M,T)a.
\end{align*}
Then we can infer from the Doob's maximal inequality and the It\^{o} isometry that
\begin{align*}
	\,&\p\left\{\sup_{t\in[0,\tau_\e^T\wedge a]}\left|\int_0^{t}
	\sum_{j=1}^{\infty}Z_{1,s,j} \dd W_j\right|>\frac12\right\} \\
	\leq \, &   C\E\left(\int_0^{\tau_\e^T\wedge a}\sum_{j=1}^{\infty}Z_{1,s,j} \dd W_j\right)^2\notag\\
	\leq \, &   C\E\int_0^{\tau_\e^T\wedge a}
	\left[f^2(\|u_\e\|_{W^{1,\infty}})(1+\|u_\e\|_{H^{s}})^2\|u_\e\|^2_{H^{s}}\right]
	\dd t\notag\\
	\leq \, &   C\E\int_0^{\tau_\e^T\wedge a}C(M,T)\dd t\leq C(M,T)a.
\end{align*}
Hence we have
\begin{align*}
	\p\left\{\sup_{t\in[0,\tau_\e^T\wedge a]}\|u_{\e}(t)\|^2_{H^s}>\|J_{\e}u_0\|^2_{H^s}+1\right\}
	\leq
	C(M,T)a,
\end{align*}
which gives \eqref{Cauchy 2}.  

\textbf{(iii) Applying Lemma \ref{Cauchy lemma GZ}.} By
Lemma \ref{Cauchy lemma GZ}, we can take limit in some subsequence of $\{u_{\e_n}\}$ to build a solution $u$ to \eqref{SGCH problem} such that $u\in C([0,\tau];H^s)$ and $\sup_{t\in[0,\tau]}\|u\|_{H^s}\leq \|u_0\|_{H^s}+2.$ Uniqueness is a direct 
corollary
of Lemma \ref{pathwise uniqueness s>3/2}.
\qed

Now we can finish the proof for \ref{Local-WP-existence} in Theorem \ref{Local-WP}.

\textit{Proof for \ref{Local-WP-existence} in Theorem \ref{Local-WP}.}
As in Proposition
\ref{Local regular pathwise solution}, we let
\begin{align*}
	u_0(\omega,x)&:=\sum_{m\geq1}u_{0,m}(\omega,x)
	:=\sum_{m\geq1}u_0(\omega,x)\textbf{1}_{\{m-1\leq\|u_0\|_{H^s}<m\}}\ \  \pas
\end{align*}
For each $m\ge1$, we can infer from Proposition \ref{pathwise solution s>3/2 bounded initial data} that \eqref{SGCH problem} has a unique solution $(u_m,\tau_m)$ with $u_m(0)=u_{0,m}$ almost surely. Furthermore,
$\sup_{t\in[0,\tau_m]}\|u_m\|_{H^s}\leq \|u_{0,m}\|_{H^s}+2$ $\pas$
Using Lemma \ref{cut-combine} in a similar way as in Proposition \ref{pathwise solution s>3/2 bounded initial data}, we find that
\begin{equation*}
	\left(u=\sum_{m\geq1}\textbf{1}_{\{m-1\leq\|u_0\|_{H^s}<m\}}u_m,\ \
	\tau=\sum_{m\geq1}\textbf{1}_{\{m-1\leq\|u_0\|_{H^s}<m\}}\tau_m\right)
\end{equation*}
is a solution to \eqref{SGCH problem} satisfying \eqref{L2 moment bound} and $u(0)=u_0$ almost surely. Uniqueness is given by Lemma \ref{pathwise uniqueness s>3/2}.
\qed

\subsection{Proof for \ref{Local-WP-blow-up} in Theorem \ref{Local-WP}: Blow-up criterion}

With a local solution $(u,\tau)$ at hand, one may pass from $(u,\tau)$ to the maximal solution $(u,\tau^*)$ as in \cite{GlattHoltz-Vicol-2014-AP,Breit-Feireisl-Hofmanova-2018-Book}. In the periodic setting, i.e., $x\in\T=\R/2\pi\Z$, the blow-up criterion \eqref{Blow-up criterion} for a maximal solution has been proved in \cite{Rohde-Tang-2021-JDDE}  by using energy estimate and some stopping-time techniques. When  $x\in\R$, \eqref{Blow-up criterion} can be also obtained in the same way, and we omit the details for brevity.

\subsection{Proof for \ref{Local-WP-dependence} in Theorem \ref{Local-WP}: Stability} Let  $u_0, v_0\in L^\infty({\rm \Omega};H^s)$ be two $H^s$-valued $\mathcal{F}_0$-measurable random variables. Let $u$ and $v$ be the corresponding solutions with initial conditions $u_0$ and $v_0$. To prove  \ref{Local-WP-dependence} in Theorem \ref{Local-WP}, for any $\epsilon>0$ and $T>0$, we need to find a $\delta=\delta(\epsilon,u_0,T)>0$ and a $\tau\in(0,T]$ $\pas$  such that \eqref{stability result} holds true as long as \eqref{stability condition L-infty} is satisfied.  Without loss of generality, by \eqref{stability condition L-infty}, we can first assume 
\begin{equation}\label{u0 v0 a prior}
	\|v_0\|_{L^\infty({\rm \Omega};H^s)}\leq \|u_0\|_{L^\infty({\rm \Omega};H^s)}+1.
\end{equation}
From now on $\epsilon>0$ and $T>0$ are given.

However, as is mentioned in Remark \ref{remark on LW-1}, the term $u^ku_x$ loses one regularity and the estimate for $\E\sup_{t\in[0,\tau]}\|u(t)-v(t)\|^2_{H^s}$ will involve $\|u\|_{H^{s+1}}$ or $\|v\|_{H^{s+1}}$, which might be infinite since we only know $u,v\in H^s$. To overcome this difficulty, we will consider \eqref{smooth initial problem}. 
Let $\varepsilon\in(0, 1)$. By \ref{Local-WP-existence} in Theorem \ref{Local-WP}, there is a unique solution $u_{\varepsilon}$ (resp. $v_{\varepsilon}$) to the problem \eqref{smooth initial problem} with initial data $J_\varepsilon u_0$ (resp. $J_\varepsilon v_0$). 
Then the $H^{s+1}$-norm is well-defined for the smooth solution $u_\e$ and $v_\e$.  
Similar to \eqref{tau-e eta T s>3/2}, for any $T>0$, we define
\begin{equation} \label{tau u v stability}
	\tau^{f,T}_{\e}:=\inf\left\{t\geq0:\|f_\e\|_{H^s}\geq \|J_{\e} f_0\|_{H^s}+2\right\}\wedge T,\ \ f\in\{u,v\}.
\end{equation}
Recalling the analysis in Lemma \ref{Cauchy estimates lemma} and Proposition \ref{pathwise solution s>3/2 bounded initial data} (for the case $f=v$, we notice \eqref{u0 v0 a prior}), we can use Lemma \ref{Cauchy lemma GZ} to find that there exits a unified subsequence $\{\e_n\}$ with $\e_n\rightarrow 0$ as $n\rightarrow\infty$ such that for $f\in\{u,v\}$, there is
a sequence of stopping times $\xi^f_{\e_n}$ and a stopping time $\tau^f$ satisfying
\begin{equation}\label{xi f en convergent time}
	\xi^f_{\e_n}\leq \tau^{f,T}_{\e_n},\  n\ge1\ \ \text{and}\ \ 
	\lim_{n\rightarrow\infty}\xi^f_{\e_n}=\tau^f\in(0,T] \ \ \pas,
\end{equation}
and  
\begin{equation}\label{f-en convergent f}
	\lim_{n\rightarrow\infty}\sup_{t\in[0,\tau^f]}\|f-f_{\varepsilon_n}\|_{H^s}=0, \ 
	\sup_{t\in[0,\tau^f]}\|f\|_{H^s}\leq \|f_0\|_{H^s}+2 \ \
	\pas
\end{equation}
Moreover, for $f\in\{u,v\}$, there exists ${\rm \Omega}^f_n\uparrow {\rm \Omega}$ 
such that
\begin{equation}\label{f-en Omega-n bound}
	{\bf 1}_{{\rm \Omega}^f_n }\sup_{t\in[0,\tau^f]}\|f_{\e_n}\|_{H^s}\leq \|f_0\|_{H^s}+2\ \ \pas
\end{equation}
Next, we let ${\rm \Omega}_n:={\rm \Omega}^u_n\, \cap\, {\rm \Omega}^v_n.$ 
Then ${\rm \Omega}_n\uparrow {\rm \Omega}$. This,  \eqref{f-en convergent f}, \eqref{f-en Omega-n bound} and Lebesgue's dominated convergence theorem yield
\begin{align*}
	\lim_{n\rightarrow \infty}\E\sup_{t\in[0,\tau^f]}\|f-{\bf 1}_{{\rm \Omega}_n}f_{\varepsilon_n}\|^2_{H^s}
	=\,  & 0,\ \ f\in\{u,v\}.
\end{align*}
Therefore, we have, when $n$ is large enough, that
\begin{align}
	\E\sup_{t\in[0,\tau^f]}\|f-{\bf 1}_{{\rm \Omega}_n}f_{\varepsilon_n}\|^2_{H^s}
	< \frac{\epsilon}{9},\ \ f\in\{u,v\}.
	\label{f-fe e/3}
\end{align}
Now we consider $\E\sup_{t\in[0,\tau^u\wedge \tau^v]}
\|{\bf 1}_{{\rm \Omega}_n}u_{\e_n}-{\bf 1}_{{\rm \Omega}_n}v_{\e_n}\|^2_{H^s}$. It follows from \eqref{xi f en convergent time} that for all $n\ge1$,
\begin{align}
	&\,\E\sup_{t\in[0,\tau^u\wedge \tau^v]}
	\|{\bf 1}_{{\rm \Omega}_n}u_{\e_n}-{\bf 1}_{{\rm \Omega}_n}v_{\e_n}\|^2_{H^s}\notag\\
	\leq \, &  
	\E\sup_{t\in[0,\xi^u_{\e_n}\wedge \xi^v_{\e_n}\wedge\tau^u\wedge \tau^v]}
	\|{\bf 1}_{{\rm \Omega}_n}u_{\e_n}-{\bf 1}_{{\rm \Omega}_n}v_{\e_n}\|^2_{H^s}\notag\\
	&+\E\sup_{t\in[\xi^u_{\e_n}\wedge \xi^v_{\e_n}\wedge\tau^u\wedge \tau^v,\tau^u\wedge \tau^v]}
	\|{\bf 1}_{{\rm \Omega}_n}u_{\e_n}-{\bf 1}_{{\rm \Omega}_n}v_{\e_n}\|^2_{H^s}\notag\\
	\leq \, &  
	\E\sup_{t\in[0,\tau^{u,T}_{\e_n}\wedge \tau^{v,T}_{\e_n}]}\|u_{\e_n}(t)-v_{\e_n}(t)\|^2_{H^s}\notag\\
	&+\E\sup_{t\in[\xi^u_{\e_n}\wedge \xi^v_{\e_n}\wedge\tau^u\wedge \tau^v,\tau^u\wedge \tau^v]}
	\|{\bf 1}_{{\rm \Omega}_n}u_{\e_n}-{\bf 1}_{{\rm \Omega}_n}v_{\e_n}\|^2_{H^s}.\label{u-en v-en stability}
\end{align}
By \eqref{f-en Omega-n bound}, 
\begin{align*}
	\sup_{t\in[\xi^u_{\e_n}\wedge \xi^v_{\e_n}\wedge\tau^u\wedge \tau^v,\tau^u\wedge \tau^v]}
	\|{\bf 1}_{{\rm \Omega}_n}u_{\e_n}-{\bf 1}_{{\rm \Omega}_n}v_{\e_n}\|^2_{H^s}
	\leq \,
	32\left(\|u_0\|_{H^s}^2+\|v_0\|_{H^s}^2+1\right).
\end{align*}
Consequently, by Lebesgue's dominated convergence theorem and \eqref{xi f en convergent time}, we have  for  $n\gg 1$ that,
\begin{align}\label{u-en v-en e/6}
	\E\sup_{t\in[\xi^u_{\e_n}\wedge \xi^v_{\e_n}\wedge\tau^u\wedge \tau^v,\tau^u\wedge \tau^v]}
	\|{\bf 1}_{{\rm \Omega}_n}u_{\e_n}-{\bf 1}_{{\rm \Omega}_n}v_{\e_n}\|^2_{H^s}
	< \frac{\epsilon}{18}.
\end{align}
Now we estimate $\E\sup_{t\in[0,\tau^{u,T}_{\e_n}\wedge \tau^{v,T}_{\e_n}]}\|u_{\e_n}(t)-v_{\e_n}(t)\|^2_{H^s}$. 
Similar to \eqref{Cauchy lemma step 1}, by  using \eqref{u0 v0 a prior}, one can show that for $s'\in\left(\frac{1}{2},\min\left\{s-1,\frac{3}{2}\right\}\right)$,  
\begin{align}
	&\,\E\sup_{t\in[0,\tau^{u,T}_{\e_n}\wedge \tau^{v,T}_{\e_n}]}\|u_{\e_n}(t)-v_{\e_n}(t)\|^2_{H^s}\notag\\
	\leq \, &   C \E\left\{\|J_{\e_n} u_0-J_{\e_n} v_0\|^2_{H^s}
	+\|J_{\e_n} u_0-J_{\e_n} v_0\|^2_{H^{s'}}\|J_{\e_n} u_0\|^2_{H^{s+1}}\right\}\notag\\
	&+C\E\sup_{t\in[0,\tau^{u,T}_{\e_n}\wedge \tau^{v,T}_{\e_n}]}\|u_{\e_n}(t)-v_{\e_n}(t)\|^2_{H^{s'}}\notag\\
	\leq \, &   C \E\Big\{\|u_0-v_0\|^2_{H^s}
	+\frac{1}{\e_n^2}\| u_0-v_0\|^2_{H^{s'}}\|u_0\|^2_{H^{s}}\Big\}\notag\\
	&+C\E\sup_{t\in[0,\tau^{u,T}_{\e_n}\wedge \tau^{v,T}_{\e_n}]}\|u_{\e_n}(t)-v_{\e_n}(t)\|^2_{H^{s'}},
	\label{u-en v-en difference}
\end{align}
where 
$C=C\left(\|u_0\|_{L^\infty({\rm \Omega};H^s)},T\right)$ and Lemma \ref{mollifier properties} is used in the last step. Since $u_0\in L^\infty({\rm \Omega};H^s)$, 
by Lemmas \ref{Uniqueness bounded data} and \ref{mollifier properties} again, we have
\begin{align}
	&\E\sup_{t\in[0,\tau^{u,T}_{\e_n}\wedge \tau^{v,T}_{\e_n}]}\|u_{\e_n}(t)-v_{\e_n}(t)\|^2_{H^s}\notag\\
	\leq \, &   C\E\left\{\|u_0-v_0\|^2_{H^s}
	+\frac{1}{\e_n^2}\|u_0-v_0\|^2_{H^{s'}}\|u_0\|^2_{H^{s}}\right\}
	+C\E\|J_{\e_n} u_0-J_{\e_n} v_0\|^2_{H^{s'}}\notag\\
	\leq \, &   C\E\|u_0-v_0\|^2_{H^s}
	+C\frac{1}{\e_n^2}\E\|u_0-v_0\|^2_{H^{s'}}
	+C\E\|u_0-v_0\|^2_{H^{s'}},
	\label{u-en v-en n to be fixed}
\end{align}
where $C=C\left(\|u_0\|_{L^\infty({\rm \Omega};H^s)},T\right)$ as before.
Fix $n=n_0\gg 1$ such that
\eqref{f-fe e/3} and \eqref{u-en v-en e/6} are satisfied, i.e., 
\begin{equation}\label{f-en u-en v-en n fixed stability}
	\left\{\begin{aligned}
		&\E\sup_{t\in[0,\tau^f]}\|f-{\bf 1}_{{\rm \Omega}_{n_0}}f_{\varepsilon_{n_0}}\|^2_{H^s}
		<\frac{\epsilon}{9},\ \ f\in\{u,v\},\\
		&\E\sup_{t\in[\xi^u_{\e_{n_0}}\wedge \xi^v_{\e_{n_0}}\wedge\tau^u\wedge \tau^v,\tau^u\wedge \tau^v]}
		\|{\bf 1}_{{\rm \Omega}_{n_0}}u_{\e_{n_0}}-{\bf 1}_{{\rm \Omega}_{n_0}}v_{\e_{n_0}}\|^2_{H^s}
		< \frac{\epsilon}{18}.
	\end{aligned}\right.
\end{equation}
Then, for \eqref{u-en v-en n to be fixed} with $n=n_0$, we can find a $\delta=\delta(\epsilon,u_0,T)\in(0,1)$ such that \eqref{u0 v0 a prior} is satisfied and
\begin{align}
	&\E\sup_{t\in[0,\tau^{u,T}_{\e_{n_0}}\wedge \tau^{v,T}_{\e_{n_0}}]}\|u_{\e_{n_0}}(t)-v_{\e_{n_0}}(t)\|^2_{H^s}<\frac{\epsilon}{18},\ \ \text{if}\ \ \|u_0-v_0\|_{L^\infty({\rm \Omega};H^s)}<\delta.\label{u-en v-en n fixed e/6}
\end{align}
As a result,
for \eqref{u-en v-en stability} with fixed $n=n_0$, we 
use \eqref{f-en u-en v-en n fixed stability}$_2$ and \eqref{u-en v-en n fixed e/6} to derive that 
\begin{align*}
	\E\sup_{t\in[0,\tau^u\wedge \tau^v]}\|{\bf 1}_{{\rm \Omega}_{n_0}}u_{\e_{n_0}}-{\bf 1}_{{\rm \Omega}_{n_0}}v_{\e_{n_0}}\|^2_{H^s}
	\leq \, &   \frac{\epsilon}{18}+\frac{\epsilon}{18}=\ \frac{\epsilon}{9},\ \ \text{if}\ \  \|u_0-v_0\|_{L^\infty({\rm \Omega};H^s)}<\delta.
\end{align*}
This inequality and \eqref{f-en u-en v-en n fixed stability}$_1$ yield that
\begin{align*}
	&\,\E\sup_{t\in[0,\tau^u\wedge \tau^v]}\|u-v\|^2_{H^s}\\
	\leq \, &
	3\sum_{f\in\{u,v\}}\E\sup_{t\in[0,\tau^u\wedge \tau^v]}\|f-{\bf 1}_{{\rm \Omega}_{n_0}}f_{\e_{n_0}}\|^2_{H^s}\\
	&+3\E\sup_{t\in[0,\tau^u\wedge \tau^v]}\|{\bf 1}_{{\rm \Omega}_{n_0}}u_{\e_{n_0}}-{\bf 1}_{{\rm \Omega}_{n_0}}v_{\e_{n_0}}\|^2_{H^s}\\
	\leq \, &   \frac{\epsilon}{3}+\frac{\epsilon}{3}+\frac{\epsilon}{3}=\epsilon,\ \ \text{if}\ \ \|u_0-v_0\|_{L^\infty({\rm \Omega};H^s)}<\delta.
\end{align*}
Hence we obtain \eqref{stability result} with $\tau=\tau^u\wedge\tau^v$. Due to \eqref{xi f en convergent time}, $\tau\in(0,T]$ almost surely. 

\begin{Remark}\label{Continuous dependence reamrk}
	Here we remark that the restriction ${\bf 1}_{{\rm \Omega}_{n}}$ is needed  to estimate $$\E\sup_{t\in[0,\tau^f]}\|f-{\bf 1}_{{\rm \Omega}_{n}}f_{\varepsilon_{n}}\|^2_{H^s}$$ for $f\in\{u,v\}$. This is because we only have $\lim_{n\rightarrow\infty}\sup_{t\in[0,\tau^f]}\|f-f_{\varepsilon_n}\|_{H^s}=0$ $\pas$  (cf. \ref{GZ lemma-convergence} in Lemma \ref{Cauchy lemma GZ}), and we need to interchange limit and expectation. By \ref{GZ lemma-moment} in Lemma \ref{Cauchy lemma GZ}, 
	$$\sup_{t\in[0,\tau^f]}\|f-{\bf 1}_{{\rm \Omega}_{n}}f_{\varepsilon_{n}}\|^2_{H^s}\leq2\sup_{t\in[0,\tau^f]}\|f\|^2_{H^s}+2{\bf 1}_{{\rm \Omega}_{n}}\sup_{t\in[0,\tau^f]}\|f_{\varepsilon_{n}}\|^2_{H^s}\leq 4\|f_0\|^2_{H^s}+16.$$ Hence Lebesgue's dominated convergence theorem can be used. In the deterministic case, one can directly consider $\|f-f_{\varepsilon_{n}}\|^2_{H^s}.$
\end{Remark}

\section{Weak instability}\label{sect:weak instability}
Now we prove Theorem \ref{Weak instability}. As is mentioned in Remark \ref{remark on instability 1}, since we can not  get an explicit expression of the solution to \eqref{SGCH problem}, we start with constructing some approximative solutions from which \eqref{sup sin t} can be established.

\subsection{Approximative solutions and actual solutions}

Following the approach in \cite{Himonas-Holliman-2014-ADE, Rohde-Tang-2021-JDDE},  now we construct the approximative solutions. We fix two functions $\phi,\tilde{\phi}\in  C_c^{\infty}$ such that
\begin{equation}\label{phi phi-tilde}
	\phi(x)=\left\{\begin{aligned}
		&1,\ \text{if}\ |x|<1,\\
		&0,\ \text{if}\ |x|\geq 2,
	\end{aligned}\right.
	\ \ \text{and} \ \ \tilde{\phi}(x)=1\ \text{if}\ x\in {\rm supp}\ \phi.
\end{equation}
Let $k\ge1$ and
\begin{equation}\label{m k relation}
	m\in\{-1,1\}\ \text{if}\ k \ \text{is odd and} \ m\in\{0,1\}\ \text{if}\ k \ \text{is even}.
\end{equation}
Then we consider the following sequence of approximative solutions
\begin{equation}\label{approximation solution u-m n}
	u_{m,n}=u_l+u_h,
\end{equation}
where  $u_h= u_{h,m,n}$ is the high-frequency part defined by
\begin{equation}\label{high-frequency approximation solutions}
	u_h = u_{h,m,n}(t,x) = n^{-\frac{\delta}{2}-s}\phi\left(\frac{x}{n^{\delta}}\right)\cos(nx-mt),\ \ n\in\N,
\end{equation}
and $u_l=u_{l,m,n}$ is the low-frequency part constructed such that $u_l$ is  the  solution to the following  problem:
\begin{equation}\label{low-frequency equation}
	\left\{\begin{aligned}
		&\partial_tu_l+u_l^k\partial_xu_l+F(u_l)=0,\quad x\in\R, \ t>0,\  k\geq 1,\\
		&u_l(0,x)=mn^{-\frac{1}{k}}\tilde{\phi}\left(\frac{x}{n^{\delta}}\right),\quad  x\in\R.
	\end{aligned}
	\right.
\end{equation}
The parameter $\delta>0$ in \eqref{high-frequency approximation  solutions} and \eqref{low-frequency equation} will be determined later for different $k\ge1$. Particularly, when $m=0$, we have $u_l=0$. In this case the approximative solution $u_{0,n}$ has \textit{no} low-frequency part and
\begin{equation*}
	u_{0,n}(t,x)=n^{-\frac{\delta}{2}-s}\phi\left(\frac{x}{n^{\delta}}\right)\cos(nx).
\end{equation*}

Next, we consider the problem \eqref{SGCH problem} with initial data $u_{m,n}(0,x)$, i.e.,
\begin{equation}\label{SGCH m n}
	\left\{\begin{aligned}
		&\, {\rm d}u+ [u^k\partial_xu+F(u)]\dd t=h(t,u) \dd \W,\ \ t>0,\ x\in\R,\\
		&u(0,x)=mn^{-\frac{1}{k}}\tilde{\phi}\left(\frac{x}{n^{\delta}}\right)+n^{-\frac{\delta}{2}-s}\phi\left(\frac{x}{n^{\delta}}\right)\cos(nx),\ \ x\in\R,
	\end{aligned} \right.
\end{equation}
where $F(\cdot)$ is defined by \eqref{F decomposition}.  Since $h$ satisfies \ref{Ass2-growth-Lip}, similar to the proof for Theorem \ref{Local-WP}, we see that for each fixed $n\in\N$, \eqref{SGCH m n} has a unique solution $(u^{m,n},\tau^{m,n})$ such that $u^{m,n}\in C\left([0,\tau^{m,n}];H^s\right)$ $\pas$  with $s>5/2$.

\subsection{Estimates on the errors}

Substituting \eqref{approximation solution u-m n} into \eqref{SGCH problem}, we define the error $\EE(\omega,t,x)$ as
\begin{align*}
	\EE(\omega,t,x)
	:=\,  & u_{m,n}(t,x)-u_{m,n}(0,x)\\
	&+ \int_0^t\left[u_{m,n}^k\partial_xu_{m,n}+F(u_{m,n})\right]\dd t'-\int_0^th(t',u_{m,n}) \dd \W\ \ \pas
\end{align*}
For simplicity, we let
\begin{equation}\label{Z nonlinear}
	\ZZ_{q}=\ZZ_q(u_h,u_l)=\left\{
	\begin{aligned}
		&\sum_{j=1}^{q}C_{q}^{j}u_l^{q-j}u_h^j,\ \text{if}\ q\geq1,\\
		&0,\ \text{if}\  q=0,
	\end{aligned}\right.
\end{equation}
where $C_{q}^{j}$ is the binomial coefficient. 
By using \eqref{approximation solution u-m n}, \eqref{low-frequency equation} and \eqref{Z nonlinear}, $\EE(\omega,t,x)$ can be reformulated as
\begin{align}
	\,& \EE(\omega,t,x)\notag\\
	=\,  & u_l(t,x)-u_l(0,x)+\int_0^tu_{l}^k\partial_xu_l\dd t'+\int_0^tF(u_l)\dd t'\notag\\
	&+u_h(t,x)-u_h(0,x)+\int_0^t\left[u_l^k\partial_xu_h+\ZZ_k(\partial_xu_l+\partial_xu_h)\right]\dd t'\notag\\
	&+\int_0^t\left[F(u_l+u_h)-F(u_l)\right]\dd t'-\int_0^th(t',u_{m,n}) \dd \W\notag\\
	=\,  & u_h(t,x)-u_h(0,x)+\int_0^t\left[u_l^k\partial_xu_h+\ZZ_k(\partial_xu_l+\partial_xu_h)\right]\dd t'\notag\\
	&+\int_0^t\left[F(u_l+u_h)-F(u_l)\right]\dd t'-\int_0^th(t',u_{m,n}) \dd \W \ \ \pas
	\label{error EE}
\end{align}

\subsubsection{Estimates on the low-frequency part}
The following lemma gives a decay estimate for the low-frequency part of $u_{m,n}$, that is, $u_l$.
\begin{Lemma}\label{lemma ul Hr}
	Let $k\ge1$, $|m|=1$ or $m=0$, $s>3/2$, $\delta\in(0,2/k)$ and $n\gg 1$. Then there is a $T_l>0$ such that for all $n\gg 1$, the initial value problem \eqref{low-frequency equation} has a unique smooth solution $u_l=u_{l,m,n}\in C([0,T_l];H^s)$ such that $T_l$ does not depend on $n$. Besides, for all $r>0$, there is a constant $C=C_{r,\tilde{\phi},T_l}>0$ such that $u_l$ satisfies 
	\begin{equation}\label{ul Hr estimate}
		\|u_l(t)\|_{H^r}\leq C|m|n^{\frac{\delta}{2}-\frac{1}{k}},\ \  t\in[0,T_l].
	\end{equation}
\end{Lemma}

\proof
When $m=0$, as is mentioned above, $u_l\equiv0$ for all $t\ge0$. It remains to prove the case $|m|=1$.
For any fixed $n\ge 1$, since $u_l(0,x)\in H^\infty$, by applying Theorem \ref{Local-WP} with  $h=0$ and deterministic initial data,  we see that for any $s>3/2$, \eqref{low-frequency equation} has a unique (deterministic) solution $u_l=u_{l,m,n}\in C\left([0,T_{m,n}];H^s\right)$. 
Different from the stochastic case, here we will show that there is a lower bound of the existence time, i.e., there is a $T_l>0$ such that for all $n\gg 1$, $u_l=u_{l,m,n}$ exists on $[0,T_l]$ and satisfies \eqref{ul Hr estimate}.

\textbf{Step 1: Estimate $\|u_l(0,x)\|_{H^r}$.} 
When $n\gg1$, we have  
\begin{align*}
	\|u_l(0,x)\|^2_{H^r}
	=\,  & m^2n^{2\delta-\frac{2}{k}}
	\int_\R(1+|\xi|^2)^r\left|\widehat{\tilde{\phi}}(n^{\delta}\xi)\right|^2\, {\rm d}\xi\\
	=\,  & m^2n^{\delta-\frac{2}{k}}
	\int_\R\left(1+\left|\frac{z}{n^{\delta}}\right|^2\right)^r\left|\widehat{\tilde{\phi}}(z)\right|^2\, {\rm d}z
	\leq C m^2n^{\delta-\frac{2}{k}}
\end{align*}
for some constant $C=C_{r,\tilde{\phi}}>0$. As a result, we have
\begin{align*}
	\|u_l(0,x)\|_{H^r}
	\leq C |m|n^{\frac{\delta}{2}-\frac{1}{k}}.
\end{align*}

\textbf{Step 2: Prove \eqref{ul Hr estimate} for $r>3/2$.} 
In this case, we apply Lemmas \ref{Kato-Ponce commutator estimate} and \ref{F lemma}, $H^r\hookrightarrow 
\Wlip$ to find 
\begin{align*}
	\,&\frac{1}{2}\frac{\, {\rm d}}{\dd t}\|u_l\|^2_{H^r}\\
	\leq  \, &    \left|\left(D^{r}u_l,D^{r}\left(u_l^k\partial_xu_l\right)\right)_{L^2}\right|+\left|\left(D^{r}u_l,D^{r}F(u_l)\right)_{L^2}\right|\\
	\leq  \, &    \left|\left([D^{r},u_l^k]\partial_xu_l,D^{r}u_l\right)_{L^2}\right|+\left|\left(u_l^kD^{r}\partial_xu_l,D^{r}u_l\right)_{L^2}\right|+\left\|u_l\right\|_{H^r}\left\|F(u_l)\right\|_{H^r}\\
	\lesssim \, &  \|u_l^k\|_{H^r}\|\partial_xu_l\|_{L^{\infty}}\|u_l\|_{H^r}
	+\|\partial_xu_l\|_{L^{\infty}}\|u_l\|_{L^\infty}^{k-1}\|u_l\|_{H^r}^2
	+\|u_l\|^{k}_{\Wlip}\|u_l\|^2_{H^r}\\
	\leq  \, &     C\|u_l\|^{k+2}_{H^r},\ \ C=C_r>0.
\end{align*}
Solving the above inequality gives
\begin{equation*}
	\|u_l\|_{H^r}\leq \frac{\|u_l(0)\|_{H^r}}{\left(1-Ckt\|u_l(0)\|^k_{H^r}\right)^{\frac{1}{k}}},\ \ 
	0\le t<\frac{1}{Ck\|u_l(0)\|^k_{H^r}}.
\end{equation*}
Therefore, we arrive at
\begin{equation}
	\|u_l\|_{H^r}\leq 2 \|u_l(0)\|_{H^r},\ \ t\in[0,T_{m,n}],\ \ T_{m,n}=\frac{1}{2Ck\|u_l(0)\|^k_{H^r}}.\label{solution bound T u_0}
\end{equation}
By \textbf{Step 1}, we have 
$T_{m,n}\gtrsim \frac{1}{2Ckn^{k(\frac{\delta}{2}-\frac{1}{k})}}\rightarrow\infty,\ \text{as}\ n\rightarrow\infty$.
Consequently, we can find a common time interval $[0,T_{l}]$ such that 
\begin{equation*}
	\|u_l\|_{H^r}\leq 2 \|u_l(0)\|_{H^r}\leq C|m|n^{\frac{\delta}{2}-\frac{1}{k}},\ \ t\in[0,T_l],
\end{equation*}
which is \eqref{ul Hr estimate}.

\textbf{Step 3: Prove \eqref{ul Hr estimate} for $0<r\leq 3/2$.} Similarly, by applying Lemmas \ref{Kato-Ponce commutator estimate} and \ref{F lemma}, we have
\begin{align*}
	\,&\frac{1}{2}\frac{\, {\rm d}}{\dd t}\|u_l\|^2_{H^r}\\
	\leq  \, &    \left|\left(D^{r}u_l,D^{r}\left(u_l^k\partial_xu_l\right)\right)_{L^2}\right|+\left|\left(D^{r}u_l,D^{r}F(u_l)\right)_{L^2}\right|\\
	\leq  \, &    \left|\left([D^{r},u_l^k]\partial_xu_l,D^{r}u_l\right)_{L^2}\right|+\left|\left(u_l^kD^{r}\partial_xu_l,D^{r}u_l\right)_{L^2}\right|+\left\|u_l\right\|_{H^r}\left\|F(u_l)\right\|_{H^r}\\
	\lesssim \, &  \|u_l^k\|_{H^r}\|\partial_xu_l\|_{L^{\infty}}\|u_l\|_{H^r}
	+\|\partial_xu_l\|_{L^{\infty}}\|u_l\|_{L^\infty}^{k-1}\|u_l\|_{H^r}^2\\
	&+\left\|u_l\right\|_{H^r}\|u_l\|^{k}_{\Wlip}\left(\|u_l\|_{H^r}+\|\partial_x u_l\|_{H^{r}}\right).
\end{align*}
It follows from the embedding $H^{r+\frac32}\hookrightarrow H^{r+1}$ and 
$H^{r+\frac32}\hookrightarrow W^{1,\infty}$ 
that
\begin{align*}
	\frac{1}{2}\frac{\, {\rm d}}{\dd t}\|u_l\|^2_{H^r}\lesssim \, &  \|u_l^k\|_{H^r}\|\partial_xu_l\|_{L^{\infty}}\|u_l\|_{H^r}
	+\|\partial_xu_l\|_{L^{\infty}}\|u_l\|_{L^\infty}^{k-1}\|u_l\|_{H^r}^2\\
	&+\left\|u_l\right\|_{H^r}\|u_l\|^{k}_{\Wlip}\left(\|u_l\|_{H^r}+\|\partial_x u_l\|_{H^{r}}\right)\\
	\lesssim \, &  \|u_l\|^k_{\Wlip}\|u_l\|^2_{H^r}
	+\|u_l\|^{k}_{\Wlip}\left\|u_l\right\|_{H^r}\left\|u_l\right\|_{H^{r+1}}\\
	\lesssim \, &  \|u_l\|_{H^{r+\frac{3}{2}}}^{k}\|u_l\|_{H^{r}}^2+\|u_l\|_{H^{r}}\|u_l\|_{H^{r+\frac{3}{2}}}^{k+1}.
\end{align*}
Using the conclusion of \textbf{Step 2} for $r+\frac{3}{2}>\frac{3}{2}$, we have
\begin{align*}
	\frac{\, {\rm d}}{\dd t}\|u_l\|_{H^r}\lesssim\|u_l\|_{H^{r}}\|u_l(0)\|_{H^{r+\frac{3}{2}}}^{k}+\|u_l(0)\|_{H^{r+\frac{3}{2}}}^{k+1},\ \ t\in[0,T_l],
\end{align*}
and hence
\begin{align*}
	\|u_l(t)\|_{H^r}\lesssim \|u_l(0)\|_{H^r}+\|u_l(0)\|_{H^{r+\frac{3}{2}}}^{k+1}T_l+\int_0^t\|u_l\|_{H^{r}}\|u_l(0)\|_{H^{r+\frac{3}{2}}}^{k}\dd t',\ \ t\in[0,T_l].
\end{align*}
Applying Gr\"{o}nwall's inequality to the above inequality, we have
\begin{align*}
	\|u_l\|_{H^r}\lesssim\left(\|u_l(0)\|_{H^r}+\|u_l(0)\|_{H^{r+\frac{3}{2}}}^{k+1}T_l\right)\exp\left\{\|u_l(0)\|_{H^{r+\frac{3}{2}}}^{k}T_l\right\},\ \ t\in[0,T_l].
\end{align*}
Since $\delta\in(0,2/k)$, we can infer from \textbf{Step 1} that $\exp\left\{\|u_l(0)\|_{H^{r+\frac{3}{2}}}^{k}T_l\right\} <C(T_l)$ for some constant $C(T_l)>0$ and
$\|u_l(0)\|_{H^{r+\frac{3}{2}}}^{k+1}\leq \|u_l(0)\|_{H^{r+\frac{3}{2}}}\leq C|m|n^{\frac{\delta}{2}-\frac{1}{k}}$.
Hence we see that there is a constant $C=C_{r,\tilde{\phi},T_l}>0$ such that
\begin{align*}
	\|u_l\|_{H^r}\leq C|m|n^{\frac{\delta}{2}-\frac{1}{k}},\ \ t\in[0,T_l],
\end{align*}
which is \eqref{ul Hr estimate}. 
\qed

Recall the approximative solution defined by \eqref{approximation solution u-m n}.
The above result means that the $H^s$-norm of the low-frequency part $u_l$ is decaying. For the high-frequency part $u_h$, as in Lemma \ref{lemma uh Hr}, its $H^s$-norm  is bounded.

\subsubsection{Estimate on $\EE$}
Recall the error $\EE$ given in \eqref{error EE}. By using \eqref{phi phi-tilde} and \eqref{m k relation}, we have $m=m^k$ and $\phi=\tilde{\phi}^k\phi$ for all $k\ge1$. Then by \eqref{high-frequency approximation solutions} and $u_l(0,x)$ in \eqref{low-frequency equation}, we see that as long as $m\neq0$, 
\begin{align*}
	\,&u_{h}(t,x)-u_{h}(0,x)\notag\\
	=\,  & n^{-\frac{\delta}{2}-s}\phi\left(\frac{x}{n^{\delta}}\right)\cos(nx-mt)-n^{-\frac{\delta}{2}-s}\phi\left(\frac{x}{n^{\delta}}\right)\cos(nx)\\
	=\,  & m^{-1}m^k\tilde{\phi}^k\left(\frac{x}{n^{\delta}}\right)n^{-\frac{\delta}{2}-s}\phi\left(\frac{x}{n^{\delta}}\right)\cos(nx-mt)\\
	&-m^{-1}m^k
	\tilde{\phi}^k\left(\frac{x}{n^{\delta}}\right)
	n^{-\frac{\delta}{2}-s}\phi\left(\frac{x}{n^{\delta}}\right)\cos(nx)\notag\\
	=\,  & 
	m^{-1}u_l^k(0,x)n^{1-\frac{\delta}{2}-s}\phi\left(\frac{x}{n^{\delta}}\right)\cos(nx-mt)\\
	&-m^{-1}u_l^k(0,x)n^{1-\frac{\delta}{2}-s}\phi\left(\frac{x}{n^{\delta}}\right)\cos(nx)\\
	=\,  & \int_0^tu_l^k(0,x)n^{1-\frac{\delta}{2}-s}\phi\left(\frac{x}{n^{\delta}}\right)\sin(nx-mt')\dd t'.
\end{align*}
If $m=0$, then $u_l=0$ and we also have
\begin{align*}
	u_{h}(t,x)-u_{h}(0,x)=0=\int_0^t 0\  \dd t'
	=\int_0^tu_l^k(0,x)n^{1-\frac{\delta}{2}-s}\phi\left(\frac{x}{n^{\delta}}\right)\sin(nx-mt')\dd t'.
\end{align*}
To sum up, we find that for all $k\ge1$, $m$ satisfying \eqref{m k relation}, $u_h$ given by \eqref{high-frequency approximation solutions} and $u_l(0,x)$ in \eqref{low-frequency equation},
\begin{align}\label{uh(t)-uh(0)}
	u_{h}(t,x)-u_{h}(0,x)
	=\int_0^tu_l^k(0,x)n^{1-\frac{\delta}{2}-s}\phi\left(\frac{x}{n^{\delta}}\right)\sin(nx-mt')\dd t'.
\end{align}
On the other hand, for all $k\ge1$,
\begin{align}
	\int_0^tu_{l}^k\partial_xu_h\dd t'=\, & -\int_0^tu_{l}^k(t')n^{1-\frac{\delta}{2}-s}\phi\left(\frac{x}{n^{\delta}}\right)\sin(nx-mt')\dd t'\notag\\
	&+\int_0^tu_{l}^k(t')n^{-\frac{3\delta}{2}-s}\partial_x\phi\left(\frac{x}{n^{\delta}}\right)\cos(nx-mt')\dd t'\label{ul k uh_x}.
\end{align}
Combining \eqref{uh(t)-uh(0)}, \eqref{ul k uh_x} and \eqref{F decomposition} into \eqref{error EE} yields
\begin{align}
	\EE(\omega,t,x)
	=\,  & \sum_{i=1}^4\int_0^tE_i\dd t'-\int_0^th(t',u_{m,n}) \dd \W,\ \  k\ge1\ \ \pas,\label{Error}
\end{align}
where
\begin{align*}
	E_1:=\, & [u_{l}^k(0)-u_{l}^k(t)]n^{1-\frac{\delta}{2}-s}\phi\left(\frac{x}{n^{\delta}}\right)\sin(nx-mt)\notag\\
	&+u_{l}^k(t)n^{-\frac{3\delta}{2}-s}\partial_x\phi\left(\frac{x}{n^{\delta}}\right)\cos(nx-mt)
	+\ZZ_k(\partial_xu_l+\partial_xu_h),\\
	E_2:=\, & F_1(u_l+u_h)-F_1(u_l)=D^{-2}\partial_x\ZZ_{k+1},\\
	E_3:=\, & F_2(u_l+u_h)-F_2(u_l)\notag\\
	=\, &   \frac{2k-1}{2}D^{-2}\partial_x\Big\{u_l^{k-1}
	\left[2(\partial_xu_l)(\partial_xu_h)+(\partial_xu_h)^2\right]
	+\ZZ_{k-1}(\partial_xu_l+\partial_xu_h)^2\Big\},\\
	E_4:=\,  &F_3(u_l+u_h)-F_3(u_l)\notag\\
	=\,  & \frac{k-1}{2}D^{-2}\Big\{u_l^{k-2}\big[3(\partial_xu_l)^2(\partial_xu_h)+3(\partial_xu_l)(\partial_xu_h)^2+(\partial_xu_h)^3\big]\nonumber\\
	&\hspace{1.8cm}+\ZZ_{k-2}(\partial_xu_l+\partial_xu_h)^3\Big\}. 
\end{align*}
We remark here that $E_4$ disappears  when $k=1$.
Recalling $\rho_0\in(1/2,1)$ in Hypothesis \ref{Assumption-2}, now we shall estimate the $H^{\rho_0}$-norm of the error $\EE$. Actually, we will show that the $H^{\rho_0}$-norm of $\EE$ is decaying.
\begin{Lemma}\label{Lemma Error estimate}
	Let $T_l>0$ be given in Lemma \ref{lemma ul Hr}, and $\rho_0\in(1/2,1)$ be given in \ref{Assumption-2}. Let $n\gg 1$, $s>5/2$. Let
	\begin{equation}\label{delta}
		\left\{\begin{aligned}
			&\frac{2}{3}<\delta<1,\ \text{when}\ k=1,\\
			&\frac2k-\frac{2}{2k-1}<\delta<\frac1k,\ \text{when}\ k\ge 2,
		\end{aligned}\right.
	\end{equation}
	and 
	\begin{equation}\label{rs}
		0>r_s=-s-1+\rho_0+k\delta,\ \ k\ge1.
	\end{equation}
	Then the error $\EE$ given by \eqref{Error} satisfies that for some  $C=C(T_l)>0$,
	\begin{align*}
		\E\sup_{t\in[0,T_l]}\|\EE(t)\|_{H^{\rho_0}}^2\leq Cn^{2r_s}.
	\end{align*}

\end{Lemma}

\proof
The proof is technical and it is given in Appendix \ref{sec:appendix}.
\qed

\subsubsection{Estimate on $u_{m,n}-u^{m,n}$}

Recall the approximative solutions $u_{m,n}$ given by \eqref{approximation solution u-m n}. Then we have the following estimates on the difference between the actual solutions and the approximative  solutions.
\begin{Lemma}\label{difference estimate lemma}
	Let $k\ge1$, $s>5/2$ and $\rho_0$ be given in \ref{Assumption-2}. Let \eqref{delta} hold true and $r_s<0$ be given in \eqref{rs}.      For any $R>1$,  we define
	\begin{align}\label{tau actual solution mnR}
		\tau^{m,n}_R:=\inf\{t>0:\|u^{m,n}\|_{H^{s}}>R\}.
	\end{align}
	Then for $n\gg 1$,
	\begin{align}
		&\E\sup_{t\in[0,T_l\wedge\tau^{m,n}_R]}\|(u_{m,n}-u^{m,n})(t)\|_{H^{\rho_0}}^2 \leq Cn^{2r_s},\label{Error sigma}
	\end{align}
	\begin{align}
		&\E\sup_{t\in[0,T_l\wedge\tau^{m,n}_R]}\|(u_{m,n}-u^{m,n})(t)\|_{H^{2s-\rho_0}}^2 \leq Cn^{2s-2\rho_0},\label{Error 2s minues sigma}
	\end{align}
	where $T_l>0$ is given in Lemma \ref{lemma ul Hr} and $C=C(R,T_l)>0$.
\end{Lemma}
\proof
Let $v=v_{m,n}=u_{m,n}-u^{m,n}$. Then $v$ satisfies $v(0)=0$ and
\begin{align*}
	v(t)+\int_0^t\left(\frac{1}{k+1}\partial_x(Pv)+F(u_{m,n})-F(u^{m,n})\right)&\dd t'\\
	=-\int_0^t h(t',u^{m,n}) &\dd \W+\sum_{i=1}^4\int_0^tE_i\dd t',
\end{align*}
where 
\begin{align*}
	P=P_{m,n}=u_{m,n}^k+u_{m,n}^{k-1}u^{m,n}+\cdots+u_{m,n}(u^{m,n})^{k-1}+(u^{m,n})^{k},\quad k\ge1.
\end{align*}
On $[0,T_l]$, by the It\^{o} formula, we have that
\begin{align*}
	\|v(t)\|_{H^{\rho_0}}^2=\, & -2\int_0^t\left(h(t',u^{m,n}) \dd \W,v\right)_{H^{\rho_0}}+2\sum_{i=1}^4\int_0^t(E_i,v)_{H^{\rho_0}}\dd t'\\
	&-\frac{2}{k+1}\int_0^t(\partial_x(Pv),v)_{H^{\rho_0}}\dd t'
	-2\int_0^t([F(u_{m,n})-F(u^{m,n})],v)_{H^{\rho_0}}\dd t'\\
	&+\int_0^t\|h(t',u^{m,n})\|_{\LL_2(\U;H^{\rho_0})}^2\dd t'.
\end{align*}
Taking supremum with respect to $t\in[0,T_l\wedge\tau^{m,n}_R]$, and then using the BDG inequality yield
\begin{align*}
	&\,\E\sup_{t\in[0,T_l\wedge\tau^{m,n}_R]}\|v(t)\|_{H^{\rho_0}}^2\\
	\leq  \, &    C\E\left(\int_0^{T_l\wedge\tau^{m,n}_R}\|v\|_{H^{\rho_0}}^2\|h(t,u^{m,n})\|_{\LL_2(\U;H^{\rho_0})}^2\dd t\right)^{1/2}\\
	&+2\sum_{i=1}^4\E\int_0^{T_l\wedge\tau^{m,n}_R}\left|(E_i,v)_{H^{\rho_0}}\right| \dd t\\
	&+\frac{2}{k+1}\E\int_0^{T_l\wedge\tau^{m,n}_R}\left|(\partial_x(Pv),v)_{H^{\rho_0}}\right| \dd t\\
	&+2\E\int_0^{T_l\wedge\tau^{m,n}_R}\big|\big(F(u_{m,n})-F(u^{m,n}),v\big)_{H^{\rho_0}}\big| \dd t\\
	&+\E\int_0^{T_l\wedge\tau^{m,n}_R}\|h(t,u^{m,n})\|_{\LL_2(\U;H^{\rho_0})}^2\dd t.
\end{align*}
It follows from Lemmas \ref{lemma uh Hr} and \ref{lemma ul Hr} that $\|u_{m,n}\|_{H^{s}}\lesssim 1$ on $[0,\tau^{m,n}_R\wedge T_l]$. Hence we can infer from Hypothesis {\rm\ref{Assumption-2}} that 
\begin{align*}
	\|h(t,u^{m,n})\|_{\LL_2(\U;H^{\rho_0})}^2
	\lesssim \, &  \|h(t,u_{m,n})\|_{\LL_2(\U;H^{\rho_0})}^2+\|h(t,u_{m,n})-h(t,u^{m,n})\|_{\LL_2(\U;H^{\rho_0})}^2\\
	\lesssim \, &  \left({\rm e}^\frac{-1}{\|u_{m,n}\|_{H^{\rho_0}}}\right)^2
	+g^2_3(CR)\|v\|_{H^{\rho_0}}^2,\ \ t\in[0,\tau^{m,n}_R\wedge T_l]\ \ \pas,
\end{align*}
where $g_3(\cdot)$ is given in \ref{Ass2-rho-0}. As a result, for any fixed $s>5/2$, by applying Lemmas \ref{lemma uh Hr} and \ref{lemma ul Hr} again, we can pick $N=N(s,k)\gg 1$ to derive
\begin{align*}
	\|h(t,u^{m,n})\|_{\LL_2(\U;H^{\rho_0})}^2
	\lesssim \, &  \|u_{m,n}\|^{2N}_{H^{\rho_0}}
	+g^2_3(CR)  \|v\|_{H^{\rho_0}}^2\\
	\lesssim \, & \left(n^{N(-s+\rho_0)}+n^{N(\frac{\delta}{2}-\frac{1}{k})}\right)^2
	+g^2_3(CR)  \|v\|_{H^{\rho_0}}^2\\
	\lesssim \, & n^{2r_s}+g^2_3(CR)  \|v\|_{H^{\rho_0}}^2,\ \ t\in[0,\tau^{m,n}_R\wedge T_l]\ \ \pas
\end{align*}
Consequently, we can infer from the above inequalities that
\begin{align*}
	&\,\E\sup_{t\in[0,T_l\wedge\tau^{m,n}_R]}\|v(t)\|_{H^{\rho_0}}^2\\
	\leq  \, &    \frac12\E\sup_{t\in[0,T_l\wedge\tau^{m,n}_R]}\|v(t)\|_{H^{\rho_0}}^2+2\sum_{i=1}^4\E\int_0^{T_l\wedge\tau^{m,n}_R}\left|(E_i,v)_{H^{\rho_0}}\right| \dd t\\
	&+\frac{2}{k+1}\E\int_0^{T_l\wedge\tau^{m,n}_R}\left|(\partial_x(Pv),v)_{H^{\rho_0}}\right| \dd t\\
	&+2\E\int_0^{T_l\wedge\tau^{m,n}_R}
	\left|\big(F(u_{m,n})-F(u^{m,n}),v\big)_{H^{\rho_0}}\right| \dd t\\
	&+C(T_l)n^{2r_s}+C_R\int_0^{T_l}\E\sup_{t'\in[0,t\wedge\tau^{m,n}_R]}\|v(t')\|_{H^{\rho_0}}^2\dd t.
\end{align*}
Via Lemma \ref{Lemma Error estimate}, we have
\begin{align*}
	2\sum_{i=1}^4\left|(E_i,v)_{H^{\rho_0}}\right|
	\leq\, &2\sum_{i=1}^4\|E_i\|_{H^{\rho_0}}\|v\|_{H^{\rho_0}}\\
	\lesssim\,&\sum_{i=1}^4\|E_i\|_{H^{\rho_0}}^2+\|v\|_{H^{\rho_0}}^2\lesssim  C(T_l)n^{2r_s}+\|v\|_{H^{\rho_0}}^2.
\end{align*}
Using Lemma \ref{Taylor} and integration by parts, we obtain that
\begin{align*}
	\,&\left|(D^{\rho_0}\partial_x(Pv),D^{\rho_0}v)_{L^2}\right|\\
	=\,  & \left|([D^{\rho_0}\partial_x,P]v,D^{\rho_0}v)_{L^2}
	+(PD^{\rho_0}\partial_xv,D^{\rho_0}v)_{L^2}\right|\\
	\lesssim \, &  \|P\|_{H^s}\|v\|_{H^{\rho_0}}^2+\|P_x\|_{L^{\infty}}\|v\|_{H^{\rho_0}}^2 
	\lesssim(\|u^{m,n}\|_{H^s}+\|u_{m,n}\|_{H^s})^k\|v\|_{H^{\rho_0}}^2.
\end{align*}
Then, we use Lemma \ref{F lemma} to find that 
\begin{align*}
	\left|\big(F(u_{m,n})-F(u^{m,n}),v\big)_{H^{\rho_0}}\right|
	\lesssim \, &  \|F(u_{m,n})-F(u^{m,n})\|_{H^{\rho_0}}\|v\|_{H^{\rho_0}}\\
	\lesssim \, &  \|F(u_{m,n})-F(u^{m,n})\|_{H^{\rho_0}}^2+\|v\|_{H^{\rho_0}}^2\\
	\lesssim \, &  (\|u^{m,n}\|_{H^{s}}+\|u_{m,n}\|_{H^{s}})^{2k}\|v\|_{H^{\rho_0}}^2+\|v\|_{H^{\rho_0}}^2,
\end{align*}
To sum up, by \eqref{tau actual solution mnR}, Lemmas \ref{lemma ul Hr} and \ref{lemma uh Hr}, we arrive at
\begin{align*}
	\E\sup_{t\in[0,T_l\wedge\tau^{m,n}_R]}\|v(t)\|_{H^{\rho_0}}^2
	\leq C(T_l)n^{2r_s}+C_R\int_0^{T_l}\E\sup_{t'\in[0,t\wedge\tau^{m,n}_R]}\|v(t')\|_{H^{\rho_0}}^2\dd t.
\end{align*}
Using the Gr\"{o}nwall inequality, we obtain \eqref{Error sigma}. 

Now we prove \eqref{Error 2s minues sigma}. Since $2s-\rho_0>s>\frac{5}{2}$ and $u^{m,n}$ is the unique solution to \eqref{SGCH m n}, similar to \eqref{E 2 u-e}, we can use \eqref{tau actual solution mnR} and \ref{Ass2-growth-Lip} to find for each fixed $n\in\N$ that
\begin{align*}
	&\,\E\sup_{t\in[0,T_l\wedge\tau^{m,n}_R]}\|u^{m,n}(t)\|_{H^{2s-\rho_0}}^2\\
	\leq\,& C\E\|u_{m,n}(0)\|_{H^{2s-\rho_0}}^2
	+C_R\int_0^{T_l}\E\sup_{t'\in[0,t\wedge\tau^{m,n}_R]}\|u^{m,n}(t')\|_{H^{2s-\rho_0}}^2\dd t.
\end{align*}
Using the Gr\"{o}nwall inequality and Lemmas \ref{lemma ul Hr} and \ref{lemma uh Hr}, we find a constant $C=C(R,T_l)$ such that for all $n\ge1$,
\begin{align*}
	\E\sup_{t\in[0,T_l\wedge\tau^{m,n}_R]}\|u^{m,n}(t)\|_{H^{2s-\rho_0}}^2
	\leq \,& C\E\|u_{m,n}(0)\|_{H^{2s-\rho_0}}^2\\
	\leq\,& C(n^{\frac\delta2-\frac1k}+n^{s-\rho_0})^2\leq Cn^{2s-2\rho_0}.
\end{align*}
Hence, by Lemmas \ref{lemma ul Hr} and \ref{lemma uh Hr} again, we arrive at 
\begin{align*}
	&\,\E\sup_{t\in[0,T_l\wedge\tau^{m,n}_R]}\|(u^{m,n}-u_{m,n})(t)\|_{H^{2s-\rho_0}}^2\\
	\leq\, & 2\E\sup_{t\in[0,T_l\wedge\tau^{m,n}_R]}\|u^{m,n}(t)\|_{H^{2s-\rho_0}}^2
	+2\E\sup_{t\in[0,T_l\wedge\tau^{m,n}_R]}\|u_{m,n}(t)\|_{H^{2s-\rho_0}}^2 \\
	\leq\, & Cn^{2s-2\rho_0},\ \ n\geq1.
\end{align*}
The proof is therefore completed.
\qed

\subsection{Finish the proof for Theorem \ref{Weak instability}}
To begin with, we observe the following property:
\begin{Lemma}\label{exiting time infty lemma}
	Let \ref{Ass2-growth-Lip} hold true. Suppose that for some $R_0\gg 1$, the $R_0$-exiting time of the zero solution to \eqref{SGCH problem} is strongly stable. Then we have
	\begin{equation}\label{tau actual solution mnR infty}
		\lim_{n\rightarrow\infty}\tau^{m,n}_{R_0}=\infty\ \pas
	\end{equation}
\end{Lemma}
\proof
By \ref{Ass2-growth-Lip}, the unique solution with zero initial data to \eqref{SGCH problem} is zero.
On the other hand, we notice that for all $s'<s$, $\lim_{n\rightarrow\infty}\|u_{m,n}(0)\|_{H^{s'}}=\lim_{n\rightarrow\infty}\|u_{m,n}(0)-0\|_{H^{s'}}=0$.   Since the $R_0$-exiting time of the zero solution is $\infty$, we see that \eqref{tau actual solution mnR infty} holds provided that the $R_0$-exiting time of the zero solution to \eqref{SGCH problem} is strongly stable.
\qed

\textbf{Proof for Theorem \ref{Weak instability}}
Our strategy is to show that  if  the $R_0$-exiting time is strongly stable at the zero solution  for some $R_0\gg1$, then $\{u^{-1,n}\}$ and $\{u^{1,n}\}$ (if $k$ is odd) or $\{u^{0,n}\}$ and $\{u^{1,n}\}$ (if $k$ is even) are two sequences of solutions such that \eqref{tau 1 2 n}, \eqref{sup u}, \eqref{same initail data} and \eqref{sup sin t} are satisfied.

For each $n>1$ and for fixed $R_0\gg 1$, Lemmas \ref{lemma ul Hr}, \ref{lemma uh Hr} and \eqref{tau actual solution mnR} give $\p\{\tau^{m,n}_{R_0}>0\}=1$, and  Lemma \ref{exiting time infty lemma} implies \eqref{tau 1 2 n}. Then, it follows from
\eqref{tau actual solution mnR} that $u^{m,n}\in C([0,\tau^{m,n}_{R_0}];H^s)$ $\pas$  and \eqref{sup u} holds true. 
Next, we check \eqref{same initail data}.  By interpolation, we have
\begin{align*}
	&\,\E\sup_{t\in[0,T_l\wedge\tau^{m,n}_{R_0}]}\|u_{m, n}-u^{m, n}\|_{H^{s}}\\
	\leq  \, &     C 
	\left( \E\sup_{t\in[0,T_l\wedge\tau^{m,n}_{R_0}]}\|u_{m, n}-u^{m, n}\|_{H^{\rho_0}}\right)^{\frac12}
	\left( \E\sup_{t\in[0,T_l\wedge\tau^{m,n}_{R_0}]}\|u_{m, n}-u^{m, n}\|_{H^{2s-\rho_0}}\right)^{\frac12}\\
	\leq  \, &     C 
	\left( \E\sup_{t\in[0,T_l\wedge\tau^{m,n}_{R_0}]}\|u_{m, n}-u^{m, n}\|^2_{H^{\rho_0}}\right)^{\frac14}
	\left( \E\sup_{t\in[0,T_l\wedge\tau^{m,n}_{R_0}]}\|u_{m, n}-u^{m, n}\|^2_{H^{2s-\rho_0}}\right)^{\frac14}.
\end{align*} 
For $T_l>0$, combining Lemma \ref{difference estimate lemma} and the above estimate yields
\begin{align*}
	\E\sup_{t\in[0,T_l\wedge\tau^{m,n}_{R_0}]}
	\|u_{m, n}-u^{m, n}\|_{H^{s}}
	\leq C(R_0,T_l) n^{\frac{1}{4}\cdot 2r_{s}}\cdot n^{\frac{1}{4}\cdot(2s-2\rho_0)}=C(R_0,T_l)n^{r'_s},
\end{align*}
where $r_s$ is defined by \eqref{delta} and
$
r'_s=r_s\cdot \frac{1}{2}+(s-\rho_0)\cdot\frac{1}{2}=\frac{k\delta-1}{2}<0.
$
Consequently, we can deduce that
\begin{align}
	\lim_{n\rightarrow\infty}
	\E\sup_{t\in[0,T_l\wedge\tau^{m,n}_{R_0}]}\|u_{m, n}-u^{m, n}\|_{H^{s}}
	=0.\label{u difference tends to 0}
\end{align}
When $k$ is odd,
\begin{align*} 
	\|u^{-1,n}(0)-u^{1,n}(0)\|_{H^{s}}
	=\,&\|u_{-1,n}(0)-u_{1,n}(0)\|_{H^{s}}\\
	=\,& 2\Big\|n^{-\frac{1}{k}}\tilde{\phi}\left(\frac{x}{n^{\delta}}\right)\Big\|_{H^{s}}\rightarrow 0, {~\rm as~}n\rightarrow \infty.
\end{align*}
When $k$ is even
\begin{align*} 
	\|u^{0,n}(0)-u^{1,n}(0)\|_{H^{s}}=\,&\|u_{0,n}(0)-u_{1,n}(0)\|_{H^{s}}\\
	=\,&\Big\|n^{-\frac{1}{k}}\tilde{\phi}\left(\frac{x}{n^{\delta}}\right)\Big\|_{H^{s}}\rightarrow 0, {~\rm as~}n\rightarrow \infty.
\end{align*}
The above two estimates imply that \eqref{same initail data} holds true.

Now we prove \eqref{sup sin t}. Let $T_l>0$ be given in Lemma \ref{lemma ul Hr}. When $k$ is odd, we use \eqref{u difference tends to 0} to derive
\begin{align*}
	&\,\liminf_{n\rightarrow \infty}
	\E\sup_{t\in[0,T_l\wedge\tau^{-1,n}_{R_0}\wedge\tau^{1,n}_{R_0}]}
	\|u^{-1,n}(t)-u^{1,n}(t)\|_{H^s}\\
	\gtrsim \,  & 
	\liminf_{n\rightarrow \infty}
	\E\sup_{t\in[0,T_l\wedge\tau^{-1,n}_{R_0}\wedge\tau^{1,n}_{R_0}]}
	\|u_{-1,n}(t)-u_{1,n}(t)\|_{H^s}\\
	&-\lim_{n\rightarrow \infty}
	\E\sup_{t\in[0,T_l\wedge\tau^{-1,n}_{R_0}\wedge\tau^{1,n}_{R_0}]}
	\|u_{-1,n}(t)-u^{-1,n}(t)\|_{H^s}\\
	&-\lim_{n\rightarrow \infty}\E\sup_{t\in[0,T_l\wedge\tau_{-1,n}^R\wedge\tau_{1,n}^R]}
	\|u_{1,n}(t)-u^{1,n}(t)\|_{H^s}\\
	\gtrsim \,  & \liminf_{n\rightarrow \infty}
	\E\sup_{t\in[0,T_l\wedge\tau^{-1,n}_{R_0}\wedge\tau^{1,n}_{R_0}]}
	\|u_{-1,n}(t)-u_{1,n}(t)\|_{H^s}.
\end{align*}
It follows from the construction of $u_{m,n}$, Fatou's lemma,  Lemmas \ref{lemma ul Hr}, \ref{lemma uh Hr} and \ref{exiting time infty lemma} that
\begin{align}
	&\, \liminf_{n\rightarrow \infty}
	\E\sup_{t\in[0,T_l\wedge\tau^{-1,n}_{R_0}\wedge\tau^{1,n}_{R_0}]}
	\|u_{-1,n}(t)-u_{1,n}(t)\|_{H^s}\notag\\
	=\,  & \liminf_{n\rightarrow \infty}
	\E\sup_{t\in[0,T_l\wedge\tau^{-1,n}_{R_0}\wedge\tau^{1,n}_{R_0}]}
	\Big\|-2n^{-\frac{\delta}{2}-s}\phi\left(\frac{x}{n^{\delta}}\right)\sin(nx)\sin(t)\notag\\
	&\hspace{4.2cm}+[u_{l,-1,n}(t)-u_{l,1,n}(t)]\Big\|_{H^s}\notag\\
	\gtrsim \,  & \liminf_{n\rightarrow \infty}\E\sup_{t\in[0,T_l\wedge\tau^{-1,n}_{R_0}\wedge\tau^{1,n}_{R_0}]}n^{-\frac{\delta}{2}-s}\Big\|\phi\left(\frac{x}{n^{\delta}}\right)\sin(nx)\Big\|_{H^s}|\sin t|-\liminf_{n\rightarrow \infty} n^{\frac{\delta}{2}-\frac{1}{k}}\notag\\
	\gtrsim \,  &  \sup_{t\in[0,T_l]}|\sin t|,\label{check nonuniform k is odd}
\end{align}
which is \eqref{sup sin t} in the case that $k$ is odd. 
When $k$ is even, one has
\begin{align*}
	\|u_{0,n}(t)-u_{1,n}(t)\|_{H^s}&=\Big\|-2n^{-\frac{\delta}{2}-s}\phi\left(\frac{x}{n^{\delta}}\right)\sin(nx-t/2)\sin(t/2)-u_{l,1,n}(t)\Big\|_{H^s}\notag\\
	&\gtrsim n^{-\frac{\delta}{2}-s}\Big\|\phi\left(\frac{x}{n^{\delta}}\right)\sin(nx-t/2)\Big\|_{H^s}|\sin(t/2)|-n^{\frac{\delta}{2}-\frac{1}{k}}.
\end{align*}
Similar to \eqref{check nonuniform k is odd}, we can also obtain
\eqref{sup sin t} in the case that $k$ is even. The proof is completed.
\qed

\subsection{Example}\label{Example:Weak instability} Now we give an example of noise structure satisfying Hypothesis {\rm\ref{Assumption-2}}. For simplicity, we consider the case that $h(t,u) \dd \W
=b(t,u) \, {\rm d}W$, where $W$ is a standard 1-D Brownian motion. Let $m\ge1$ and $f(\cdot)$ be a continuous and bounded function, then
$$b(t,u)= f(t) {\rm e}^{-\frac{1}{\|u\|_{H^{\rho_0}}}}u^m,$$
satisfies Hypothesis {\rm\ref{Assumption-2}}.

\section{Noise prevents blow up}\label{sect:noise vs blow-up}

\subsection{Proof for Theorem \ref{Non breaking}}Our approach is motivated by  \cite{Ren-Tang-Wang-2020-Arxiv,Brzezniak-etal-2005-PTRF}. 
Let $s>5/2$ and $u_0\in H^s$ be an $H^s$-valued $\mathcal{F}_0$-measurable random variable with $\E\|u_0\|^2_{H^s}<\infty$. 
With \ref{Ass3-cut-off} and \ref{Ass3-Lip} at hand,  one can follow the steps in the proof for Theorem \ref{Local-WP} to obtain a unique solution $u$ to \eqref{SGCH non blow up Eq} such that $u\in C([0,\tau^*);H^s)$ $\pas$ and
\begin{equation}\label{tau*-tau*}
	\textbf{1}_{\left\{\limsup_{t\rightarrow \tau^*}\|u(t)\|_{H^{s}}=\infty\right\}}=\textbf{1}_{\left\{\limsup_{t\rightarrow \tau^*}\|u(t)\|_{W^{1,\infty}}=\infty\right\}}\ \pas
\end{equation}
Here we remark that \ref{Ass3-Lip} is the condition of locally Lipschitz continuous in $H^{\sigma}$ with $\sigma>3/2$, hence uniqueness can only be considered for solution in $H^s$ with $s>5/2$. This is because, if two solutions to \eqref{SGCH non blow up Eq} belong to $H^s$, the difference between them can be only estimated in $H^{s'}$ for $s'\leq s-1$   (Recalling \eqref{Verify Cauchy 1}, $H^{s+1}$-norm appears).

Define
\begin{align*}
	\tau_{m}:=\inf\left\{t\geq0: \|u(t)\|_{H^{s-1}}\geq m\right\},\ \ m\ge1\ \ \text{and} \ \ \widetilde{\tau^*}:=\lim\limits_{m\rightarrow\infty}\tau_m.
\end{align*}
Due to \eqref{tau*-tau*}, we have $\tau_{m}<\widetilde{\tau^*}=\tau^*$ $\pas$ and hence we only need to show
\begin{equation}
	\widetilde{\tau^*}=\infty\ \ \pas.\label{hat tau*}
\end{equation}
For $V\in\mathcal{V}$, applying the It\^{o} formula to
$\|u(t)\|^2_{H^{s-1}}$ and then to $V(\|u\|^2_{H^{s-1}})$, we find
\begin{align*}
	\, {\rm d}V(\|u\|^2_{H^{s-1}})
	=\,  & 2V'(\|u\|^2_{H^{s-1}})\left(q(t,u), u\right)_{H^{s-1}}\, {\rm d}W\notag\\
	&+V'(\|u\|^2_{H^{s-1}})
	\left\{-2\left(u^ku_x,u\right)_{H^{s-1}}
	-2\left(F(u), u\right)_{H^{s-1}}\right\}\dd t\notag\\
	&+V'(\|u\|^2_{H^{s-1}})\|q(t,u)\|^2_{H^{s-1}}\dd t\\
	&+2V''(\|u\|^2_{H^{s-1}})\left|\left(q(t,u), u\right)_{H^{s-1}}\right|^2\dd t.
\end{align*}
Next, we recall $\tau_{m}<\widetilde{\tau^*}=\tau^*$ and $s-1>3/2$, take expectation and then use Hypothesis {\rm \ref{Assumption-q}} and Lemma \ref{uux+F u Hs inner product} to find that
\begin{align*}
	&\E V(\|u(t\wedge \tau_m)\|^2_{H^{s-1}})  \\
	=\,  & \E V(\|u_0\|^2_{H^{s-1}})+\E\int_0^{t\wedge \tau_m}V'(\|u\|^2_{H^{s-1}})
	\left\{-2\left(u^ku_x, u\right)_{H^{s-1}}
	-2\left(F(u), u\right)_{H^{s-1}}\right\}\dd t'\notag\\
	&+\E\int_0^{t\wedge \tau_m}V'(\|u\|^2_{H^{s-1}})\|q(t',u)\|^2_{H^{s-1}}
	\dd t'\\
	&+\E\int_0^{t\wedge \tau_m}
	2V''(\|u\|^2_{H^{s-1}})\left|\left(q(t',u),u\right)_{H^{s-1}}\right|^2\dd t'\\
	\leq\,  & \E V(\| u_0\|^2_{H^{s-1}})+
	\E\int_0^{t\wedge \tau_m}\mathcal{H}_{s-1}(t',u)\dd t'\notag\\
	\leq  \, &     \E V(\| u_0\|^2_{H^{s-1}})+N_1 t-\E\int_0^{t\wedge \tau_m}
	N_2\frac{\left\{V'(\|u\|^2_{H^{s-1}})\left|\left(q(t',u), u\right)_{H^{s-1}}\right|\right\}^2}{ 1+V(\|u\|^2_{H^{s-1}})}
	\dd t',
\end{align*}
where $\mathcal{H}_{\sigma}(t,u)$ ($u\in H^{\sigma}$ and $\sigma>3/2$) is defined in Hypothesis \ref{Ass3-growth}.
Then we can infer from the above estimate that there is a constant $C(u_0,N_1,N_2,t)>0$ such that
\begin{align}
	\E\int_0^{t\wedge \tau_m}
	\frac{\left\{V'(\|u\|^2_{H^{s-1}})\left|\left(q(t',u), u\right)_{H^{s-1}}\right|\right\}^2}
	{ 1+V(\|u\|^2_{H^{s-1}})}\dd t'
	\leq C(u_0,N_1,N_2,t).\label{to use BDG 2}
\end{align}
Next,  for any $T>0$, it follows from the BDG inequality that
\begin{align*}
	&\E\sup_{t\in[0,{T\wedge \tau_m}]}V(\|u\|^2_{H^{s-1}})- \E V(\|u_0\|^2_{H^{s-1}})\\
	\leq  \, &  
	C\E\left(\int_0^{T\wedge \tau_m}
	\left\{V'(\|u\|^2_{H^{s-1}})\left|\left(q(t,u), u\right)_{H^{s-1}}\right|\right\}^2
	\dd t\right)^\frac12\\
	&+N_1T+N_2\E\int_0^{T\wedge \tau_m}
	\frac{\left\{V'(\|u\|^2_{H^{s-1}})\left|\left(q(t,u), u\right)_{H^{s-1}}\right|\right\}^2}{ 1+V(\|u\|^2_{H^{s-1}})}\dd t\\
	\leq  \, &   
	\frac12\E\sup_{t\in[0,{T\wedge \tau_m}]}\left(1+V(\|u\|^2_{H^{s-1}})\right)
	+C\E\int_0^{T\wedge \tau_m}
	\frac{\left\{V'(\|u\|^2_{H^{s-1}})
		\left|\left(q(t,u), u\right)_{H^{s-1}}\right|\right\}^2}
	{ 1+V(\|u\|^2_{H^{s-1}})}
	\dd t\\
	&+N_1T+N_2\E\int_0^{T\wedge \tau_m}
	\frac{\left\{V'(\|u\|^2_{H^{s-1}})\left|\left(q(t,u), u\right)_{H^{s-1}}\right|\right\}^2}{ 1+V(\|u\|^2_{H^{s-1}})}\dd t.
\end{align*}
Thus we use  \eqref{to use BDG 2}  to obtain
\begin{align*}
	&\E\sup_{t\in[0,{T\wedge \tau_m}]}V(\| u\|^2_{H^{s-1}}) \\
	\leq  \, &     1+2\E V(\|u_0\|^2_{H^{s-1}})
	+C\E\int_0^{T\wedge \tau_m}
	\frac{\left\{V'(\|u\|^2_{H^{s-1}})
		\left|\left(q(t,u), u\right)_{H^{s-1}}\right|\right\}^2}
	{ 1+V(\|u\|^2_{H^{s-1}})}
	\dd t\\
	&+2N_1T+2N_2\E\int_0^{T\wedge \tau_m}\frac{\left\{V'(\|u\|^2_{H^{s-1}})\left|\left(q(t,u), u\right)_{H^{s-1}}\right|\right\}^2}{ 1+V(\|u\|^2_{H^{s-1}})}
	\dd t\\
	\leq  \, &    
	C(u_0,N_1,T)+C(N_2) \E\int_0^{T\wedge \tau_m}
	\frac{\left\{V'(\|u\|^2_{H^{s-1}})\left|\left(q(t,u), u\right)_{H^{s-1}}\right|\right\}^2}{ 1+V(\|u\|^2_{H^{s-1}})}
	\dd t\\
	\leq  \, &    C(u_0,N_1,N_2,T).
\end{align*}
As a result, for all $m\ge1$,
\begin{align*}
	\p\{\widetilde{\tau^*}<T\}\leq\, &\p\{\tau_{m}<T\}\\
	\leq\,&
	\p\left\{\sup_{t\in[0,T\wedge\tau_m]}V(\|u\|^2_{H^{s-1}})\geq V(m^2)\right\}
	\leq \,\frac{C(u_0,N_1,N_2,T)}{V(m^2)}.
\end{align*}
Since $\p\{\widetilde{\tau^*}<T\}$ does not depend on $m$,
sending $m\rightarrow\infty$ gives rise to
$\p\{\tau^*<T\}=0$. Since $T>0$ is arbitrary, we obtain \eqref{hat tau*}, which completes the proof for
Theorem \ref{Non breaking}.

\subsection{Example}\label{Example:large noise}

As in \eqref{Blow-up criterion}, for the solution to \eqref{SGCH problem}, its $H^s$-norm blows up if and only if its $W^{1,\infty}$-norm blows up. On the other hand, \ref{Ass3-growth} means that the growth of $2\lambda_s\|u\|^k_{\Wlip}\|u\|^2_{H^s}$ can be canceled by $2V''(\|u\|^2_{H^s})|(q(t,u), u)_{H^s}|$.
Motivated by these two observations, we consider the following examples where the ${W^{1,\infty}}$-norm of $u$ will be involved, that is,
\begin{equation}\label{sigma alpha u}
	q(t,u)=\beta(t,\|u\|_{\Wlip})u,
\end{equation}
where $\beta(t,x)$ satisfies the following conditions:

\begin{Hypothesis}\label{Assumption-beta}
	We assume that
	\begin{itemize}
		\item The function $\beta(t,x)\in C\left([0,\infty)\times [0,\infty)\right)$ such that for any $x\ge0$, $\beta(\cdot,x)$ is bounded as a function of $t$, and for all $t\ge0$,  $\beta(t,\cdot)$ is locally Lipschitz continuous as a function of $x$;

		\item The function $\beta(t,x)\neq0$ for all $(t,x)\in  [0,\infty)\times [0,\infty)$, and 
		$\limsup_{x\rightarrow +\infty}\frac{2\lambda_s x^k}{\beta^2(t,x)}<1$ for all $t\ge0$, 
		where $\lambda_s>0$ is given in Lemma \ref{uux+F u Hs inner product}.
	\end{itemize}
	
\end{Hypothesis}

Now we give a concrete example $\beta(t,x)$ satisfying Hypothesis \ref{Assumption-beta}. 
Let $b:[0,\infty)\to[0,\infty)$ be a continuous function and there are constants $b_*,b^*>0$ such that $b_*\leq b^2(t)\leq b^*<\infty$ for all $t$. For all $k\ge1$, if
$$\text{either}\ \theta>k/2,\  b^*>b_*>0\ \ {\rm or}\ \ \theta=k/2,\  b^*>b_*>2\lambda_s,$$
then
$\beta(t,x)=b(t)(1+x)^\theta$ 
satisfies Hypothesis \ref{Assumption-beta}. Moreover, by the following two lemmas, we will see that $q(t,u)=b(t)(1+\|u\|_{\Wlip})^\theta u$ satisfies Hypothesis {\rm \ref{Assumption-q}}.

\begin{Lemma}\label{beta lemma}
	Let $\lambda_s$ be given in Lemma \ref{uux+F u Hs inner product}. Let $K>0$. If
	Hypothesis \ref{Assumption-beta} holds true, then there is an $M_1>0$ such that for any $M_2>0$ and all $0<x\leq Ky<\infty$,
	\begin{equation}\label{a K1 K2}
		\frac{2\lambda_s x^ky^2+
			\beta^2(t,x)y^2}{1+y^2}
		-\frac{2\beta^2(t,x)y^4}{(1+y^2)^2}\leq M_1-M_2\frac{2\beta^2(t,x)y^4}{(1+y^2)^2(1+\log(1+y^2))}.
	\end{equation}
\end{Lemma}
\proof
By Hypothesis \ref{Assumption-beta},  we have
\begin{align*}
	&\limsup_{x\rightarrow +\infty}\frac{2\lambda_s x^ky^2+
		\beta^2(t,x)y^2}{1+y^2}
	-\frac{2\beta^2(t,x)y^4}{(1+y^2)^2}+M_2\frac{2\beta^2(t,x)y^4}{(1+y^2)^2(1+\log(1+y^2))}\\
	\leq  \, &     \limsup_{x\rightarrow +\infty}\left(\frac{2\lambda_s x^k}{\beta^2(t,x)}+
	1
	-\frac{2\left(\frac{x}{K}\right)^4}{\left(1+\left(\frac{x}{K}\right)^2\right)^2}
	+M_2\frac{2}{\left(1+\log\left(1+\left(\frac{x}{K}\right)^2\right)\right)}\right)\beta^2(t,x)<0,
\end{align*}
which implies  \eqref{a K1 K2}. 
\qed

\begin{Lemma}
	If $\beta(t,x)$ satisfies  Hypothesis \ref{Assumption-beta}, then  $q(t,u)$ defined by \eqref{sigma alpha u} satisfies
	Hypothesis {\rm \ref{Assumption-q}}.
\end{Lemma}

\proof 
It follows from Lemma \ref{beta lemma} that \ref{Ass3-growth} holds true with the choice $V(x)=\log(1+x)\in\mathcal{V}$. Since $H^s\hookrightarrow \Wlip$ with $s>3/2$, it is obvious that the other requirements in Hypothesis {\rm \ref{Assumption-q}} are verified. 
\qed

\appendix\section{Auxiliary results}\label{Section:Preliminaries}
In this appendix we formulate and prove some estimates employed in the above proofs.  We first recall the Friedrichs mollifier $J_{\e}$ defined as 
\begin{equation}\label{Friedrichs mollifier}
	[J_{\varepsilon}f](x)=[j_{\varepsilon}\star f](x),\ \ \e\in(0,1),
\end{equation}
where
$\star$ stands for the convolution, $j_{\varepsilon}(x)=\frac{1}{\e}j(\frac{x}{\e})$ and $j(x)$ is a Schwartz function satisfying $\widehat{j}(\xi):\R\to[0,1]$ and $\widehat{j}(\xi)=1$ for $\xi\in[-1,1]$. 
From the above construction, we have
\begin{Lemma}[\cite{Tang-2018-SIMA,Li-Liu-Tang-2021-SPA}]\label{mollifier properties}
	For all $\e\in(0,1)$, $s,r\in\R$ and $u\in H^s$, $J_\e$ constructed in \eqref{Friedrichs mollifier} satisfies
	\begin{align*}
		\|I-J_{\varepsilon}\|_{\LL(H^s;H^r)}\lesssim\, \varepsilon^{s-r},\ \ 
		\|u-J_\e u\|_{H^r}\sim\, o(\varepsilon^{s-r}),\ \ r\leq s,
	\end{align*}
	\begin{align*}
		\|J_\e\|_{\LL(H^s;H^r)}\sim\,& O(\e^{s-r}),\ \ r>s,
	\end{align*}
	and
	\begin{align*}
		[D^s,J_\e]=0,\ \ 
		(J_\e f, g)_{L^2}=(f, J_\e g)_{L^2},\ \ 
		\|J_\e\|_{\LL(L^\infty;L^\infty)}\lesssim1,\ \ 
		\|J_\e\|_{\LL(H^s;H^s)}\leq 1,
	\end{align*}
	where $\LL(\X;\Y)$ is the space of bounded linear operators from  $\X$ to $\Y$.
\end{Lemma}

\begin{Lemma}[\cite{Taylor-2011-PDEbook3}]\label{Te commutator} 
	Let  $f,g$ be two functions such that $g\in W^{1,\infty}$ and $f\in L^2$. Then for some $C>0$,
	\begin{align*}
		\|[J_{\varepsilon}, g]f_x\|_{L^2}
		\leq C\|g_x\|_{L^\infty}\|f\|_{L^2}.
	\end{align*}
\end{Lemma}

\begin{Lemma}[\cite{Kato-Ponce-1988-CPAM}]\label{Kato-Ponce commutator estimate}
	If $f\in H^s\bigcap W^{1,\infty},\ g\in H^{s-1}\bigcap L^{\infty}$ for $s>0$, then there exists a constant $C_s>0$ such that
	$$
	\left\|\left[D^s,f\right]g\right\|_{L^2}\leq C_s(\|D^sf\|_{L^2}\|g\|_{L^{\infty}}+\|\partial_xf\|_{L^{\infty}}\|D^{s-1}g\|_{L^2}).
	$$
	Besides, if $s>0$, then we have for all $f,g \in H^s\bigcap L^{\infty}$,
	$$\|fg\|_{H^s}\leq C_s(\|f\|_{H^s}\|g\|_{L^{\infty}}+\|f\|_{L^{\infty}}\|g\|_{H^s}).$$
\end{Lemma}

\begin{Lemma}[Proposition 4.2, \cite{Taylor-2003-PAMS}]\label{Taylor}
	Let $\rho>3/2$ and $0\leq \eta+1\leq \rho$.  We have for some $c>0$,
	$$\|[D^{\eta}\partial_x, f]v\|_{L^2}\leq c\|f\|_{H^{\rho}}\|v\|_{H^{\eta}}\ \   \forall\,  f\in H^{\rho}, v\in H^{\eta}.$$
\end{Lemma}

\begin{Lemma}\label{F lemma}
	For $F(\cdot)$ defined in \eqref{F decomposition}, we have for all $k\ge1$ the following estimates:
	\begin{align*}
		\|F(v)\|_{H^s}
		\lesssim \, &   \|v\|_{W^{1,\infty}}^k\|v\|_{H^s},
		\ \ s>3/2,\\
		\|F(v)\|_{H^s}
		\lesssim \, &    \|v\|^{k}_{\Wlip}\left(\|v\|_{H^s}+\|v_x\|_{H^{s}}\right),
		\ \ 0<s\leq3/2,\\
		\|F(u)-F(v)\|_{H^s}\lesssim \, &  \left(\|u\|_{H^s}+\|v\|_{H^s}\right)^k\|u-v\|_{H^s},
		\ \ s>3/2,\\
		\|F(u)-F(v)\|_{H^s}\lesssim \, &  \left(\|u\|_{H^{s+1}}+\|v\|_{H^{s+1}}\right)^k\|u-v\|_{H^s},
		\ \  1/2<s\leq3/2.
	\end{align*}
\end{Lemma}
\proof We only estimate $\|F(v)\|_{H^s}$ for $0<s\leq3/2$ since the other cases can be found in \cite{Rohde-Tang-2021-JDDE,Tang-Zhao-Liu-2014-AA,Tang-Shi-Liu-2015-MM}.
When $s>0$, by using \eqref{F decomposition} and Lemma \ref{Kato-Ponce commutator estimate}, we derive
\begin{align}
	\|F_1(v)\|_{H^s}\lesssim\|v^{k+1}\|_{H^s}\lesssim\|v\|_{L^\infty}^{k}\|v\|_{H^s},\ \ k\ge 1.\label{F1 s<3/2}
\end{align}
When $k\ge2$, we have
\begin{align*}
	\left\|F_2(v)\right\|_{H^{s}}\lesssim \, &  \left\|v^{k-1}v_x^2\right\|_{H^s}\\
	\lesssim \, &  \|v\|_{H^s}\|v\|_{L^\infty}^{k-2}\|v_x\|_{L^\infty}^2
	+\left\|v\right\|_{L^{\infty}}^{k-1}\|v_x\|_{H^{s}}\|v_x\|_{L^\infty}\\
	\lesssim \, &  \|v\|^{k}_{\Wlip}\left(\|v\|_{H^s}+\|v_x\|_{H^{s}}\right).
\end{align*}
When $k=1$, $F_2(v)=\frac{1}{2}(1-\partial_{x}^2)^{-1}\partial_x\left(v_x^2\right)$ and hence
\begin{align*}
	\left\|F_2(v)\right\|_{H^{s}}
	\lesssim\left\|v_x^2\right\|_{H^s}
	\lesssim
	\left\|\partial_xv\right\|_{L^{\infty}}\|v_x\|_{H^{s}}
	\lesssim\|v\|_{\Wlip}\left(\|v\|_{H^s}+\|v_x\|_{H^{s}}\right).
\end{align*}
Combining the above two cases for $F_2$, we arrive at
\begin{align}
	\left\|F_2(v)\right\|_{H^{s}}\lesssim\|v\|^k_{\Wlip}\left(\|v\|_{H^s}+\|v_x\|_{H^{s}}\right),\ \ k\ge 1.\label{F2 s<3/2}
\end{align}
Now we consider $F_3$. When $k\ge 3$, we have
\begin{align*}
	\left\|F_3(v)\right\|_{H^{s}}\lesssim \, &  \left\|v^{k-2}v_x^3\right\|_{H^s}\\
	\lesssim \, &  \left\|v\right\|_{H^{s}}\|v\|_{L^\infty}^{k-3}\|v_x\|_{L^\infty}^3
	+\left\|v\right\|_{L^\infty}^{k-2}\|v_x\|_{H^{s}}\|v_x\|_{L^\infty}^2\\
	\lesssim \, &  \|v\|^{k}_{\Wlip}\left(\|v\|_{H^s}+\|v_x\|_{H^{s}}\right).
\end{align*}
When $k=2$, we have $F_3(v)=\frac{1}{2}(1-\partial_{x}^2)^{-1}\left(v_x^3\right)$ and then
\begin{align*}
	\left\|F_3(v)\right\|_{H^{s}}\lesssim \, &  \left\|v_x^3\right\|_{H^s}
	\lesssim \|v_x\|_{H^{s}}\|v_x\|_{L^\infty}^2 
	\lesssim\|v\|^{2}_{\Wlip}\left(\|v\|_{H^s}+\|v_x\|_{H^{s}}\right).
\end{align*}
Combining the above two cases for $F_3$ with noticing that $F_3=0$ for $k=1$, we find
\begin{align}
	\left\|F_3(v)\right\|_{H^{s}}\lesssim\|v\|^k_{\Wlip}\left(\|v\|_{H^s}+\|v_x\|_{H^{s}}\right),\ \ k\ge 1.\label{F3 s<3/2}
\end{align}
Then the desired estimate is a consequence of \eqref{F1 s<3/2}, \eqref{F2 s<3/2} and \eqref{F3 s<3/2}.
\qed

\begin{Lemma}\label{uux+F u Hs inner product}
	Let $s>3/2$, $k\ge1$, $F(\cdot)$ be given in \eqref{F decomposition}  and $J_\e$ be the mollifier defined in \eqref{Friedrichs mollifier}. There exists a constant $\lambda_{s}>0$ such that for all $\e>0$,
	\begin{align*}
		\left|\left(D^sJ_\e 
		\left[u^ku_x\right],D^sJ_\e u\right)_{L^2}\right|+\left|\left(D^sJ_\e F(u), D^sJ_\e u\right)_{L^2}\right|
		\leq  \, &      \lambda_{s}\|u\|^{k}_{W^{1,\infty}}\|u\|^2_{H^s},\ \ u\in H^s,\ \ s>3/2.
	\end{align*}
	If $u\in H^{s+1}$, then $u^ku_x\in H^{s}$, and the above estimate also holds true without $J_\e$.
\end{Lemma}

\proof
We only prove the case that $u\in H^s$. 
It follows from Lemmas  \ref{mollifier properties},  \ref{Te commutator} and \ref{Kato-Ponce commutator estimate}, integration by parts and $H^s\hookrightarrow W^{1,\infty}$ that
\begin{align*}
	&\big|\left(D^sJ_\e 
	\left[u^ku_x\right],D^sJ_\e u\right)_{L^2}\big|\notag\\
	\leq\,  & 
	\big|\left(\left[D^s,
	u^k\right]u_x,D^sJ^2_\e u\right)_{L^2}\big|+
	\big|\left([J_\e,u^k]D^su_x, D^sJ_\e  u\right)_{L^2}\big|
	+\big|\left(u^kD^sJ_\e u_x, D^sJ_\e u\right)_{L^2}\big|\notag\\
	\leq  \, &      C(s)\|u\|^{k}_{W^{1,\infty}}\|u\|^2_{H^s}.
\end{align*}
From Lemma \ref{F lemma}, we also have
\begin{align*}
	\big|\left(D^sJ_\e F(u), D^sJ_\e u\right)_{L^2}\big|
	\leq  C(s)\|u\|^{k}_{W^{1,\infty}}\|u\|^2_{H^s}.
\end{align*}
Combining the above two inequalities gives rise to the desired estimate of the lemma.
\qed

The following technique has been used in \cite{GlattHoltz-Ziane-2009-ADE,Breit-Feireisl-Hofmanova-2018-Book,Alonso-Rohde-Tang-2021-JNLS}. Here we formulate such a technique result in an abstract way.
\begin{Lemma}\label{cut-combine}
	Suppose $u_0$ is an $H^s$-valued $\mathcal{F}_0$-measurable random variable, and suppose \ref{Ass1-growth assumption} holds true. Let $I$ be a   countable index set and let $\{{\rm \Omega}_i\}_{i\in I}$ satisfy
	\begin{equation}
		{\rm \Omega}_i\subset {\rm \Omega},\ \p\left\{\cup_{i\in I} {\rm \Omega}_i\right\}=1 \ \text{and}\ \ 
		{\rm \Omega}_i\cap {\rm \Omega}_j=\emptyset\ \text{for all}\ i,j\in I,\ i\neq j.\label{Omega i condition}
	\end{equation}
	If $(u_i,\tau_i)$ with $i\in I$ is a solution  to \eqref{SGCH problem} with initial value $\textbf{1}_{{\rm \Omega}_i}u_0$, then 
	\begin{equation}
		\left(u=\sum_{i\in I}\textbf{1}_{{\rm \Omega}_i}u_i,\ \
		\tau=\sum_{i\in I}\textbf{1}_{{\rm \Omega}_i}\tau_i\right)\label{combine-solution}
	\end{equation}
	is a solution to \eqref{SGCH problem} with initial data $u_0$.
\end{Lemma}
\proof
Since $(u_i,\tau_i)$ is a solution  to \eqref{SGCH problem} with initial value $u_0\textbf{1}_{{\rm \Omega}_i}$, we find
\begin{align*}
	u_i(t\wedge \tau_i)-\textbf{1}_{{\rm \Omega}_i}u_{0}
	=\,  & -\int_0^{t\wedge \tau_i}\left[u_i^k\partial_x u_i+F(u_i)\right] \dd t'
	+\int_0^{t\wedge \tau_i}h(t,u_i) \dd \W\ \ \pas
\end{align*}
Therefore, we restrict the above equation to ${\rm \Omega}_i$ and we  obtain
\begin{align*}
	\textbf{1}_{{\rm \Omega}_i}u_i(t\wedge \tau_i)-\textbf{1}_{{\rm \Omega}_i}u_{0}
	=\,  & -\int_0^{t\wedge \textbf{1}_{{\rm \Omega}_i}\tau_i}\textbf{1}_{{\rm \Omega}_i}\left[u_i^k\partial_x u_i+F(u_i)\right] \dd t'
	+\int_0^{t\wedge \textbf{1}_{{\rm \Omega}_i}\tau_i}\textbf{1}_{{\rm \Omega}_i} h(t,u_i) \dd \W\ \ \pas
\end{align*}
It is clear that almost surely, 
$$\textbf{1}_{{\rm \Omega}_i} h(t,u_i)=h(t,\textbf{1}_{{\rm \Omega}_i} u_i)-\textbf{1}_{{\rm \Omega}^C_i}h(t,\mathbf{0}),\ \textbf{1}_{{\rm \Omega}_i} \left[u^k_i\partial_x u_i+F(u_i)\right]=\left[(\textbf{1}_{{\rm \Omega}_i} u_i)^k\partial_x \left(\textbf{1}_{{\rm \Omega}_i} u_i\right)+F(\textbf{1}_{{\rm \Omega}_i}u_i)\right].$$
By \ref{Ass1-growth assumption}, we have $\|h(t,\mathbf{0})\|_{\LL_2(\U;H^s)}<\infty$. Then, from the above three equations, we have that almost surely
\begin{align*}
	\textbf{1}_{{\rm \Omega}_i}u_i(t\wedge \tau_i)-\textbf{1}_{{\rm \Omega}_i}u_{0} 
	=\,  & \textbf{1}_{{\rm \Omega}_i}u_i(t\wedge \textbf{1}_{{\rm \Omega}_i}\tau_i)-\textbf{1}_{{\rm \Omega}_i}u_{0}\notag\\
	=\,  & -\int_0^{t\wedge \textbf{1}_{{\rm \Omega}_i}\tau_i}\left[(\textbf{1}_{{\rm \Omega}_i} u_i)^k\partial_x (\textbf{1}_{{\rm \Omega}_i} u_i)+F(\textbf{1}_{{\rm \Omega}_i}u_i)\right]\dd t'\\
	&+\int_0^{t\wedge \textbf{1}_{{\rm \Omega}_i}\tau_i}h(t,\textbf{1}_{{\rm \Omega}_i} u_i) \dd \W,
\end{align*}
which means $(\textbf{1}_{{\rm \Omega}_i}u_i,\textbf{1}_{{\rm \Omega}_i}\tau_i)$  also solves \eqref{SGCH problem} with initial data $\textbf{1}_{{\rm \Omega}_i}u_{0}$. By summing up both sides of the above equation with noticing \eqref{Omega i condition}, we derive   that
\eqref{combine-solution} 
is a  solution to \eqref{SGCH problem} with initial data $u_0$ almost surely. Indeed, for the initial data, we have  $u_0=\sum_{i\in I}\textbf{1}_{{\rm \Omega}_i}u_0$ $\pas$  For the nonlinear  term $u^k\partial_x u$, by \eqref{Omega i condition}, we have that $\pas$,
\begin{align*}
	\,&\sum_{i\in I}\int_0^{t\wedge\textbf{1}_{{\rm \Omega}_i}\tau_i} (\textbf{1}_{{\rm \Omega}_i}u_i)^k\partial_x(\textbf{1}_{{\rm \Omega}_i}u_i) \dd t'\\
	=\,  & \sum_{i\in I}\int_0^{t\wedge\sum_{j\in I}\textbf{1}_{{\rm \Omega}_j}\tau_j} (\textbf{1}_{{\rm \Omega}_i}u_i)^k\partial_x\left(\sum_{l\in I}\textbf{1}_{{\rm \Omega}_l}u_{l}\right) \dd t'\notag\\
	=\,  & \int_0^{t\wedge\tau} \sum_{i\in I}(\textbf{1}_{{\rm \Omega}_i}u_i)^k\partial_x\left(\sum_{l\in I}\textbf{1}_{{\rm \Omega}_l}u_{l}\right) \dd t'\notag\\
	=\,  & \int_0^{t\wedge\tau} \left(\sum_{i\in I}\textbf{1}_{{\rm \Omega}_i}u_i\right)^k
	\partial_x\left(\sum_{l\in I}\textbf{1}_{{\rm \Omega}_l}u_{l}\right)\dd t' 
	=\,    \int^{t\wedge\tau}_0 u^k\partial_x u\dd t'.
\end{align*}
The other terms can also be  justified in the same way, here we omit the details. 
\qed

Finally, we recall the following estimate on the product of a  Schwartz function and a trigonometric function.
\begin{Lemma}[\cite{Himonas-Kenig-2009-DIE,Koch-Tzvetkov-2005-IMRN}]\label{lemma uh Hr}
	Let $\mathscr{S}(\R)$ be the set of Schwartz functions. 
	Let $\delta>0$ and $\alpha\in\R$. Then for any $r\geq 0$ and $\psi\in\mathscr{S}(\R)$, we have that
	\begin{equation}\label{uh Hr estimate}
		\lim_{n\rightarrow\infty}n^{-\frac{\delta}{2}-r}\left\|\psi\left(\frac{x}{n^{\delta}}\right)\cos(nx-\alpha)\right\|_{H^r}=\frac{1}{\sqrt{2}}\|\psi\|_{L^2}.
	\end{equation}
	Relation \eqref{uh Hr estimate} is also true if $\cos$ is replaced by $\sin$.
\end{Lemma}

\section{Proof for Lemma 
	\ref{Lemma Error estimate}}\label{sec:appendix}

As $u_{m,n}=u_l+u_h$ is explictly given, we will firstly estimate $E_i$ ($i=1,2,3,4$). Let $T_l>0$ be given in Lemma  \ref{lemma ul Hr} such that $u_l$ exits on $[0,T_l]$ for all $n\gg 1$ and \eqref{ul Hr estimate} is satisfied.

\textbf{ (i) Estimating  $\|E_1\|_{H^{\rho_0}}$.} We apply the embedding $H^{\rho_0}\hookrightarrow L^\infty$, Lemmas \ref{lemma ul Hr} and \ref{lemma uh Hr} to obtain
\begin{align}
	&\,\|E_1\|_{H^{\rho_0}}\notag\\
	\leq\,
	&\left\|\left[u_{l}^k(0)-u_{l}^k(t)\right]n^{1-\frac{\delta}{2}-s}\phi\left(\frac{x}{n^{\delta}}\right)\sin(nx-mt)+u_{l}^k(t)n^{-\frac{3\delta}{2}-s}\partial_x\phi\left(\frac{x}{n^{\delta}}\right)\cos(nx-mt)\right\|_{H^{\rho_0}}\notag\\
	&+\left\|\ZZ_k\partial_xu_l\right\|_{H^{\rho_0}}+\left\|\ZZ_k\partial_xu_h\right\|_{H^{\rho_0}}\notag\\
	\lesssim \, &   n^{1-\frac{\delta}{2}-s}\left\|u_l^k(0)-u_l^k(t)\right\|_{H^{\rho_0}}\left\|\phi\left(\frac{x}{n^{\delta}}\right)\sin(nx-mt)\right\|_{H^{\rho_0}}\notag\\
	&+n^{-\frac{3\delta}{2}-s}\left\|u_l(t)\right\|^k_{H^{\rho_0}}
	\left\|\partial_x\phi\left(\frac{x}{n^{\delta}}\right)\cos(nx-mt')\right\|_{H^{\rho_0}}+\left\|\ZZ_k\partial_xu_l\right\|_{H^{\rho_0}}+\left\|\ZZ_k\partial_xu_h\right\|_{H^{\rho_0}}\notag\\
	\lesssim \, &  n^{1-s+\rho_0}\left\|u_l^k(0)-u_l^k(t)\right\|_{H^{\rho_0}}
	+n^{-s-1+\rho_0+(\frac{k}{2}-1)\delta}
	+\left\|\ZZ_k\partial_xu_l\right\|_{H^{\rho_0}}+\left\|\ZZ_k\partial_xu_h\right\|_{H^{\rho_0}}\notag\\
	\lesssim \, &  n^{1-s+\rho_0}\left\|u_l^k(0)-u_l^k(t)\right\|_{H^{\rho_0}}
	+n^{r_s}
	+\left\|\ZZ_k\partial_xu_l\right\|_{H^{\rho_0}}+\left\|\ZZ_k\partial_xu_h\right\|_{H^{\rho_0}}, \ \ t\in[0,T_l].\label{Error E1 estimate}
\end{align}
Next, we estimate $\|u_l^k(0)-u_l^k(t)\|_{H^{\rho_0}}$. Using the fundamental theorem of calculus and the algebra property, we have that for all $k\ge 1$ and $t\in[0,T_l]$,
\begin{align*}
	\left\|u_l^k(0)-u_l^k(t)\right\|_{H^{\rho_0}}=\left\|k\int_0^tu_l^{k-1}(t')\partial_tu_l(t')\dd t'\right\|_{H^{\rho_0}}\lesssim\int_0^t\|u_l(t')\|^{k-1}_{H^{\rho_0}}\|\partial_tu_l(t')\|_{H^{\rho_0}}\dd t'.
\end{align*}
Using \eqref{low-frequency equation} with $t\in[0,T_l]$, \eqref{F decomposition}, Lemmas \ref{F lemma} and \ref{lemma ul Hr} and the embedding $H^{\rho_0+1}\hookrightarrow \Wlip$, we get
\begin{align*}
	\left\|u_l^k(0)-u_l^k(t)\right\|_{H^{\rho_0}}
	\lesssim \, &  \int_0^t\|u_l\|^{k-1}_{H^{\rho_0}}
	\left(\|u_l^{k}\partial_xu_l\|_{H^{\rho_0}}+\|F(u_l)\|_{H^{\rho_0+1}}\right)\dd t'\notag\\
	\lesssim \, &  \int_0^t\|u_l\|^{k-1}_{H^{\rho_0+1}}
	\|u_l\|_{H^{\rho_0+1}}^{k+1}\dd t'
	\lesssim \,  n^{(\frac\delta2-\frac1k)2k}T_l,\ \ k\geq 1,\ \ t\in[0,T_l],
\end{align*}
which implies
\begin{align}
	n^{1-s+\rho_0}\left\|u_l^k(0)-u_l^k(t)\right\|_{H^{\rho_0}}\lesssim n^{1-s+\rho_0+k\delta-2}T_l=n^{r_s}T_l,\ \ k\geq 1,\ \ t\in[0,T_l].\label{E1 part 1}
\end{align}
Again, applying the algebra property and using Lemmas \ref{lemma ul Hr}, \ref{lemma uh Hr} and \eqref{Z nonlinear}, we have that for all $k\geq 1$ and $t\in[0,T_l]$,
\begin{align}
	\left\|\ZZ_k\partial_xu_l\right\|_{H^{\rho_0}}
	\lesssim \, &  \sum_{j=1}^k\|u_l\|_{H^{\rho_0}}^{k-j}\|u_h\|_{H^{\rho_0}}^j\|u_l\|_{H^{\rho_0+1}}\notag\\
	\lesssim \, &  \left(\sum_{j=1}^kn^{(\frac\delta2-\frac1k)(k-j)}n^{(-s+\rho_0)j}\right)n^{\frac\delta2-\frac1k}
	\notag\\
	=\,  & \sum_{j=1}^k n^{j(-s+\rho_0-\frac\delta2+\frac1k)-1-\frac1k+(k+1)\frac\delta2} 
	\lesssim\   n^{-s-1+\rho_0+k\frac\delta2} \lesssim n^{r_s}. \label{E1 part 3}
\end{align}
Here  we used the facts that  $-s+\rho_0-\frac\delta2+\frac1k<-s+1-\frac\delta2+1=-s+2-\frac\delta2<-\frac{1}{2}-\frac\delta2<0$ for all $k\ge1$, which means that the term corresponding to $j=1$ dominates.

For the last term $\ZZ_k\partial_xu_h$, by using \eqref{high-frequency approximation solutions}, \eqref{Z nonlinear}, Lemma \ref{Kato-Ponce commutator estimate}, Lemmas \ref{lemma ul Hr} and \ref{lemma uh Hr}, we obtain that 
\begin{align*} 
	&\left\|\ZZ_k\partial_xu_h\right\|_{H^{\rho_0}}\\
	\lesssim \, &  \sum_{j=1}^k\|u_l\|_{H^{\rho_0}}^{k-j}\|u_h^j\partial_xu_h\|_{H^{\rho_0}}\notag\\
	\lesssim \, &  \sum_{j=1}^k\|u_l\|_{H^{\rho_0}}^{k-j}\Big[\|u_h^j\|_{H^{\rho_0}}\|\partial_xu_h\|_{L^{\infty}}+\|u_h^j\|_{L^{\infty}}\|u_h\|_{H^{\rho_0+1}}\Big]\notag\\
	\lesssim \, &  \sum_{j=1}^kn^{(\frac\delta2-\frac1k)(k-j)}
	\Bigg\{n^{(-s+\rho_0)j}\bigg[\Big\|-n^{1-\frac\delta2-s}\phi\left(\frac{x}{n^{\delta}}\right)\sin(nx-mt)\Big\|_{L^{\infty}}\notag\\
	&\hspace*{4.1cm}+\Big\|n^{-\frac{3\delta}{2}-s}\partial_x\phi\left(\frac{x}{n^{\delta}}\right)\cos(nx-mt)\Big\|_{L^{\infty}}\bigg]\notag\\
	&\hspace*{2.7cm}+\Big\|n^{-\frac\delta2-s}\phi\left(\frac{x}{n^{\delta}}\right)\cos(nx-mt)\Big\|_{L^{\infty}}^jn^{-s+\rho_0+1}\Bigg\}\notag\\
	\lesssim \, &  \sum_{j=1}^kn^{(\frac\delta2-\frac1k)(k-j)}\left\{n^{(-s+\rho_0)j}\left[n^{-s-\frac\delta2+1} +n^{-\frac{3\delta}{2}-s}\right]+n^{(-\frac\delta2-s)j}n^{-s+\rho_0+1}\right\}\notag\\
	\lesssim \, &  
	\sum_{j=1}^kn^{j(-s+\rho_0-\frac\delta2+\frac1k)-s+(k-1)\frac\delta2}+\sum_{j=1}^kn^{j(-s-\delta+\frac1k)-s+\rho_0+k\frac\delta2},\ \ k\geq 1,\ \ t\in[0,T_l].
\end{align*}
When $k\ge1$, $-s+\rho_0-\frac\delta2+\frac1k<0$ and $-s-\delta+\frac1k<0$, therefore, both sums are bounded by $n^{-2s+\rho_0+\frac1k+(k-2)\frac\delta2}$. Furthermore, when $k\ge 1$, $-2s+\rho_0+\frac1k+(k-2)\frac\delta2\leq r_s$, which means
\begin{align} 
	\left\|\ZZ_k\partial_xu_h\right\|_{H^{\rho_0}}
	\lesssim \, &   n^{-2s+\rho_0+\frac1k+(k-2)\frac\delta2}\lesssim n^{r_s}, \ \ k\ge1,\ \ t\in[0,T_l].\label{E1 part 4}
\end{align}
Finally, inserting \eqref{E1 part 1}, \eqref{E1 part 3} and \eqref{E1 part 4} into \eqref{Error E1 estimate}, we arrive at 
\begin{align}\label{error E1 estimate}
	\|E_1\|_{H^{\rho_0}}\lesssim n^{r_s}, \ \ k\ge1,\ \ t\in[0,T_l].
\end{align}

\textbf{(ii) Estimating $\|E_2\|_{H^{\rho_0}}$.}  For $E_2$, we first recall \eqref{Z nonlinear}.
Applying the embedding $H^{\rho_0}\hookrightarrow L^\infty$ and Lemma \ref{lemma ul Hr}, and then taking the dominated term $j=1$, we find that for all $k\geq 1$ and  $t\in[0,T_l]$,
\begin{align}
	\|E_2\|_{H^{\rho_0}}\lesssim \, &  \left\|\sum_{j=1}^{k+1}C_{k+1}^{j}u_l^{k+1-j}u_h^j\right\|_{H^{\rho_0}} \notag\\
	\lesssim \,  & \sum_{j=1}^{k+1}n^{(k+1-j)(\frac\delta2-\frac1k)}n^{(-s+\rho_0)j} 
	\lesssim \,   n^{-s-1+\rho_0+k\frac\delta2}\lesssim n^{r_s}.\label{error E2 estimate}
\end{align}

\textbf{(iii) Estimating $\|E_3\|_{H^{\rho_0}}$.} As in the estimate for \eqref{E1 part 4}, we have obtained that $$\|\partial_xu_h\|_{L^{\infty}}\lesssim n^{-s-\frac\delta2+1} +n^{-\frac{3\delta}{2}-s}\lesssim n^{-s-\frac\delta2+1},\  \ t\in[0,T_l].$$
When $k=1$,  $Z_{k-1}=0$ and then we find
\begin{align*}
	\|E_3\|_{H^{\rho_0}}
	\lesssim \, &  \|2(\partial_xu_l)(\partial_xu_h)+(\partial_xu_h)^2\|_{H^{\rho_0-1}}\notag\\
	\lesssim \, &  \|2(\partial_xu_l)(\partial_xu_h)+(\partial_xu_h)^2\|_{L^{2}}\notag\\
	\lesssim \, &  \|u_l\|_{H^2}\|\partial_xu_h\|_{L^{\infty}}
	+\|\partial_x u_h\|_{L^2}\|\partial_xu_h\|_{L^{\infty}}\notag\\
	\lesssim \, &  \|u_l\|_{H^2}\|\partial_xu_h\|_{L^{\infty}}
	+\|u_h\|_{H^1}\|\partial_xu_h\|_{L^{\infty}}\notag\\
	\lesssim \, &  n^{\frac{\delta}{2}-1}n^{1-\frac{\delta}{2}-s}
	+n^{-s+1}n^{1-\frac{\delta}{2}-s} 
	\lesssim \,  n^{-s}+n^{-2s+2-\frac{\delta}{2}},\ \ t\in[0,T_l].
\end{align*}
Since $\delta>0$, 
$-2s+2-\frac{\delta}{2}-r_s=-s+3-\rho_0-\frac{3}{2}\delta<-\frac{5}{2}+3-\frac12-\frac{3}{2}\delta=-\frac{3}{2}\delta<0$, hence
\begin{align}
	\|E_3\|_{H^{\rho_0}}
	\lesssim \, &  n^{-s}+n^{-2s+2-\frac{\delta}{2}}
	\lesssim n^{r_s},\ \ k=1,\ \ t\in[0,T_l].\label{error E3 estimate k=1}
\end{align}
When $k\geq 2$, we can use the above estimate, Lemma \ref{lemma ul Hr}, the facts $\|f\|_{H^{\rho_0-1}}\leq\|f\|_{L^2}$ and $\|fg\|_{L^2}\leq\|f\|_{L^2}\|g\|_{L^{\infty}}$ and take the dominate term $j=1$ to obtain
\begin{align}
	\,&\|E_3\|_{H^{\rho_0}}\notag\\
	\lesssim \, &  \left\|u_l^{k-1}[2(\partial_xu_l)(\partial_xu_h)+(\partial_xu_h)^2]\right\|_{H^{\rho_0-1}}
	+\left\|\ZZ_{k-1}(\partial_xu_l+\partial_xu_h)^2\right\|_{H^{\rho_0-1}}\notag\\
	\lesssim \, &  \left\|u_l^{k-1}[2(\partial_xu_l)(\partial_xu_h)+(\partial_xu_h)^2]\right\|_{L^{2}}+\left\|\left(\sum_{j=1}^{k-1}C_{k-1}^ju_l^{k-1-j}u_h^j\right)(\partial_xu_l+\partial_xu_h)^2\right\|_{L^{2}}\notag\\
	\lesssim \, &  \|u_l\|_{H^2}^{k}\|\partial_xu_h\|_{L^{\infty}}+\|u_l\|_{H^2}^{k-1}\|\partial_xu_h\|_{L^{\infty}}^2+\sum_{j=1}^{k-1}\|u_l\|_{H^2}^{k+1-j}\|u_h\|_{L^{\infty}}^j\notag\\
	&+\sum_{j=1}^{k-1}\|u_l\|_{H^2}^{k-j}\|u_h^j\partial_xu_h\|_{L^{\infty}} 
	+\sum_{j=1}^{k-1}\|u_l\|_{H^2}^{k-1-j}\|u_h^j(\partial_xu_h)^2\|_{L^{\infty}}\notag\\
	\lesssim \, &  n^{k(\frac{\delta}{2}-\frac{1}{k})}n^{1-\frac{\delta}{2}-s}+n^{(k-1)(\frac{\delta}{2}-\frac{1}{k})}n^{2-\delta-2s}+\sum_{j=1}^{k-1}n^{(k+1-j)(\frac{\delta}{2}-\frac{1}{k})}n^{j(-\frac{\delta}{2}-s)}\notag\\
	&+\sum_{j=1}^{k-1}n^{(k-j)(\frac{\delta}{2}-\frac{1}{k})}n^{j(-\frac{\delta}{2}-s)+1-\frac{\delta}{2}-s}+\sum_{j=1}^{k-1}n^{(k-1-j)(\frac{\delta}{2}-\frac{1}{k})}n^{j(-\frac{\delta}{2}-s)+2-\delta-2s}\notag\\
	=\,  & n^{-s+(k-1)\frac{\delta}{2}}+n^{-2s+1+\frac{1}{k}+(k-3)\frac{\delta}{2}}+\sum_{j=1}^{k-1}n^{j(-s-\delta+\frac{1}{k})-1-\frac{1}{k}+(k+1)\frac{\delta}{2}}\notag\\
	&+\sum_{j=1}^{k-1}n^{j(-s-\delta+\frac{1}{k})-s+(k-1)\frac{\delta}{2}}+\sum_{j=1}^{k-1}n^{j(-s-\delta+\frac{1}{k})-2s+1+\frac{1}{k}+(k-3)\frac{\delta}{2}}\notag\\
	\lesssim \, &  n^{r_s},\ \ k\ge 2,\ \ t\in[0,T_l].\label{error E3 estimate k geq 2}
\end{align}
Combining \eqref{error E3 estimate k geq 2} and \eqref{error E3 estimate k=1}, we  have the following conclusion for $E_3$:
\begin{align}
	\|E_3\|_{H^{\rho_0}}\lesssim n^{r_s},\ \ k\ge 1,\ \ t\in[0,T_l].\label{error E3 estimate}
\end{align}

\textbf{(iv) Estimating  $\|E_4\|_{H^{\rho_0}}$.} For $k=1$, $E_4=0$ since $F_3$ disappears. When $k=2$, $\ZZ_{k-2}=0$ and then
\begin{align*}
	\,&\|E_4\|_{H^{\rho_0}}\\
	\lesssim \, &  \left\|3(\partial_xu_l)^2(\partial_xu_h)+3(\partial_xu_l)(\partial_xu_h)^2+(\partial_xu_h)^3\right\|_{H^{\rho_0-2}}\\
	\lesssim \, &  \|u_l\|_{H^{2}}^{2}\|\partial_xu_h\|_{L^{\infty}}+\|u_l\|_{H^{2}}\|\partial_xu_h\|_{L^{\infty}}^2+\|u_h\|_{H^1}\|\partial_xu_h\|_{L^{\infty}}^2\\
	\lesssim \, &  n^{2(\frac{\delta}{2}-\frac12)}n^{1-\frac{\delta}{2}-s}+n^{\frac{\delta}{2}-\frac12}n^{2(1-\frac{\delta}{2}-s)}+n^{1-s}n^{2(1-\frac{\delta}{2}-s)}\\
	=\,  & n^{-s+\frac{\delta}{2}}+n^{-2s-\frac{\delta}{2}+\frac32}+n^{3-\delta-3s}
	\lesssim \,     n^{r_s},\ \ k=2,\ \ t\in[0,T_l].
\end{align*}
Finally, for $k\geq 3$ and $t\in[0,T_l]$,
\begin{align*}
	&\,\|E_4\|_{H^{\rho_0}}\\
	\lesssim \, &  \left\|u_l^{k-2}[3(\partial_xu_l)^2(\partial_xu_h)+3(\partial_xu_l)(\partial_xu_h)^2+(\partial_xu_h)^3]\right\|_{H^{\rho_0-2}}+\left\|\ZZ_{k-2}(\partial_xu_l+\partial_xu_h)^3\right\|_{H^{\rho_0-2}}\\
	\lesssim \, &  \left\|u_l^{k-2}[3(\partial_xu_l)^2(\partial_xu_h)+3(\partial_xu_l)(\partial_xu_h)^2+(\partial_xu_h)^3]\right\|_{L^{2}}\\
	&+\left\|\left(\sum_{j=1}^{k-2}C_{k-2}^ju_l^{k-2-j}u_h^j\right)(\partial_xu_l+\partial_xu_h)^3\right\|_{L^{2}}\\
	\lesssim \, &  \|u_l\|_{H^{2}}^{k}\|\partial_xu_h\|_{L^{\infty}}+\|u_l\|_{H^{2}}^{k-1}\|\partial_xu_h\|_{L^{\infty}}^2+\|u_l\|_{H^{2}}^{k-2}\|\partial_xu_h\|_{L^{\infty}}^3+\sum_{j=1}^{k-2}\|u_l\|_{H^{2}}^{k+1-j}\|u_h\|_{L^{\infty}}^j\\
	&+\sum_{j=1}^{k-2}\|u_l\|_{H^{2}}^{k-j}\|u_h^j\partial_xu_h\|_{L^{\infty}}+\sum_{j=1}^{k-2}\|u_l\|_{H^{2}}^{k-1-j}\|u_h^j(\partial_xu_h)^2\|_{L^{\infty}}\\
	&+\sum_{j=1}^{k-2}\|u_l\|_{H^{2}}^{k-2-j}\|u_h^j(\partial_xu_h)^3\|_{L^{\infty}}.
\end{align*}
Repeating the  analysis in \eqref{error E3 estimate k geq 2}, one has that for $k\geq 3$ and $t\in[0,T_l]$,
\begin{equation*}
	n^{r_s}\gtrsim
	\left\{
	\begin{aligned}
		&\|u_l\|_{H^{2}}^{k}\|\partial_xu_h\|_{L^{\infty}},\  
		\|u_l\|_{H^{2}}^{k-1}\|\partial_xu_h\|_{L^{\infty}}^2,\\ 
		&\sum_{j=1}^{k-2}\|u_l\|_{H^{2}}^{k+1-j}\|u_h\|_{L^{\infty}}^j,\
		\sum_{j=1}^{k-2}\|u_l\|_{H^{2}}^{k-j}\|u_h^j\partial_xu_h\|_{L^{\infty}},\\ &\sum_{j=1}^{k-2}\|u_l\|_{H^{2}}^{k-1-j}\|u_h^j(\partial_xu_h)^2\|_{L^{\infty}},
	\end{aligned}
	\right.
\end{equation*}
and therefore it suffices to estimate the different terms, which are
\begin{align*}
	\|u_l\|_{H^{2}}^{k-2}\|\partial_xu_h\|_{L^{\infty}}^3
	\lesssim \, &  n^{(k-2)(\frac{\delta}{2}-\frac{1}{k})}n^{3(1-\frac{\delta}{2}-s)}\\
	=\,  & n^{-3s+2+\frac{2}{k}+(k-5)\frac{\delta}{2}}\lesssim n^{r_s}, \ \ k\ge 3,\ \ t\in [0,T_l],\\
	\sum_{j=1}^{k-2}\|u_l\|_{H^{2}}^{k-2-j}\|u_h^j(\partial_xu_h)^3\|_{L^{\infty}}
	\lesssim \, &  \sum_{j=1}^{k-2}n^{(k-2-j)(\frac{\delta}{2}-\frac{1}{k})}n^{j(-\frac{\delta}{2}-s)+3(1-\frac{\delta}{2}-s)}\\
	\lesssim \, &  \sum_{j=1}^{k-2}n^{j(-s-\delta+\frac1k)+2+\frac2k-3s+(k-5)\frac{\delta}{2}}\\
	\lesssim \, &   n^{-4s+\frac3k+2+(k-7)\frac{\delta}{2}}\lesssim n^{r_s}, \ \ k\ge 3, 
	\ \ t\in [0,T_l].
\end{align*}
Combining the above estimations, we get 
\begin{align}\label{error E4 estimate}
	\|E_4\|_{H^{\rho_0}}\lesssim n^{r_s},\ \ k\ge1,\ \ t\in[0,T_l].
\end{align}

\textbf{(v) Estimating $\|\EE\|_{H^{\rho_0}}$.}
Let $T_l>0$ be given in Lemma \ref{lemma ul Hr} such that $T_l$ does not depend on $n$.  Let $t\in[0,T_l]$, by virtue of the It\^{o} formula and \eqref{Error}, we derive that
\begin{align*}
	&\, \|\EE(t,x)\|_{H^{\rho_0}}^2\\
	\leq \,  & \left|-2\int_0^t(h(t',u_{m,n}) \dd \W,\EE)_{H^{\rho_0}}\right|+2\sum_{i=1}^{4}\int_0^t\left|(E_i,\EE)_{H^{\rho_0}}\right| \dd t'
	+\int_0^t\|h(t',u_{m,n})\|_{\LL_2(\U;H^{\rho_0})}^2\dd t'.
\end{align*}
Taking supremum with respect to $t\in[0,T_l]$ and then using the BDG inequality give rise to
\begin{align*}
	&\E\sup_{t\in[0,T_l]}\|\EE(t)\|_{H^{\rho_0}}^2\notag\\
	\leq  \, &    \frac 12\E\sup_{t\in[0,T_l]}\|\EE(t)\|_{H^{\rho_0}}^2
	+C\E\int_0^{T_l}\|h(t,u_{m,n})\|_{\LL_2(\U;H^{\rho_0})}^2\dd t\\
	&+C\int_0^{T_l}
	\Big[\sum_{i=1}^4\E\|E_i\|_{H^{\rho_0}}^2+\E\|\EE(t)\|_{H^{\rho_0}}^2\Big]\dd t.
\end{align*}
For any fixed $s>\frac52$, since $\|u_{m,n}\|_{H^s}\lesssim 1$, on account of \eqref{Ass2-small assumption}, Lemmas \ref{lemma uh Hr} and \ref{lemma ul Hr}, we can pick $N=N(s,k)\gg 1$ such that
\begin{align*}
	\|h(t,u_{m,n})\|_{\LL_2(\U;H^{\rho_0})}^2
	\lesssim \left({\rm e}^\frac{-1}{\|u_{m,n}\|_{H^{\rho_0}}}\right)^2 
	\lesssim  \left(n^{N(-s+\rho_0)}+n^{N(\frac{\delta}{2}-\frac{1}{k})}\right)^2 
	\lesssim  n^{2r_s}.
\end{align*}
This,  \eqref{error E1 estimate}, \eqref{error E2 estimate}, \eqref{error E3 estimate} and  \eqref{error E4 estimate} yield
\begin{align*}
	&\E\sup_{t\in[0,T_l]}\|\EE(t)\|_{H^{\rho_0}}^2\notag\\
	\leq  \, &    \frac 12\E\sup_{t\in[0,T_l]}\|\EE(t)\|_{H^{\rho_0}}^2
	+C\E\int_0^{T_l}\|h(t,u_{m,n})\|_{\LL_2(\U;H^{\rho_0})}^2\dd t\\
	&+C\int_0^{T_l}
	\Big[\sum_{i=1}^4
	\E\|E_i\|_{H^{\rho_0}}^2+\E\|\EE(t)\|_{H^{\rho_0}}^2\Big]\dd t\\
	\leq  \, &    \frac 12\E\sup_{t\in[0,T_l]}\|\EE(t)\|_{H^{\rho_0}}^2+C(T_l)n^{2r_s}+C\int_0^{T_l}\E\sup_{t'\in[0,t]}\|\EE(t')\|_{H^{\rho_0}}^2\dd t
\end{align*}
Obviously, for each $n\geq 1$, $\E\sup_{t\in[0,T_l]}\|\EE(t)\|_{H^{\rho_0}}^2$ is finite. Then by the Gr\"{o}nwall inequality, we have
\begin{align*}
	\E\sup_{t\in[0,T_l]}\|\EE(t)\|_{H^{\rho_0}}^2\leq Cn^{2r_s},\ C=C(T_l).
\end{align*}

\section*{Acknowledgement}
The authors would like to express their gratitude to the anonymous referee for the valuable suggestions, which have led to substantial improvements  in this paper. H.T would also like to record his indebtedness to Professor Giulia Di Nunno.

\section*{Declaration}
The authors declare that they have no conflict of interest and data sharing is not applicable to this article  since no datasets were generated or analyzed during the current study.

\end{document}